\documentclass[final,leqno]{siamltex704}

\usepackage{graphicx}      
\usepackage{amsmath}       
\usepackage{amssymb}       
\usepackage{amsfonts}      
\usepackage{rotating}
\usepackage{epsfig}
\usepackage{exscale,cmmib57}

\newcommand{\mypsdraft}{\psdraft}
\renewcommand{\mypsdraft}{\psfull}
\newcommand{\mypsfull}{\psfull}

\renewcommand{\ldots}{\dotsc}
\def\dst{\displaystyle}

\newtheorem{prop}[theorem]{Proposition}
\newtheorem{remark}[theorem]{\textit{Remark}}
\newtheorem{algorithm}[theorem]{Algorithm}
\newtheorem{example}[theorem]{\textit{Example}}

\newcommand{\der}[2]{\frac{\partial #1}{\partial #2}}
\newcommand{\veps}{\varepsilon}

\newcommand{\sdeg}{\operatorname{sdeg}}

\newcommand{\phm}{\phantom{-}}
\newcommand{\mbb}[1]{\mathbb{#1}}
\newcommand{\mb}[1]{\mathbf{#1}}
\newcommand{\mc}[1]{\mathcal{#1}}
\newcommand{\jd}{\displaystyle}
\newcommand{\js}{\scriptstyle}

\newcommand{\e}[1]{{(#1)}}

\newcommand{\brak}[1]{\langle #1 \rangle}

\newcommand{\wtil}{\widetilde}

\newcommand{\tem}[1]{\text{\it #1}}

\title{Practical Error Estimates for Reynolds' Lubrication
Approximation and its Higher Order Corrections\thanks{
This work was supported in part
    by the Director, Office of Science, Advanced Scientific Computing
    Research, U.S. Department of Energy contract 
    DE-AC02-05CH11231.
}}

\author{Jon Wilkening\thanks{Department of Mathematics and Lawrence Berkeley National
    Laboratory, University of California, Berkeley, CA  94720
    (wilken@math.berkeley.edu).}}

\begin{document}
\maketitle



\begin{abstract}
Reynolds' lubrication approximation is used extensively to study flows
between moving machine parts, in narrow channels, and in thin films.
The solution of Reynolds' equation may be thought of as the zeroth
order term in an expansion of the solution of the Stokes equations in
powers of the aspect ratio $\varepsilon$ of the domain.  In this paper, we
show how to compute the terms in this expansion to arbitrary order on
a two-dimensional, $x$-periodic domain and derive rigorous, a priori
error bounds for the difference between the exact solution and the
truncated expansion solution.  Unlike previous studies of this sort,
the constants in our error bounds either are independent of the
function $h(x)$ describing the geometry or depend on $h$ and its
derivatives in an explicit, intuitive way.  Specifically, if the
expansion is truncated at order $2k$, the error is
$O(\varepsilon^{2k+2})$, and $h$ enters into the error bound only through its
first and third inverse moments $\int_{0}^{1}h(x)^{-m}\,dx$, $m=1,3$, and
via the max norms $\|\frac{1}{\ell!}h^{\ell-1}\partial_{x}^{\ell}h\|_{\infty}$, $1\le\ell\le2k+2$.  We validate our estimates by
comparing with finite element solutions and present numerical evidence
that suggests that even when $h$ is real analytic and periodic, the
expansion solution forms an asymptotic series rather than a convergent
series.
\end{abstract}

\begin{keywords}
  incompressible flow, lubrication theory, asymptotic expansion,
  Stokes equations, thin domain, a priori error estimates
\end{keywords}

\begin{AMS}
  76D08, 35C20, 41A80
\end{AMS}

\begin{DOI}
10.1137/070695447
\end{DOI}

\pagestyle{myheadings}
\thispagestyle{plain}
\markboth{JON WILKENING}{ERROR ESTIMATES FOR REYNOLDS' APPROXIMATION}

\section{Introduction}

Reynolds' lubrication equation
\cite{Rey:86,poz:intro,langlois,elrod:60} is used extensively in
engineering applications to study flows between moving machine parts,
e.g.,  in journal bearings or computer disk drives.  It is also used in
microfluid and bio-fluid mechanics to model creeping flows through narrow
channels and in thin films.  Although there is a vast literature
(including several textbooks) on viscous flows in thin geometries, the
equations are normally derived either directly from physical arguments
\cite{langlois} or using formal asymptotic arguments \cite{elrod:60}.
This is acceptable in most circumstances as the original equations
(Stokes or Navier--Stokes) have also been derived from physical
considerations, and by now the lubrication equations have been used
frequently enough that one can draw on experience and intuition to
determine whether they will work well for a given problem.

On the other hand, as soon as the geometry of interest develops (or
approaches) a singularity, or if we wish to compute several terms in
the asymptotic expansion of the solution in powers of the aspect ratio
$\veps$, we rapidly leave the space of problems for which we can use
experience as a guide; thus, it would be helpful to have a rigorous
proof of convergence to serve as a guide to identify the features of
the geometry that could potentially invalidate the approximation.  For
example, in \cite{snail}, Wilkening and Hosoi used lubrication
theory to study the optimal wave shapes that an animal such as a
gastropod should use as it propagates ripples along its muscular foot
to crawl over a thin layer of viscous fluid.  In certain limits of
this constrained optimization problem, the optimal wave shape develops
a kink or cusp in the vicinity of the region closest to the substrate,
and there is a competing mechanism controlling the size of the
modeling error (singularity formation versus nearness to the substrate).
We found that shape optimization within (zeroth order) lubrication
theory drives the geometry out of the realm of applicability of
the lubrication model; however,
by computing higher order corrections and monitoring the errors (using
the results of this paper), we learned that cusp-like singularities
are appropriately penalized by the full Stokes equations, yielding
nonsingular optimal solutions; see \cite{snail} for further details.

\subsection{Previous work}

In most of the following papers, the Stokes or Navier--Stokes equations
are solved in a domain $\Omega_\veps$ bounded below by a flat
substrate and above by a curved boundary $y=\veps h(x)$ in two
dimensions, or $z=\veps h(x,y)$ in three dimensions, where $\veps$ is
a small parameter and the function $h$ is fixed.  These solutions are
then compared to the solution of Reynolds' equation (or to a truncated
expansion solution of the Stokes or Navier--Stokes equations), and the
error is shown to converge to zero in the limit as
$\veps\rightarrow0$.

In 1983, Cimatti \cite{cimatti} used a stream function formulation
to compare the solution of Reynolds' equation to that of the Stokes
equation in two dimensions.  The key idea of the proof, which all
subsequent studies (including this one) also use, is that the
Poincar\'e--Friedrichs inequality holds uniformly as
$\veps\rightarrow0$ for the rescaled biharmonic equation (where the
domain $\Omega=\Omega_{\veps=1}$ is held fixed and the equations
contain the small parameter).
Cimatti assumes $h$ has four weak derivatives (whereas, we require only 
$h\in C^{1,1}$) and shows that for any compact set $K\subset\Omega$,
\begin{equation}
  \label{eqn:cimatti1}
  \|\veps u-\bar{u}\|_{L^2(\Omega)}\le C\veps, \qquad
  \max\left(\|\veps^3p_x - \bar{p}_x\|_{L^2(K)},
  \|\veps^2p_y\|_{L^2(K)}\right)
  \le C\veps^{1/2},
\end{equation}
where $u$ is the $x$-component of velocity, $p$ is the pressure, a bar
denotes the solution of Reynolds' equation, and $C$ is independent of
$\veps$ but depends on $h$ in the first inequality and on $h$ and $K$
in the second.  The scaling here in not standard: he imposes the
boundary condition $\veps u(x,0)=\bar{u}(x,0)=\const$, which accounts
for the extra factor of $\veps$ in each of the left-hand sides of
(\ref{eqn:cimatti1}).  There are a few problems with Cimatti's
analysis, notably the dependence of $C$ on $L$ (the ``arbitrary
cutoff'' used to make the unbounded domain bounded) and the fact that
some of his arguments seem to require $\veps$ to be small in
comparison to $C^{-1}$; however, his basic approach is interesting and
inspired much of the work that followed in this subject.

In 1986, Bayada and Chambat \cite{bayada} generalized Cimatti's
work to three dimensions. They analyze the Stokes equations directly
rather than using a stream function formulation, assume less
regularity of $h$ (apparently only $h\in C^1$), and state their
results in terms of limits (i.e., the quantities $u^\veps_i$, $\veps
\partial_xu^\veps_i$, $\partial_yu^\veps_i$, and $p^\veps$ in the
solution of the Stokes equations converge in $L^2$ to the
corresponding quantities in the solution of Reynolds equations as
$\veps\rightarrow0$); hence, they do not give rates of convergence.
In a later paper \cite{bayada2}, they also studied the asymptotics of
the solution at a junction between a three-dimensional Stokes flow and
a thin film flow.

In 1990, Nazarov \cite{nazarov} generalized previous work to 
the case of the Navier--Stokes equations and also showed how
to treat higher order corrections in an asymptotic expansion in
the small parameter $\veps$.  He proved that if $h(x,y)$ is
smooth, then there is a constant $C$ depending on $h$, $N$, and
the boundary conditions such that
\begin{equation}
\left\|\mb{u}-\mb{u}^N\right\|_{H^1} + \left\|p-\veps^{-1}p^N\right\|_{L^2}\le C\veps^{N-1/2},
\end{equation}
where $(\mb{u},p)$ is the solution of the Navier--Stokes equations,
$\mb{u}^N$ and $p^N$ are the terms of the asymptotic expansion
truncated at the $N$th order (including a boundary layer expansion
near the lateral edges of the thin domain), and the norms are taken on
the thin domain $\Omega_\veps$ (rather than the rescaled domain
$\Omega$).  As a corollary, if the expansion is computed with
``superfluous'' terms that are afterwards treated as remainders, he
obtains the optimal estimate
\begin{equation}
  \label{eqn:nazarov}
  \left\|\mb{u}-\mb{u}^N\right\|_{L^2} +
  \veps^{1/2}\left\|\left(\pi_{1/2}^\veps\nabla\right)\left(\mb{u}-\mb{u}^N\right)\right\|_{L^2}
  +\left\|p-\veps^{-1}p^N\right\|_{L^2}\le C\veps^{N+1}.
\end{equation}
Nazarov's paper is concise to the point of being impenetrable at
times.  We interpret $\pi_{1/2}^\veps\nabla=(\partial_x,
\partial_y,\veps^{1/2}\partial_z)$, but this symbol was not defined
and may actually be a variable coefficient operator that incorporates
the boundary conditions in its definition.  We are also unsure of the
definition of $p$ and $p^N$, as we would have expected $p-\veps^{-2}p^N$
to appear together.\enlargethispage{9pt}

In a later paper~\cite{nazarov:point}, Nazarov studies the asymptotics
of the solution of the Stokes equations in a domain in which two
smooth surfaces meet at a point.  This problem is also studied in a
recent paper of Ciuperca, Hafidi, and Jai~in \cite{ciuperca:06}.  This singular
limit is interesting in that deriving even the first correction to
the zeroth order approximation in the asymptotic expansion remains an
open problem.

Assemien, Bayada, and Chambat \cite{assemien:94}
have studied the important question of the effect of inertia
on the asymptotic behavior of a thin film flow, which can in many
cases be significant, requiring that the Navier--Stokes equations be
used in place of the Stokes equations as the underlying model for
the asymptotic expansion.
We also mention that there is a large body of literature on the
long-time behavior of solutions of the Navier--Stokes equations on thin
domains; see, e.g., \cite{raugel:sell,temam:thin:NS}.

In 2000, Duvnjak and Maru\u{s}i\'c-Paloka \cite{paloka}
showed how to rigorously analyze the lubrication approximation of the
Navier--Stokes equations for a slipper bearing in a circular
geometry.  The focus of their paper is on formulating the problem
in cylindrical coordinates and showing how to adapt the zeroth
order case of Nazarov's proof to handle the change of variables.
Elrod's pioneering 1960 paper \cite{elrod:60} is also concerned
with the (formal) relationship between the Navier--Stokes equations and
Reynolds' equation for this geometry.

\subsection{Motivation and summary}

None of the studies described above shows how the constant $C$
bounding the error depends on the function $h(x)$ describing the
geometry.  This is because most theorems of analysis
give constants that depend on the domain~$\Omega$, which is usually
fixed.  But in our case, the data $h(x)$ of the problem actually
specifies the domain; therefore, to obtain bounds that are independent
of $h$, one must avoid or modify standard arguments for flattening the
boundary, etc.,~so as not to lose track of $h(x)$ in the analysis.
Moreover, arguments based on the closed graph theorem or Rellich's
compactness theorem must be avoided entirely, as these also depend on
the geometry.  This forces us to look for new ways to analyze old
problems using tools that furnish explicit constants.

In this paper, we consider only the two-dimensional, periodic Stokes
equations with a specific choice of boundary conditions, but we derive
error estimates that depend on $h$ in an explicit, intuitive way.  Our
main result is summarized in Theorem~\ref{thm:bound2}, which may be
stated as follows: Let $T=[0,1]_p$ be the periodic unit interval.  If
$k\ge0$, $h\in C^{2k+1,1}(T)$, $0<h_0\le h(x)\le1$ for $x\in T$, and
$\veps\le r_0/3$ (defined below), then the error in truncating the
expansion of the stream function, velocity, vorticity, and pressure (in
appropriate $\veps$-weighted Sobolev norms) at order $2k$ (keeping in
mind that only even powers of $\veps$ appear in these expansions) is
bounded by
\begin{equation}
  \label{eqn:star:def:intro}
    \sqrt{I_1}\left( |V_0| + |V_1|\right)
  \left[1 + \theta_{k} \frac{\veps}{r_{k}}\sqrt{\frac{I_3}{I_1}}
  \right] \left( \frac{\veps}{\rho_{k} r_{k}} \right)^{2k+2},
\end{equation}
where $V_0$ and $V_1$ are prescribed tangential velocities on the
lower and upper boundaries of the domain
\begin{equation}
  \label{eqn:rk:def:intro}
  r_k = \left(\max_{1\le \ell\le 2k+2}
  \left\{\left\|\frac{1}{\ell!}h^{\ell-1}
  \partial_x^\ell h\right\|_\infty^{1/\ell}\right\}\right)^{-1}, \qquad
  I_m = \int_0^1 h(x)^{-m}\,dx
\end{equation}
and $\rho_k$, $\theta_k$ are constants independent of $h$.
The bound on pressure has another term
involving $h_0$; see~(\ref{eqn:thm:bound2}) below.

The constants in (\ref{eqn:star:def:intro}) have been divided into two
types: those that are (1) given in the problem statement or easily
computable from $h$; or (2) difficult to compute but universal
(independent of $h$).
We show how to compute the constants in the latter category
($\rho_k$ and $\theta_k$) in section~\ref{sec:error}; see
Table~\ref{tbl:rho:theta}.
The constants in the former category ($r_k$ and $I_m$) help us
understand the competing mechanism of singularity formation
versus~proximity to the substrate: the curvature and higher derivatives
are allowed to diverge as long as the gap size simultaneously
approaches zero in such a way that the homogeneous products
$\frac{1}{\ell!}h^{\ell-1}\partial_x^\ell h$ remain uniformly bounded.
Although the factors $\sqrt{I_1}$ and $\sqrt{I_3/I_1}$ in
(\ref{eqn:star:def:intro}) also diverge in this limit, the norm of the
exact solution diverges at a similar rate --- so the relative error 
in the expansion solution truncated at order $2k$ is
$O(\veps^{2k+2})$, with $\rho_kr_k$ serving as an effective radius
of convergence.

The framework we have chosen for this
paper is intended to be general enough to cover many interesting
applications (such as a crawling gastropod~\cite{snail} or an
``unwrapped'' slipper bearing) but simple enough to obtain explicit
detailed estimates that reveal the dependence of the error on the
geometry $h(x)$.  We also wanted to determine whether there might
exist geometries for which the asymptotic expansion yields a
convergent series.  Although we do not have a rigorous proof, the
answer appears to be negative even for the simplest case of a real
analytic function such as $h(x)=\frac{3}{5}+\frac{2}{5}\sin 2\pi x$, 
for which the $r_k$ in (\ref{eqn:rk:def:intro}) are bounded away from
zero.  It is hoped that this work will serve as a useful
first step toward obtaining similar error estimates for
three-dimensional problems that include more general boundary
conditions, incorporate end effects near the lateral edges of the
domain (which we avoid by studying the periodic case), and include the
effect of inertia or viscoelasticity.

\subsection{Outline}

In section~\ref{sec:derive}, we derive Reynolds' lubrication
approximation in its primitive and stream function formulations.
In section~\ref{sec:expand}, we show how to compute successive terms
in an asymptotic expansion of the stream function.  In
section~\ref{sec:struc},
we prove a structure theorem describing the dependence of these
terms on $h(x)$ and its derivatives.

\looseness=1In section~\ref{sec:error}, we formulate the problem weakly and
analyze the truncation error equation using weighted Sobolev spaces
and a uniform Poincar\'e--Friedrichs argument.  The first challenge is
to find the right weighted norms on the lower and upper boundaries
(equivalent to $H^{1/2}(\Gamma_0)$ and $H^{1/2}(\Gamma_1)$ for fixed
$\veps$) to yield manageable error estimates in terms of $h$ when
we change variables to straighten out the boundaries.
In section~\ref{sec:error:est}, we reduce the problem of bounding the
truncation errors to that of bounding the second and fourth
derivatives of the two highest order terms retained in the asymptotic
expansion, namely, $\|\psi_{xx}^\e{2k}\|_0$ and
$\|h^2\psi_{xxxx}^\e{2k-2}\|_0$.  We then use the structure
theorem of section~\ref{sec:struc} to compute these norms in order to
obtain the constants $\rho_k$ and $\theta_k$ in
(\ref{eqn:star:def:intro}) for $0\le k\le 25$.  In
section~\ref{sec:pressure}, we show how to compute the error in
velocity, vorticity, and pressure from that of the stream function.
This requires that we determine how the Babu\u{s}ka--Brezzi inf-sup
constant $\beta$ depends on $h(x)$; see \cite{infsup}.

In section~\ref{sec:fe}, we validate our results by comparing to
``exact'' solutions (computed using finite elements) for a geometry
typical of engineering applications.  The result of
this comparison is that the effective radius of convergence
$r_k\rho_k$ is within a factor of 3 of optimal for $k=5$, $k=10$,
and perhaps all $k\ge5$.  These calculations
also suggest that even when $h(x)$ is real analytic, the
expansion solution is an asymptotic series rather than a convergent
series.  This is because the constants $\rho_k$ converge to
zero as $k\rightarrow\infty$.  Fortunately, $\rho_k$ initially
increases and does not become smaller than $\rho_0=0.197$
until $2k=26$, which is already outside of the practical range
of~$k$.
Finally, in Appendix~\ref{sec:impl}, we present our numerical
algorithm for computing the expansion solutions, which can be
performed symbolically using a computer algebra system such as
Mathematica or in floating point arithmetic, e.g., in
$C^{++}$.

\section{Reynolds' approximation}
\label{sec:derive}

\begin{figure}[b!]
\centerline{\includegraphics[width=.46\linewidth]{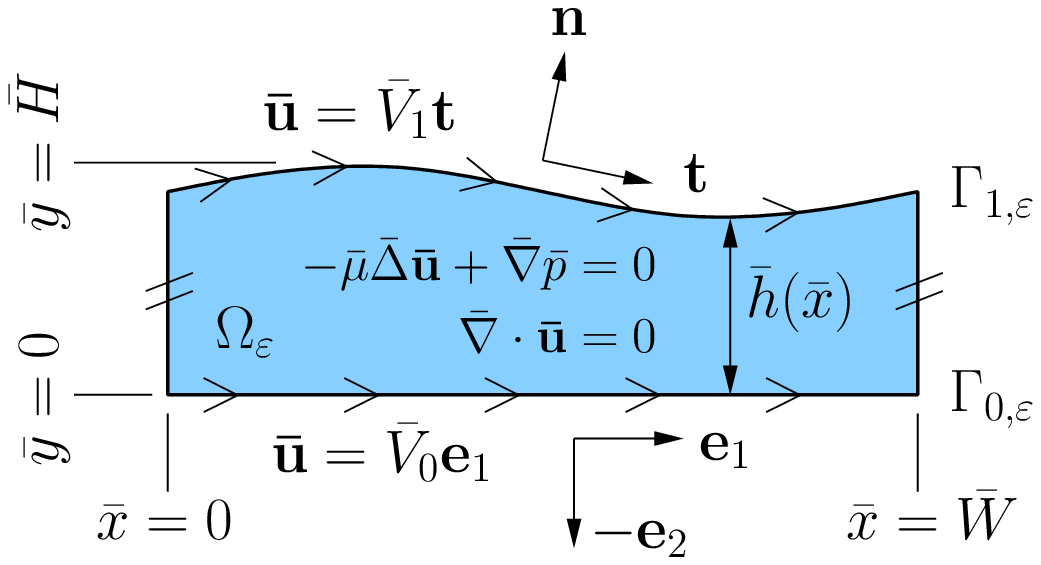} \qquad
  \includegraphics[width=.46\linewidth]{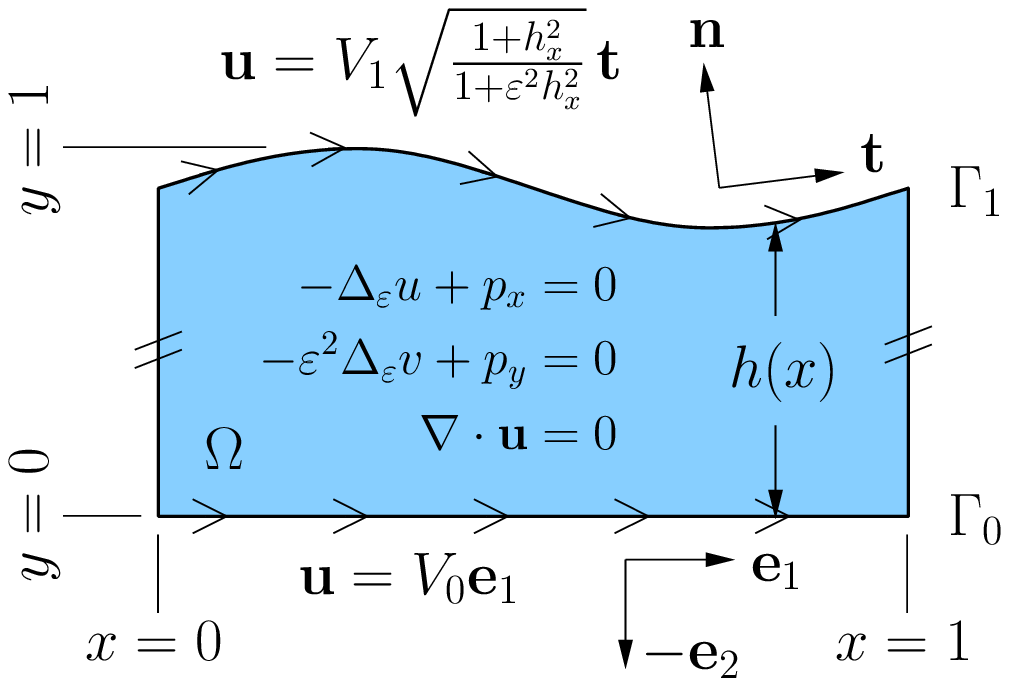}}
    \caption{ Geometry commonly encountered in
      lubrication-type problems.  Left: Physical coordinate
      system.  Right: Dimensionless coordinate system
      ($\Delta_\veps = \veps^2\partial_x^2+\partial_y^2$).
}
    \label{fig:geom}
\end{figure}

Consider the Stokes equations on a periodic domain of width $\bar{W}$
bounded below by a flat wall moving with constant speed $\bar{V}_0$
and above by an inextensible sheet moving with constant speed
$\bar{V}_1$ along a fixed curve
$\Gamma_{1,\veps}=\{(\bar{x},\bar{h}(\bar{x})):0\le \bar{x}\le
\bar{W}\}$; see Figure~\ref{fig:geom}.  A bar is used to distinguish 
a physical variable from its dimensionless counterpart.  We
nondimensionalize the variables by choosing a characteristic speed
$\bar{U}$ and height $\bar{H}$ for the problem, and set
$\bar{x}=\bar{W}x$, $\bar{y}=\bar{H}y$,
$\bar{h}(\bar{x})=\bar{H}h(x)$, $\bar{V}_i=\bar{U}V_i$,
$(\bar{u},\bar{v})=\bar{\mb{u}}=(\bar{U}u,\bar{U}\frac{\bar{H}}{\bar{W}}v)$,\vspace*{1pt} 
and $\bar{p} = \bar{\mu}\frac{\bar{U}\bar{W}}{\bar{H}^2}p$.  The 
stream function $\psi$, flux $Q$, and vorticity $\omega$ introduced
below satisfy\vspace{1.5pt} $\bar{\psi}=\bar{U}\bar{H}\psi$, $\bar{Q}=\bar{U}\bar{H}Q$, 
and $\bar{\omega}=\frac{\bar{U}}{\bar{H}}\omega$.

%

We have in mind a situation where the aspect ratio
$\veps=\bar{H}/\bar{W}$ of the physical domain is small.  By scaling
the $x$- and $y$-axes differently, we map the problem onto a nicer
geometry, which introduces terms in the equations that vanish in the
singular limit $\veps\rightarrow0$.  Specifically, we wish to find
$x$-periodic functions $u,v,p$ defined on the rescaled domain
\begin{equation}
  \Omega=\{(x,y)\;:\;0\le x\le 1,\;\;0<y<h(x) \}
\end{equation}
such that
\begin{equation}
  \label{eqn:naive}
  p_x=\veps^2u_{xx}+u_{yy}, \qquad
  p_y=\veps^4v_{xx}+\veps^2 v_{yy}, \qquad
  v_y = -u_x \qquad (\text{in } \Omega)
\end{equation}
subject to periodic boundary conditions on the left and right
sides of $\Omega$ and
\begin{equation}
  \label{eqn:naive:bc}
  (u,v)\left\vert_{\Gamma_0} = (g_0,0), \qquad
  (u,v)\right\vert_{\Gamma_1} = (g_1,h_x g_1)
\end{equation}
on the bottom and top boundaries.  Here
\begin{equation}
  \label{eqn:gdef}
  g_0(x)=V_0, \qquad g_1(x)=V_1\left[1+\veps^2h'(x)^2\right]^{-1/2}, 
\end{equation}
i.e., $g_1(x)=V_1\cos\theta(x)$, where\, $\theta=\arctan(\veps h_x)$\, is
the angle of the curve $\bar{h}(\bar{x})$ relative to the horizontal. 
Reynolds' lubrication approximation is
obtained by setting $\veps=0$ in the equations and solving
\begin{equation}
  \label{eqn:rle}
  p_x=u_{yy}, \quad\;\;
  p_y=0, \quad\;\;
  v_y=-u_x, \quad\;\;
  \mb{u}\left\vert_{\Gamma_0} = (V_0;0), \quad\;\;
  \mb{u}\right\vert_{\Gamma_1} = (1;h_x)V_1.
\end{equation}
If we write (\ref{eqn:naive}) in the form
$L(\mb{u};p)=(0;0;0)$, where $L=L^\e{0} + \veps^2L^\e{2} + \veps^4L^\e{4}$
is given by
\begin{equation}
  L =
  \begin{pmatrix}
    -\partial_y^2 & 0 & \partial_x \\
    \phm0 & 0 & \partial_y \\
    \phm\partial_x & \partial_y & 0
  \end{pmatrix} + \veps^2
  \begin{pmatrix}
    -\partial_x^2 & \phm0 & 0 \\
    \phm0 & -\partial_y^2 & 0 \\
    \phm0 & \phm0 & 0
  \end{pmatrix} + \veps^4
  \begin{pmatrix}
    0 & \phm0 & 0 \\
    0 & -\partial_x^2 & 0 \\
    0 & \phm0 & 0
  \end{pmatrix},
\end{equation}
then (\ref{eqn:rle}) is just the zeroth order system
$L^\e{0}(\mb{u};p)=(0;0;0)$ with zeroth order boundary conditions
(expanding $g_0$ and $g_1$ in (\ref{eqn:naive:bc}) in powers of
$\veps$).  The equation for $v$ decouples from the others, and we find
that $p$ is independent of $y$ and
\begin{equation}
  \label{eqn:u}
  u(x,y) = \left(\frac{y^2}{2}-\frac{h(x)y}{2}\right)p_x(x) +
  \left(1-\frac{y}{h(x)}\right)V_0 + \frac{y}{h(x)}V_1.
\end{equation}
Integrating from $0$ to $h$ and solving for $p_x$, we obtain
\begin{equation}
\label{eqn:px}
  p_x = \frac{6}{h^2}(V_0+V_1) - \frac{12}{h^3}Q,
\end{equation}
where $Q=\int_0^h u(x,y)\,dy$ is the volume flux through any cross
section of the fluid.  ($Q$~is constant since $\nabla\cdot\mb{u}=0$
and $\mb{u}$ is tangent to $\Gamma_0$ and $\Gamma_1$).
Since $p$ is periodic, $\int p_x\,dx=0$, and
we find that
\begin{equation}
  \label{eqn:Q}
  Q = \frac{V_0+V_1}{2}\frac{I_2}{I_3}	 \qquad\qquad
  \left(I_m = \int_0^1 h(x)^{-m}\,dx\right).
\end{equation}
Substituting (\ref{eqn:Q}) and (\ref{eqn:px}) into (\ref{eqn:u}) and
using $v_y=-u_x$, $v(x,0)=0$, we obtain the solution
\begin{align}
  \notag
  p_x &= \frac{6(V_0+V_1)}{h^2}\left(1-\frac{I_2}{I_3h}\right), \\
  \label{eqn:lubr1}
  u &= (V_0+V_1)\left(3\frac{I_2}{I_3h}-3\right)\left(
  \frac{y}{h} - \frac{y^2}{h^2}\right) + \left(1-\frac{y}{h}\right)V_0
  + \frac{y}{h}V_1, \\
  \notag
  v &= (V_0+V_1)\left(3\frac{I_2}{I_3h}-2\right)\left(
  \frac{y^2}{h^2} - \frac{y^3}{h^3}\right)h_x + V_1\frac{y^2}{h^2}h_x.
\end{align}
The vertical component $v$ of the velocity field is customarily omitted
from zeroth order lubrication theory as $\bar{v}=\veps\bar{U}v$ is
$O(\veps)$ on the thin geometry $\Omega_\veps$ of
Figure~\ref{fig:geom}.

We may also derive (\ref{eqn:lubr1}) using a stream function
formulation of the problem.  Our procedure for computing higher order
corrections to the lubrication approximation and our method for
estimating the error of these expansion solutions are both done in the
stream function formulation.  Let us define
\begin{equation}
  \Delta_\veps = \veps^2\partial_x^2 + \partial_y^2, \qquad\qquad
  \mb{u}_\veps = \begin{pmatrix} u \\ \veps^2 v \end{pmatrix}.
\end{equation}
In
our error estimates below, we will need to consider the inhomogeneous
problem $L(\mb{u};p)=(F_1;F_2;0)$ with boundary conditions
(\ref{eqn:naive:bc}), i.e.,
\begin{equation}
  \label{eqn:up:bvp}
\begin{aligned}
  -\Delta_\veps \mb{u}_\veps + \nabla p &= \mb{F}, \\
  \nabla\cdot\mb{u} &= 0, 
\end{aligned}
\qquad
  \mb{u}\left\vert_{\Gamma_0} = (g_0;0), \qquad
  \mb{u}\right\vert_{\Gamma_1} = (g_1;h_xg_1).
\end{equation}
Since $\mb{u}$ is incompressible, there is a stream function $\psi$
such that
\begin{equation}
  \label{eqn:u:from:psi}
  \mb{u}=\nabla\times\psi=(\psi_y,-\psi_x), \qquad
  \nabla\times\mb{u}_\veps = \veps^2v_x-u_y = -\Delta_\veps\psi.
\end{equation}
It follows from (\ref{eqn:up:bvp}) that $\psi$ satisfies the rescaled
biharmonic equation
\begin{equation}
  \label{eqn:psi:bvp}
  \Delta_\veps^2\psi =
  \psi_{yyyy} + 2\veps^2\psi_{xxyy} + \veps^4\psi_{xxxx} =
  \nabla\times\mb{F},
\end{equation}
with periodic boundary conditions in the $x$-direction and
\begin{equation}
  \label{eqn:psi:bcs}
  \left\{
  \begin{aligned}
    \psi &= 0 \\ \psi_y &= g_0
    \end{aligned}
  \right\} \text{ on } \Gamma_0, \qquad
  \left\{
  \begin{aligned}
    \psi&=Q \\ \psi_y(x,h(x))&=g_1
    \end{aligned}
  \right\} \text{ on } \Gamma_1,
\end{equation}
where $Q=\int_0^{h(0)} u(0,y)\,dy$.  Since $p$ is periodic,
$\int_0^1 p_x(x,0)\,dx=0$, i.e.,
\begin{equation}
  \label{eqn:int:psi:yyy}
  \int_0^1 \psi_{yyy}(x,0)+F_1(x,0)\,dx=0.
\end{equation}
Conversely, suppose we are able to find a flux $Q$ and a classical
solution $\psi$ of (\ref{eqn:psi:bvp}) and (\ref{eqn:psi:bcs}) such
that (\ref{eqn:int:psi:yyy}) holds.  Then we define
$\mb{u}=\nabla\times\psi$ and note that (\ref{eqn:psi:bvp}) implies
$\nabla\times(\Delta_\veps\mb{u}_\veps+\mb{F})\equiv0$, i.e., the
integral
\begin{equation}
  \label{eqn:path:indep}
  p(x,y)=\int_\gamma (\Delta_\veps\mb{u}_\veps+\mb{F})\cdot\mb{t}
  \,ds, \qquad
    \left(\parbox{2in}{\begin{center}
	$\mb{t}$ $=$ unit tangent vector along path $\gamma$ joining
    $(0,0)$ to $(x,y)$\end{center}}
    \right)
\end{equation}
is independent of the path $\gamma$.  A canonical choice for $\gamma$ is
\begin{equation}\label{eqn:path:indep2}
  p(x,y) = \int_0^x \left[\veps^2u_{xx}+u_{yy}+F_1\right](\xi,0)\,d\xi +
  \int_0^y \left[\veps^4v_{xx}+\veps^2v_{yy}+F_2\right](x,\eta)\,d\eta.
\end{equation}
Condition (\ref{eqn:int:psi:yyy}) is equivalent to requiring
$p(1,0)=p(0,0)$, from which it follows that $p(1,y)=p(0,y)$ for $0\le
y\le h(0)$, since the integrand of the second integral in
(\ref{eqn:path:indep2}) is periodic in $x$.  By construction, the
variables $\mb{u}$, $p$ satisfy (\ref{eqn:up:bvp}), where the boundary
condition on $\Gamma_1$ follows from the fact that
$\psi_x+h_x\psi_y=0$ there; hence, classical solutions of the rescaled
biharmonic equation yield classical solutions of the rescaled Stokes
equations and vice versa.  Reynolds' approximation (\ref{eqn:lubr1})
is recovered if $\mb{F}$ and $\veps$ are set to zero in
(\ref{eqn:psi:bvp})--(\ref{eqn:int:psi:yyy}) when solving for $\psi$
and $Q$; see section~\ref{sec:recur}.

\section{Higher order corrections}
\label{sec:expand}

In this section we show how to compute successive terms in the formal
expansion of the solution of the rescaled biharmonic equation
(\ref{eqn:psi:bvp}) in powers of $\veps=\bar{H}/\bar{W}$.  For this
purpose, it is convenient to manipulate the equations assuming they
are satisfied classically.  Once we obtain formulas for the higher
order approximations, we will show (in section~\ref{sec:error}) that
they satisfy a weak formulation of the problem that makes it possible
to obtain error estimates.
See \cite{kevorkian} for background on perturbation methods in partial
differential equations.

\subsection{A recursive algorithm}
\label{sec:recur}

Matching like powers of $\veps$ in the expansion
\begin{equation}
  \left[\partial_y^4 + 2\veps^2 \partial_x^2\partial_y^2 +
  \veps^4\partial_x^4\right]\left[\psi^\e0+\veps^2\psi^\e2+
  \veps^4\psi^\e4+\cdots\right]=0,
\end{equation}
we obtain the recursion
\begin{align}
  \notag
  \psi^\e0_{yyyy} &= 0, \\
  \notag
  \psi^\e2_{yyyy} &= -2\psi^\e0_{xxyy}, \\
  \label{eqn:M:recur}
  \psi^\e{2k}_{yyyy} &= -2\psi^\e{2k-2}_{xxyy} - \psi^\e{2k-4}_{xxxx},
  & &k=2,3,4,\ldots.
  \intertext{\noindent
    The boundary conditions (\ref{eqn:psi:bcs}) become
  }
  \label{eqn:B:recur}
  B\psi^\e{2k} &= \left(0,g_0^\e{2k},Q^\e{2k},g_1^\e{2k}\right), &
  &k=0,1,2,3,\ldots,
\end{align}
where $B\psi=(\psi\vert_{\Gamma_0},\psi_y\vert_{\Gamma_0},
\psi\vert_{\Gamma_1},\psi_y\vert_{\Gamma_1})$ and $g_0(x)$,
$g_1(x)$ were defined in (\ref{eqn:gdef}):
\begin{equation}
  \label{eqn:g:recur}
  g_0^\e{2k}(x)=\begin{cases} V_0, & k=0, \\ 0, & k>0, \end{cases} \qquad
  g_1^\e{2k}(x) = V_1{-1/2\choose k}h'(x)^{2k}.
\end{equation}
Condition (\ref{eqn:int:psi:yyy}) (with $F_1=0$) becomes
\begin{equation}
  \label{eqn:int:psi2k:yyy}
  \int_0^1\psi^\e{2k}_{yyy}(x,0)\,dx=0, 
  \qquad k=0,1,2,\ldots.
\end{equation}
If $\mb{F}$ were nonzero in (\ref{eqn:psi:bvp}) and depended on
$\veps$ in such a way that $\nabla\times\mb{F}$ had an expansion
in even powers of $\veps$, we could incorporate these terms into
(\ref{eqn:M:recur}) and (\ref{eqn:int:psi2k:yyy}) as well; however, we
will assume $\mb{F}=\mb{0}$ except in section~\ref{sec:error}, where
we consider the general case only to derive error estimates for the
$\mb{F}=\mb{0}$ case.  Let us denote the right-hand side of
(\ref{eqn:M:recur}) by $f^\e{2k}(x,y)$ for $k\ge0$.  The terms
$\psi^\e{2k}$, $Q^\e{2k}$ in (\ref{eqn:M:recur}) and
(\ref{eqn:B:recur}) may be computed via
\begin{equation}
  \label{eqn:psi:Q:recur}
  \left(\psi^\e{2k}, Q^\e{2k}\right) =
  G\left(f^\e{2k},\; g_0^\e{2k}, \; g_1^\e{2k}\right), \qquad
  k=0,1,2,\ldots,
\end{equation}
where $G$ is defined by Algorithm~\ref{alg:G:def} in
Figure~\ref{fig:G:def}.
%
\begin{figure}[t]
\begin{center}
\fbox{\parbox{4.9in}{
\vspace*{4pt}
\begin{algorithm} $(\psi,Q)=G(f,g_0,g_1)${\rm :} \vspace*{2pt}
\label{alg:G:def}
\begin{tabbing}
\hspace*{.25in} \= \hspace*{.2in} \=
\hspace*{.25in} \=\kill
\> $\psi_0 = \mbb{V}^4f$ \qquad $\left(\mbb{V}= \mbox{Volterra operator: }
\jd \mbb{V}f(x,y)=\int_0^y f(x,\eta)\,d\eta\right)$ \\[5pt]
\> $\jd Q=\frac{1}{2I_3}\int_0^1 \frac{2\psi_0(x,h(x))}{h(x)^3}
  +\frac{-\psi_{0,y}(x,h(x))+g_0+g_1(x)}{h(x)^2}\,dx$ \\[7pt]
\> $\psi(x,y) = \psi_0(x,y) + \left(g_0 h(x)\right)\frac{y}{h(x)}$ \\[5pt]
\>\> $+\,\left(3Q - 3\psi_0(x,h(x))
  + \psi_{0,y}(x,h(x))h(x) - 2g_0 h(x) - g_1(x)h(x)\right)
  \frac{y^2}{h(x)^2}$ \\[5pt]
\>\> $+\,\left(-2Q + 2\psi_0(x,h(x))
  - \psi_{0,y}(x,h(x))h(x) + g_0 h(x) + g_1(x)h(x)\right)
  \frac{y^3}{h(x)^3}$ \\[5pt]
\> \textbf{return} $(\psi,Q)$
\end{tabbing}
\vspace*{3pt}
\end{algorithm}}}
\end{center}
\caption{Algorithm to solve $\psi_{yyyy}=f$, $B\psi=(0,g_0,Q,g_1)$,
$\int_0^1\psi_{yyy}(x,0)\,dx=0$.}
\label{fig:G:def}
\end{figure}
In this algorithm, we solve \mbox{$\psi_{yyyy}=f$} by integrating four
times in the $y$-direction and then correct the boundary conditions
with a cubic polynomial.  The formula for $Q$ in the algorithm may be
derived from the one for $\psi$ as follows.  As $\psi_{0,yyy}(x,0)=0$,
the requirement that $\int_0^1\psi_{yyy}(x,0)\,dx=0$ is equivalent to
the condition
\begin{equation}
  0=6\int_0^1 \frac{-2Q+2\psi_0 - \psi_{0,y}h + g_0h + g_1h}{h^3}\,dx.
\end{equation}
Solving for $Q$ and using $\int h^{-3}\,dx=I_3$ gives the result.

The formulas $(u,v)=(\psi_y,-\psi_x)$,\, $\omega=\veps^2
v_x-u_y$,\, $p_x=u_{yy}+\veps^2u_{xx}$, and $p_y=\veps^2
v_{yy}+\veps^4v_{xx}$ allow us to compute the expansions of $\mb{u}$,
$\omega$, and $p$ in terms of $\psi$:
\begin{equation}
  \label{eqn:uvp:from:psi}
  \begin{aligned}
    u^\e{2k} &= \psi_y^\e{2k}, &
    v^\e{2k} &= -\psi_x^\e{2k}, 
    & \quad &k\ge0, \\[2pt]
    \omega^\e0 &= -\psi_{yy}^\e0, &
    \omega^\e{2k} &= -\psi_{xx}^\e{2k-2}-\psi_{yy}^\e{2k},
    & \quad &k\ge1, \\[2pt]
    p_x^\e0 &= \psi_{yyy}^\e0, &
    p_x^\e{2k} &=
      \psi_{xxy}^\e{2k-2} + \psi_{yyy}^\e{2k}, 
      & &k\ge1, \\[2pt]
    p_y^\e0 &= 0, \qquad
    p_y^\e2 = -\psi_{xyy}^\e0,\; &
    p_y^\e{2k} &=
      -\psi_{xxx}^\e{2k-4}-\psi_{xyy}^\e{2k-2},
      & &k\ge2, \\[2pt]
    p^\e{2k}&(x,y) =
      \int_0^x p_x^\e{2k}(\xi,0)\,d\xi + \int_0^y
      p_y^\e{2k}(x,\eta)\,d\eta, \hspace*{-2.5in} & & & 
       & k\ge0.
\end{aligned}
\end{equation}
Equation (\ref{eqn:M:recur}) implies that
$\partial_xp_y^\e{2k}=\partial_yp_x^\e{2k}$ for $k\ge0$; hence,
differentiating under the integral sign in (\ref{eqn:uvp:from:psi}),
we see that $p_x^\e{2k}$ and $p_y^\e{2k}$ actually are the partial
derivatives of $p^\e{2k}$.  Finally, our choice of $Q^\e{2k}$ ensures
$\int_0^1p_x^\e{2k}(\xi,0)\,d\xi=
\int_0^1\psi_{yyy}^\e{2k}(x,0)\,dx=0$ so that $p^\e{2k}$ is periodic.

Using Algorithm~\ref{alg:G:def} to evaluate
$(\psi^\e0,Q^\e0)=G(0,V_0,V_1)$ yields
\begin{align}
  \label{eqn:Q0}
  Q^\e0 &= \frac{V_0+V_1}{2}\frac{I_2}{I_3}, \\
  \label{eqn:psi0}
  \psi^\e0 &= (V_0h)\frac{y}{h} +
  \left(3Q^\e0 - (2V_0+V_1)h\right)\frac{y^2}{h^2} +
  \left(-2Q^\e0 + (V_0+V_1)h\right)\frac{y^3}{h^3},
\end{align}
which agrees with Reynolds' approximation (\ref{eqn:lubr1}) when
$\mb{u}^\e0$, $p^\e0$ are computed from $\psi^\e0$.  To compute
higher order terms in the expansion, we need to study the recursion
(\ref{eqn:psi:Q:recur}) more closely to determine how $h$ will
enter into the formulas for $Q^\e{2k}$ and $\psi^\e{2k}$.

\subsection{Algebraic structure of the stream function expansion}
\label{sec:struc}

In this section, we show how the terms $\psi^\e{2k}$ and $Q^\e{2k}$ in
the stream function expansion depend on $h$.  The key result of this
section is that these higher order corrections have a structure
similar to the zeroth order formulas (\ref{eqn:Q0}) and
(\ref{eqn:psi0}), but the coefficient on each~$\frac{y^n}{h^n}$ now
belongs to a more complicated polynomial algebra in the symbols $V_0$,
$V_1$, $h$, the derivatives of $h$, and certain weighted averages of
the products of $h$ and its derivatives.  We also present a concise
representation for the correction terms using matrices of rational
numbers that are independent of any particular choice of shape
function $h$. By splitting the analysis into one part that holds
universally and another that depends on $h$ in a simple way, we are
able to derive useful error estimates governing the expansion solution
truncated at any order in section~\ref{sec:error}.

Let $\mc{P}=\mbb{Q}[h,h_x,h_{xx},\ldots]$ denote the algebra of
polynomials in $h$ and its derivatives over the rationals.  A typical
element of $\mc{P}$ might be $3+\frac{2}{5}h^2h_{xx}h_{xxx}^3$.  In
$\mc{P}$, the generators $h$, $h_x$, etc.,~are treated as symbols
rather than functions.  Thus, if $h(x)$ happens to equal $1$
identically, the polynomials $1-h$ and $h_x^6$ are nonzero in
$\mc{P}$ even though they are mapped to zero when $\mc{P}$ is
(noninjectively) embedded in $C^\infty(T)$,  the space of $C^\infty$
functions on the periodic interval $T=[0,1]_p$.  If $h$ is not smooth,
its derivatives can still be manipulated symbolically and various
subspaces (involving terms with few enough derivatives) can still be
embedded in actual function spaces such as $C^k(T)$.

For any monomial $\alpha=Ch^{i_0} h_x^{i_1} h_{xx}^{i_2}
\cdots\in\mc{P}$ with $C\ne0$, we define its \emph{superdegree} to be
the number of derivatives present:
\begin{equation}
  \sdeg(\alpha)=i_1 + 2i_2 + 3i_3 + \cdots.
\end{equation}
If $\alpha\in\mc{P}$, we define its superdegree to be the maximal
superdegree of any of its terms, and set $\sdeg(0)=-\infty$.  Since
$\mbb{Q}$ is a field, $\sdeg(\alpha\beta)= \sdeg(\alpha)+\sdeg(\beta)$
for any $\alpha,\beta\in\mc{P}$.  We say that $\alpha$ is homogeneous
of superdegree $k$ if each of its terms has superdegree $k$.

Let $\mc{H}\subset\mc{P}$ denote the subalgebra generated by the set
$\{h^{k-1}\partial_x^kh:k\ge1\}$, i.e.,
\begin{equation}
  \label{eqn:H:def}
  \mc{H}=\mbb{Q}[\{h_x,hh_{xx},h^2h_{xxx},\ldots\}],
\end{equation}
and for $k\ge0$, let $\mc{H}_k\subset\mc{H}$ denote the subspace
\begin{equation}
  \label{eqn:Hk:def}
  \mc{H}_k = \{0\}\cup\{\alpha\in\mc{H}\; : 
\;\alpha\text{ is homogeneous of superdegree $k$}\}.
\end{equation}
Note that $\mc{H}_k$ is finite-dimensional for all $k$, and
$\mc{H}_0=\mbb{Q}$ is the set of constant polynomials.  We will
denote the dimension of $\mc{H}_k$ by
\begin{equation}
  d_k = \dim(\mc{H}_k).
\end{equation}

\begin{figure}[t]
\begin{center}
\fbox{\parbox{4.9in}{
\vspace*{4pt}
\begin{algorithm} $(${\sc basis generation}$)$
\label{alg:Phi}
\begin{center}
\parbox[b][1.15in][c]{.75\linewidth}{
\begin{tabbing}
\hspace*{.25in} \= \hspace*{.25in} \=
\hspace*{1.5in} \= \hspace*{.25in} \= \kill
 \textbf{for} $k=0,\ldots,k_0$ \\
\> $\Phi_k = \{t_1^k\}$, \>\> $(\mbox{or } d_k=1)$ \\
 \textbf{for} $j=2,\ldots,k_0$ \\
\> \textbf{for} $k=j,\ldots,k_0$ \\
\>\> $\Phi_k = \Phi_k\cup t_j\Phi_{k-j}$, \> $(\mbox{or } d_k=d_k+d_{k-j})$ \\
 \textbf{return} $\{\Phi_0,\ldots,\Phi_{k_0}\}$
\end{tabbing}
}
\qquad
\parbox[b][1.15in][c]{.18\linewidth}{
\includegraphics[width=\linewidth]{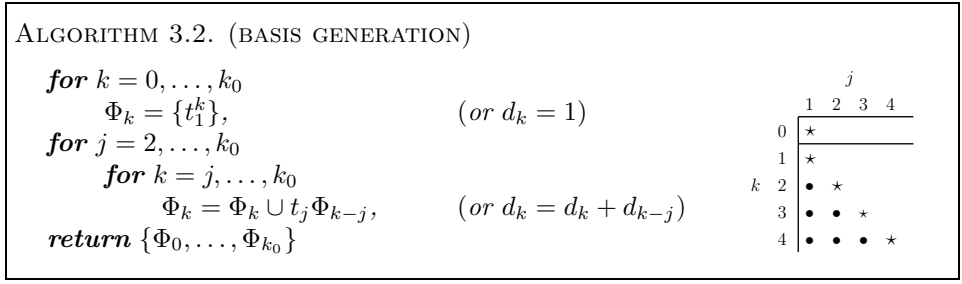}
}
\end{center}
\end{algorithm}
}}
\end{center}
\caption{Algorithm to find a canonical basis $\Phi_k$ for each space
$\mc{H}_k$ in the range $0\le k\le k_0$.  Here $t_1\leftrightarrow
h_x, \ldots, t_k\leftrightarrow \frac{1}{k!}h^{k-1}\partial_x^kh$.}
\label{fig:alg:Phi}
\end{figure}

Given an integer $k_0\ge0$, we can use Algorithm~\ref{alg:Phi} in
Figure~\ref{fig:alg:Phi} to construct a canonical basis $\Phi_k =
\{\varphi^\e{k}_1,\ldots,\varphi^\e{k}_{d_k}\}$ for each $\mc{H}_k$  
with $k$ in the range $0\le k\le k_0$.  For notational convenience,
let $t_j$ stand for $\frac{1}{j!}h^{j-1}\partial_x^jh$.  As the outer
loop (on $j$) progresses, $\Phi_k$ contains a basis for the subspace
of $\mc{H}_k$ that involves only the symbols $t_1, \ldots, t_j$.  Let
us denote these auxiliary sets by
\begin{equation}
  \Phi_{kj} =
  \{t_1^{i_1}\ldots t_j^{i_j} \; : \;
  i_1+2i_2+\cdots+ji_j = k\}, \qquad 1\le j\le k.
\end{equation}
Then
$\Phi_{k1}=\{t_1^k\}$, $\Phi_{kk}=\Phi_k$, and $\Phi_{kj} =
\Phi_{k,j-1}\cup t_j\Phi_{k-j,j}$ for $2\le j\le k$.  In other words,
$\Phi_{kj}$ consists of $\Phi_{k,j-1}$ together with all products of
the variables $t_1, \ldots, t_j$ of superdegree $k$ that contain at
least one power of $t_j$.  The first several $\Phi_k$ and $d_k$ are
given by
\begin{align}
\label{eqn:Phi:024}
    &\Phi_0 = \{1\}, \quad
    \Phi_1 = \{h_x\}, \quad
    \Phi_2 = \left\{h_x^2,\; \frac{hh_{xx}}{2}\right\}, \quad
    \Phi_3 = \left\{h_x^3,\, \frac{hh_xh_{xx}}{2}, \frac{h^2h_{xxx}}{6}\right\},\\
 \notag   &\Phi_4 = \left\{h_x^4,\, \frac{hh_x^2h_{xx}}{2},\,
    \frac{h^2h_{xx}^2}{4},\, \frac{h^2h_xh_{xxx}}{6},\,
    \frac{h^3h_{xxxx}}{24}\right\}, \\
\notag  &(d_0,\ldots,d_{10}) = \{1,1,2,3,5,7,11,15,22,30,42\},
    \quad d_{20} = 627, \quad d_{50} = 204226.
\end{align}
\noindent We have found empirically that the first 75000 terms satisfy
$\frac{1}{2}(\frac{13^{\sqrt{k}}}{6k+1})<d_k<
\frac{13^{\sqrt{k}}}{6k+1}$.  In fact, we have recently learned
of the Hardy--Ramanujan formula
\begin{equation}
  d_k \sim \frac{\exp\left(\pi\sqrt{2k/3}\right)}{4k\sqrt{3}}
  \quad\text{ as } \quad k\rightarrow\infty
\end{equation}
for the number of partitions of the integer $k$.  Thus, rather
than 13, the base is in fact $e^{\pi\sqrt{2/3}}=13.001954$.

We can now describe the structure of the stream function expansion
in terms of the shape function $h$.  In the following theorem,
$\mc{V}_1\mc{H}_{2k}$ is the tensor product of $\mc{V}_1$ and
$\mc{H}_{2k}$, where
\begin{equation}
  \mc{V}_1 = \{0\}\cup\{\alpha\in\mbb{Q}[V_0,V_1]
\,:\,\alpha\text{ is homogeneous of degree 1}\}
\end{equation}
is the space of rational linear combinations of $V_0$ and $V_1$.
Recall from (\ref{eqn:Q}) above that $I_m=\int_0^1 h(x)^{-m}\,dx$.


\begin{theorem}
\label{thm:struc}
The terms $Q^\e{2k},\psi^\e{2k}$ in the stream function expansion
defined by the recursion {\rm(\ref{eqn:psi:Q:recur})} and
Algorithm~{\rm\ref{alg:G:def}} have the form\vspace{5pt}
\begin{equation}
  \label{eqn:Q:form}
  \begin{aligned}
    Q^\e{2k} &= \frac{I_2}{I_3}a^\e{2k}+
    \sum_{\ell=0}^{k-1}Q^\e{2\ell}b^\e{2k-2\ell}, \qquad
    \psi^\e{2k} = \frac{I_2}{I_3}\alpha^\e{2k} +
    \sum_{\ell=0}^k Q^\e{2\ell}\beta^\e{2k-2\ell},
  \end{aligned}
\end{equation}
where
\begin{equation}
  \label{eqn:alpha:formula}
  \alpha^\e{2k}(x,y) = \sum_{n=1}^{2k+3}
  \alpha_n^\e{2k}(x)\frac{y^n}{h(x)^n},
  \qquad
  \beta^\e{2k}(x,y) = \sum_{n=1}^{2k+3}
  \beta_n^\e{2k}(x)\frac{y^n}{h(x)^n},
\end{equation}
and
\begin{equation}
  \label{eqn:alpha:form}
  \alpha_n^\e{2k}\in \frac{I_3}{I_2}h\mc{V}_1\mc{H}_{2k}, \qquad
  \beta_n^\e{2k}\in \mc{H}_{2k}.
\end{equation}
Moreover,
$a^\e{2k}=\frac{1}{2I_3}\int_0^1\frac{\alpha^\e{2k}_3(x)}{h(x)^3}\,dx$ and\,
$b^\e{2k}=\frac{1}{2I_3}\int_0^1\frac{\beta^\e{2k}_3(x)}{h(x)^3}\,dx$.
\end{theorem}

\begin{remark} \label{rk:mat:rep} \upshape
In addition to pinning down the way in which $h$ appears in the
formulas for the stream function expansion, this theorem allows us to
represent $\psi^\e{2k}$ and $Q^\e{2k}$ using matrices of rational
numbers.  Explicitly, (\ref{eqn:alpha:formula}) and
(\ref{eqn:alpha:form}) hold iff there are matrices
$A^\e{2k}$, $B^\e{2k}$ with entries
in $\mc{V}_1$ and $\mbb{Q}$, respectively,
with rows indexed from 0 to $(2k+3)$ and
columns indexed from $1$ to $d_{2k}$, and containing\break\vspace*{-12pt}\pagebreak

\noindent only zeros
in row 0, such that
\begin{equation}
  \label{eqn:mat:rep}
  \begin{aligned}
  \alpha^\e{2k}(x,y) &= \left(Y_{2k}(x,y)\right)^T A^\e{2k}
  \left(\frac{I_3}{I_2}h(x)\Phi_{2k}(x)\right), \\
  \beta^\e{2k}(x,y) &= \left(Y_{2k}(x,y)\right)^T B^\e{2k}
  \Phi_{2k}(x),
  \end{aligned}
\end{equation}
where $Y_{2k}=(1,\frac{y}{h}, \ldots,
(\frac{y}{h})^{2k+3})^T$ and
$\Phi_{2k}=(\varphi^\e{2k}_1,\ldots,\varphi^\e{2k}_{d_{2k}})^T$ are
treated as column vectors.  The purpose of the zeroth row is to make
it easy to convert to orthogonal polynomials in $y/h$ if desired.  The
final statement of the theorem asserts that the formulas for
$a^\e{2k}$ and $b^\e{2k}$ are also encoded in the matrices $A^\e{2k}$
and $B^\e{2k}$.  If we adopt Matlab notation and denote row $i$ of
$A^\e{2k}$ by $A^\e{2k}(i,:)$, then
\begin{equation}
  \label{eqn:a2k:b2k}
  a^\e{2k} = \frac{1}{2}A^\e{2k}(3,:)E^\e{2k}_2, \qquad
  b^\e{2k} = \frac{1}{2}B^\e{2k}(3,:)E^\e{2k}_3,
\end{equation}
where
\begin{equation}
  \label{eqn:Eklm:def}
  E^\e{2k}_m = \left(E^\e{2k}_{m,1},\ldots,E^\e{2k}_{m,d_{2k}}\right)^T
  = \frac{1}{I_m}\int_0^1 \frac{\Phi_{2k}(x)}{h(x)^m}\,dx.
\end{equation}
Note that $E^\e{2k}_{m,j}$ is the weighted average of
$\varphi^\e{2k}_j$  
with weight function $I_m^{-1} h^{-m}$.  For example, $E^\e{0}_{m}=(1)$,
$E^\e{2}_{m}=(\frac{1}{I_m}\int_0^1 \frac{h_x^2}{h^m}\,dx,\,
\frac{1}{I_m}\int_0^1 \frac{hh_{xx}}{2h^m}\,dx)^T$,
etc.; see (\ref{eqn:Phi:024}) above.
\end{remark}

\begin{example}\label{exa:thm} \upshape
We can now represent $Q^\e0$ and $\psi^\e0$ in (\ref{eqn:Q0}) and
(\ref{eqn:psi0}) by
\begin{equation}
  \label{eqn:A0:B0}
  a^\e0 = \frac{V_0+V_1}{2}, \qquad
  A^\e0 = V_0\begin{pmatrix} \phm0 \\ \phm1 \\ -2 \\ \phm1 \end{pmatrix} +
  V_1 \begin{pmatrix} \phm0 \\ \phm0 \\ -1 \\ \phm1 \end{pmatrix}, \qquad
  B^\e0 = \begin{pmatrix} \phm0 \\ \phm0 \\ \phm3 \\ -2 \end{pmatrix}.
\end{equation}
The second order terms $Q^\e2$ and $\psi^\e2$ involve these as well as
\begin{align}
  a^\e2 &= \frac{1}{2}\left[
    V_0\begin{pmatrix} \frac{7}{15}, \,\frac{2}{15} \end{pmatrix} +
    V_1\begin{pmatrix} \frac{19}{30}, \, -\frac{8}{15}
    \end{pmatrix} \right]
 E^\e2_2,
  \qquad
  b^\e2 = \frac{1}{2}\begin{pmatrix} -\frac{6}{5}, \,-\frac{2}{5}\end{pmatrix}
  E^\e2_3, \\[5pt]
  \notag
  A^\e2 &= V_0 \begin{pmatrix}
    \js\phm0 & \js\phm0 \\ 
    \js\phm0 & \js\phm0 \\ 
    \js-8/15 & \js\phm2/15 \\
    \js\phm7/15 & \js\phm2/15 \\
    \js\phm2/3 & \js-2/3 \\
    \js-3/5 & \js\phm2/5
    \end{pmatrix} + V_1 \begin{pmatrix}
    \js\phm0 & \js\phm0 \\ 
    \js\phm0 & \js\phm0 \\ 
    \js-11/30 & \js\phm7/15 \\
    \js\phm19/30 & \js-8/15 \\
    \js\phm1/3 & \js-1/3 \\
    \js-3/5 & \js\phm2/5
    \end{pmatrix}, \quad
  B^\e2 = \begin{pmatrix}
    \js\phm0 & \js\phm0 \\ 
    \js\phm0 & \js\phm0 \\ 
    \js\phm9/5 & \js-2/5 \\
    \js-6/5 & \js-2/5 \\
    \js-3 & \js\phm2 \\
    \js\phm12/5 & \js-6/5
    \end{pmatrix}.
\end{align}
For $k\ge2$, $A^\e{2k}$ and $B^\e{2k}$ are both $(2k+4)\times d_{2k}$
matrices with rows 0 and 1 containing only zeros.  These matrices
are universal: the shape function $h$ enters into the formulas only
through $Y_{2k}$, $\Phi_{2k}$, and $E^\e{2k}_m$ in (\ref{eqn:mat:rep})
and (\ref{eqn:a2k:b2k}).  In Appendix~\ref{sec:impl}, we show how
to compute $A^\e{2k}$ and $B^\e{2k}$ directly from the lower order
matrices $A^\e{2\ell}$ and $B^\e{2\ell}$ with $0\le \ell<k$.
\end{example}

\begin{proofof} {\it of Theorem}~\ref{thm:struc}:
We saw in Example~\ref{exa:thm} above that $Q^\e0$ and $\psi^\e0$ have
the desired form.  Suppose $k_0\ge1$, and the theorem holds for $0\le
k<k_0$.  We must show that it is also true for $k=k_0$.  By
(\ref{eqn:psi:Q:recur}),
\begin{equation}
    \left(\psi^\e{2k_0},Q^\e{2k_0}\right) = \begin{cases}
    G\left(-2\psi_{xxyy}^\e{2k_0-2},\hspace*{.665in}\,
    0,\, g_1^\e{2k_0}\right), & k_0 = 1, \\[6pt]
    G\left(-2\psi_{xxyy}^\e{2k_0-2}-\psi_{xxxx}^\e{2k_0-4},\,
    0,\,g_1^\e{2k_0}\right),  &  k_0 \ge 2.
  \end{cases}
\end{equation}
We will use the second formula for both cases with the understanding
that $\psi^\e{-2}$ should be replaced by zero.  The first step of
Algorithm~\ref{alg:G:def} is to compute $\psi_0^\e{2k_0}$. Using the
induction hypothesis, we obtain
\begin{align}
  \label{eqn:psi0:1}
  \psi_0^\e{2k_0} &= \frac{I_2}{I_3}
  \left(-2\mbb{V}^4\alpha_{xxyy}^\e{2k_0-2}
  - \mbb{V}^4\alpha_{xxxx}^\e{2k_0-4}\right) \\
  \notag
  &\quad+\; \sum_{\ell=0}^{k_0-1}Q^\e{2\ell}
  \left(-2\mbb{V}^4\beta_{xxyy}^\e{2k_0-2\ell-2}\right)
  + \sum_{\ell=0}^{k_0-2} Q^\e{2\ell}
  \left(-\mbb{V}^4\beta_{xxxx}^\e{2k_0-2\ell-4}\right).
\end{align}
The upper limit of the last sum can be replaced by $k_0-1$, since we
interpret $\beta^\e{-2}$ as zero.  We would like to rewrite this in
the form
\begin{equation}
  \label{eqn:psi0:2}
  \psi_0^\e{2k_0} = \frac{I_2}{I_3}\left(\sum_{n=4}^{2k_0+3}
  \alpha_n^\e{2k_0}(x)\frac{y^n}{h^n}\right) +
  \sum_{\ell=0}^{k_0-1}Q^\e{2\ell}\left(\sum_{n=4}^{2k_0-2\ell+3}
  \beta_n^\e{2k_0-2\ell}(x)\frac{y^n}{h^n}\right).
\end{equation}
If we use the induction hypothesis and substitute
(\ref{eqn:alpha:formula}) into (\ref{eqn:psi0:1}), the operator
$\mbb{V}^4\partial_y^2$ annihilates a single power of $y$ and
antidifferentiates higher powers of $y$ twice.  Similarly, $\mbb{V}^4$
antidifferentiates all powers of $y$ four times.  Thus, for $k=k_0$
and $4\le n\le 2k+3$, we should define
\begin{equation}
  \label{eqn:alpha:recur}
  \begin{aligned}
    \alpha_n^\e{2k}(x) &= \frac{-2h^n\partial_x^2
      \left(\alpha_{n-2}^\e{2k-2}h^{-n+2}\right)}{n(n-1)} +
    \frac{-h^n\partial_x^4
      \left(\alpha_{n-4}^\e{2k-4}h^{-n+4}\right)}{n(n-1)(n-2)(n-3)} \\
  \end{aligned}
\end{equation}
with an identical formula for $\beta_n^\e{2k}$ in terms of
$\beta_{n-2}^\e{2k-2}$ and $\beta_{n-4}^\e{2k-4}$.  The second term
should be omitted when $k=1$ or $n=4$, and is zero when $n=5$.  As
part of the induction hypothesis, we may assume
that (\ref{eqn:alpha:recur}) and its $\beta$ version hold for $1\le k<k_0$
as well, so that each term in the sum over $\ell$ in (\ref{eqn:psi0:1})
also has the form described in (\ref{eqn:psi0:2}).
Note that for
$n\ge0$ and any differentiable function $\varphi(x)$,
\begin{equation}
  \partial_x(h^{-n}\varphi) = h^{-(n+1)}(h\partial_x - nh_x)\varphi.
\end{equation}
By Lemmas~\ref{lem:mult:hx} and~\ref{lem:hdx} below, $h\partial_x$ and
multiplication by $h_x$ both map $\mc{H}_k$ to $\mc{H}_{k+1}$ for all
$k\ge0$.  Thus
\begin{align}
  \label{eqn:dx:alpha:hn}
  h^n\partial_x^2\left(\alpha_{n-2}^\e{2k_0-2}h^{-n+2}\right) & \\
  \notag
=  h[h\partial_x - &\,(n - 2)h_x][h\partial_x -
    (n-3)h_x]\left(h^{-1}\alpha_{n-2}^\e{2k_0-2}\right)\in
  \frac{I_3}{I_2}\mc{V}_1 h\mc{H}_{2k_0}, \\
  \notag
  h^n\partial_x^2\left(\beta_{n-2}^\e{2k_0-2}h^{-n+2}\right) &=
  [h\partial_x - (n - 1)h_x][h\partial_x -
    (n-2)h_x]\left(\beta_{n-2}^\e{2k_0-2}\right)\in \mc{H}_{2k_0},
\end{align}
with similar formulas for
$h^n\partial_x^4(\alpha_{n-4}^\e{2k_0-4}h^{-n+4})$ and
$h^n\partial_x^4(\beta_{n-4}^\e{2k_0-4}h^{-n+4})$.
We conclude that $\alpha_n^\e{2k}$ and $\beta_n^\e{2k}$ have the form
claimed in (\ref{eqn:alpha:form}) when $k=k_0$ and $4\le n\le 2k_0+3$.
Finally, we obtain $Q^\e{2k_0}$ and $\psi^\e{2k_0}$ from
$\psi_0^\e{2k_0}$ in (\ref{eqn:psi0:2}) using
Algorithm~\ref{alg:G:def}.  They satisfy (\ref{eqn:Q:form}) and
(\ref{eqn:alpha:formula}) if we set $k=k_0$ and define
$\alpha_1^\e{2k}=0$, $\beta_1^\e{2k}=0$,
\begin{equation}
  \label{eqn:alpha2}
  \begin{aligned}
    \alpha_2^\e{2k}(x) &= \sum_{n=4}^{2k+3}(n-3)\alpha_n^\e{2k}(x) -
    V_1{-1/2\choose k}h_x^{2k}, \\
    \beta_2^\e{2k}(x) &= \sum_{n=4}^{2k+3}(n-3)\beta_n^\e{2k}(x), \\
    \alpha_3^\e{2k}(x) &= \sum_{n=4}^{2k+3}(2-n)\alpha_n^\e{2k}(x) +
    V_1{-1/2\choose k}h_x^{2k}, \\
    \beta_3^\e{2k}(x) &= \sum_{n=4}^{2k+3}(2-n)\beta_n^\e{2k}(x),
    \end{aligned}
\end{equation}
$a^\e{2k}=\frac{1}{2I_3}\int_0^1\frac{\alpha^\e{2k}_3(x)}{h(x)^3}\,dx$, 
and
$b^\e{2k}=\frac{1}{2I_3}\int_0^1\frac{\beta^\e{2k}_3(x)}{h(x)^3}\,dx$.
As part of the induction hypothesis, we may assume (\ref{eqn:alpha2})
also holds for $1\le k<k_0$.  The factors of $n$ in (\ref{eqn:alpha2})
are due to the terms $\pm\psi_{0,y}(x,h)h$ in the formula for
$\psi^\e{2k_0}$ in Algorithm~\ref{alg:G:def}.
The terms $3Q^\e{2k_0}\frac{y^2}{h^2}$
and $-2Q^\e{2k_0}\frac{y^3}{h^3}$ in the formula for $\psi^\e{2k_0}$
are accounted for in (\ref{eqn:Q:form})
by extending the upper limit of the sum over $\ell$ from $k_0-1$ to $k_0$
and noting that $\beta^\e{0}(x,y) =
3\frac{y^2}{h^2}-2\frac{y^3}{h^3}$.  Thus, $\psi^\e{2k_0}$ and $Q^\e{2k_0}$
have the desired form, and $\alpha_n^\e{2k_0}$, $\beta_n^\e{2k_0}$
belong to the appropriate spaces, as claimed.\end{proofof}

To complete this proof, we need two simple lemmas about the spaces
$\mc{H}_k$ (which also serve as the foundation for our numerical
algorithm described in Appendix~\ref{sec:impl}).

\begin{lemma}
  \label{lem:mult:hx}
  If $k\ge0$ and $\varphi\in\mc{H}_k$, then $h_x\varphi\in\mc{H}_{k+1}$.
\end{lemma}
\begin{proof}
This follows easily from the definition of $\mc{H}_k$ in (\ref{eqn:Hk:def}).\qquad\end{proof}

\begin{lemma}
\label{lem:hdx}
  If $k\ge0$ and $\varphi\in\mc{H}_k$, then $h\partial_x\varphi\in\mc{H}_{k+1}$.
\end{lemma}
\begin{proof}
If $k=0$, then $h\partial_x\varphi=0\in\mc{H}_{k+1}$.  Suppose $k_0\ge1$, 
and the result holds for $k<k_0$. Let $\varphi\in\mc{H}_{k_0}$ be a
monomial.  Then there is a $k\in\{1,\ldots,k_0\}$ and a monomial
$\beta\in\mc{H}_{k_0-k}$ such that
$\varphi=(h^{k-1}\partial_x^kh)\beta$.  But then
\begin{equation}
  h\partial_x\varphi = (k-1)h_x\varphi +
  \left(h^k\partial_x^{k+1}h\right)\beta + \left(h^{k-1}\partial_x^kh\right)(h\partial_x\beta).
\end{equation}
Evidently, all three terms belong to $\mc{H}_{k_0+1}$, the third due
to the induction hypothesis.  This result can now be applied term by
term for any polynomial $\varphi\in\mc{H}_k$.\qquad\end{proof}

\section{Error analysis}
\label{sec:error}

To estimate the error of the expansion of $\psi$ and $Q$ through order
$2k$, we show that the truncation error satisfies a weak form of the
rescaled biharmonic equation (\ref{eqn:psi:bvp}) with data ($\mb{F}$,
$g_0$, $g_1$) of order $\veps^{2k+2}$.  We also prove a uniform
coercivity result for the family of bilinear forms involved in the
weak formulation, which allows us to bound the truncation error in
terms of the data.

Throughout this section, we will treat $\Omega$ and $T=[0,1]_p$ as
$C^\infty$ manifolds by identifying the points
\begin{equation}
\begin{aligned}
  \Omega:&  & (0,y)&\sim(1,y), \qquad 0<y<h(0), \\
	T:& & 0&\sim1
\end{aligned}
\end{equation}
and adding a coordinate chart to each that ``wraps around.''  In
particular: a function in $C^k(\Omega)$ or $C^k(T)$ is understood to
have $k$ continuous periodic derivatives;
$\partial\Omega=\Gamma_0\cup\Gamma_1$; $\partial T=\varnothing$; the
support of a function $\phi\in C^k_c(\Omega)$ vanishes near $\Gamma_0$
and $\Gamma_1$ but not necessarily at $x=0$ and $x=1$; and the Sobolev
spaces $H^k(\Omega)$ and $H^k_0(\Omega)$ are the completions of
$C^k(\overline\Omega)$ and $C^k_c(\Omega)$ in the $\|\cdot\|_k$ norm 
and thus contain only $x$-periodic functions with appropriate smoothness
at $x=0,1$.

\subsection{Weak formulation of the rescaled biharmonic equation}

An interesting difference between the biharmonic equation and the
Poisson equation is that the boundary conditions in the latter are
completely specified in the problem statement, whereas one of them
(the flux $Q$) in the former problem must be determined as part of the
solution.  The integral condition (\ref{eqn:int:psi:yyy}), which
uniquely determines $Q$, must also be reformulated weakly, since it
involves more than two derivatives of $\psi$.  This can be done
\cite{girault} by slightly enlarging the space of test functions to
include functions that are constant along $\Gamma_1$ (rather than
equal to 0 there).  To this end, we define
\begin{equation}
  \Psi = \left\{\phi\in H^2(\Omega)\;\; :\;\;
  (\phi,\partial_y\phi)\left\vert_{\Gamma_0}=(0,0), \;\;
  (\phi,\partial_y\phi)\right\vert_{\Gamma_1}=(\text{const},0)\right\}.
\end{equation}
For $\phi$, $\psi$ in $H^2(\Omega)$, we define the bilinear form
\begin{equation}
  \begin{aligned}
  a_\veps(\psi,\phi) &= \int_\Omega \psi_{yy}\phi_{yy} +
  2\veps^2\psi_{xy}\phi_{xy} + \veps^4\psi_{xx}\phi_{xx}\,dA \\ &=
  a^\e0(\psi,\phi) + \veps^2a^\e2(\psi,\phi) + \veps^4a^\e4(\psi,\phi).
  \end{aligned}
\end{equation}
To obtain estimates that hold uniformly in $\veps$, it will be useful
to work with the weighted norms and seminorms
\begin{equation}
  \begin{gathered}
  \|\psi\|_0^2 = \int_\Omega \psi^2\,dA, \quad
  |\psi|_{1,\veps}^2 = \int_\Omega \psi_y^2 + (\veps\psi_x)^2\,dA, \quad
  |\psi|_{2,\veps}^2 = a_\veps(\psi,\psi), \\
\|\psi\|_{1,\veps}=\sqrt{\|\psi\|_0^2 + |\psi|_{1,\veps}^2}, \quad
\|\psi\|_{2,\veps}=\sqrt{\|\psi\|_0^2+|\psi|_{1,\veps}^2+
|\psi|_{2,\veps}^2}.
  \end{gathered}
\end{equation}
For fixed $\veps$, these norms are equivalent to the usual Sobolev
norms in which $\veps$ is set to 1 in these expressions.
We use $x$ to parametrize functions defined on $\Gamma_0$ or
$\Gamma_1$ and define the weighted boundary norm
\begin{equation}
  \label{eqn:H12}
  \|g\|^2_{1/2,\veps}
  = \sum_{k=-\infty}^\infty
  \left[1+(2\pi k\veps)^2\right]^{1/2}|c_k|^2, \qquad
  c_k = \int_0^1 g(x)e^{-2\pi i k x}\,dx.
\end{equation}
We equip the dual spaces $\Psi'$ and
$H^{-1}(\Omega)^2=[H^1_0(\Omega)^2]'$ with the weighted norms
\begin{equation}
  \label{eqn:neg:norms}
  \|l\|_{-2,\veps} = \sup_{\|\psi\|_{2,\veps}=1}|\langle l,\psi\rangle|,
  \qquad
  \|\mb{F}\|_{-1,\veps}=\sup_{\|u\|_{1,\veps}^2 +
    \|\veps v\|_{1,\veps}^2=1}
  |\langle\mb{F},(u,v)\rangle|.
\end{equation}
Since
$\|\psi\|_{2,\veps}^2 \ge \|\psi_y\|_{1,\veps}^2 +
\|\veps\psi_x\|_{1,\veps}^2$, the linear functional
$\langle l,\psi\rangle =
\langle \mb{F},\nabla\times\psi\rangle$ on $\Psi$ satisfies
$\|l\|_{-2,\veps}\le\|\mb{F}\|_{-1,\veps}$.

\begin{definition}[weak solutions]
Suppose
\begin{equation}
  h\in C^{1,1}(T), \qquad
  \mb{F}\in H^{-1}(\Omega)^2, \qquad
  g_0\in H^{1/2}(\Gamma_0), \qquad
  g_1\in H^{1/2}(\Gamma_1).
\end{equation}
We say that $(\psi,Q)\in H^2(\Omega)\times\mbb{R}$ is a
weak solution of {\rm(\ref{eqn:psi:bvp})--(\ref{eqn:int:psi:yyy})} if
\begin{equation}
  \label{eqn:weak:soln}
  a_\veps(\psi,\phi) = \langle \mb{F},\nabla\times\phi\rangle
\end{equation}
for all $\phi\in\Psi$ and the boundary conditions
\begin{equation}
  \label{eqn:B:psi2}
  B\psi = (0,g_0,Q,g_1)
\end{equation}
hold in the trace sense, where
$B\psi:=(\psi\vert_{\Gamma_0},\psi_y\vert_{\Gamma_0},
\psi\vert_{\Gamma_1},\psi_y\vert_{\Gamma_1})$.
\end{definition}

\begin{prop} Every classical solution is a weak solution.
\end{prop}

\begin{proof}
We assume $\psi\in C^4(\overline{\Omega})$ and
(\ref{eqn:psi:bvp})--(\ref{eqn:int:psi:yyy}) hold classically; this
requires additional regularity for $\mb{F}$, $g_0$, $g_1$, of course.
If we multiply (\ref{eqn:psi:bvp}) by a test function $\phi\in
C^2(\overline{\Omega})\cap\Psi$ and use the identity
$\chi(\nabla\times\mb{v}) = \nabla\times(\chi\mb{v}) +
(\nabla\times\chi)\cdot\mb{v}$,
we obtain
\begin{align}
 \label{eqn:by:parts}
    0 &= \int_\Omega \phi(-\Delta_\veps^2\psi+\nabla\times\mb{F})\,dA \\
 \notag &=\int_\Omega \left(\nabla\times[
      \phi(\Delta_\veps\mb{u}_\veps+\mb{F})] +
    (\nabla\times\phi)\cdot[
      \Delta_\veps\mb{u}_\veps+\mb{F}] \right)\,dA \\
  \notag
  &= \int_{\Gamma_0-\Gamma_1}\phi(\Delta_\veps\mb{u}_\veps+\mb{F})\cdot
  \mb{t}\,ds +
  \int_\Omega \left[(\nabla\times\phi)_\veps\cdot
  (\nabla\times\Delta_\veps\psi) + (\nabla\times\phi)\cdot\mb{F}
  \right]\,dA \\
\notag
  &= \int_{\Gamma_0-\Gamma_1}[\phi(\cdots)
  - (\Delta_\veps\psi)(\nabla\times\phi)_\veps]\cdot\mb{t}\,ds
  + \int_\Omega\left[-(\Delta_\veps\phi)(\Delta_\veps\psi) +
  (\nabla\times\phi)\cdot\mb{F}\right]\,dA,
\end{align}
where $(\nabla\times\phi)_\veps=(\phi_y,-\veps^2\phi_x)$ and
the curves $\Gamma_0$ and $\Gamma_1$ are both oriented from left
to right as in Figure~\ref{fig:geom}.  The conditions
\begin{equation}
  \label{eqn:Psi:bcs}
  \begin{aligned}
  \phi\left\vert_{\Gamma_0}=0, \qquad
  \partial_y\phi\right\vert_{\Gamma_0} = 0,\\
  \phi\left\vert_{\Gamma_1}=\text{const}, \qquad
  \partial_y\phi\right\vert_{\Gamma_1} = 0,
  \end{aligned}
\end{equation}
ensure that the boundary terms are zero: the first boundary term is
equal to
\begin{equation}
  (\phi\big\vert_{\Gamma_1})[p(1,h(1))-p(0,h(0))]=0
\end{equation}
(with $p$
as in (\ref{eqn:path:indep}), where it was shown to be periodic), and
the second is zero since $\nabla\times\phi=0$ on $\Gamma_0$ and
$\Gamma_1$.  One more integration by parts gives $\int_\Omega
(\Delta_\veps\phi)(\Delta_\veps\psi)\,dA = a_\veps(\psi,\phi)$, so
(\ref{eqn:weak:soln}) holds.  Since $C^2(\overline{\Omega})\cap\Psi$
is dense in $\Psi$ and both sides of (\ref{eqn:weak:soln}) are bounded
linear functionals of $\phi\in \Psi$, this formula holds for all
$\phi\in\Psi$.\qquad\end{proof}

\subsection{Uniform coercivity}

The following two theorems are the key to obtaining error estimates
for the expansion solutions of section~\ref{sec:expand}.

\begin{theorem}
  The bilinear form $a_\veps(\cdot,\cdot)$ is coercive on $\Psi$
  (uniformly in $\veps$) with respect to the weighted norm
  $\|\cdot\|_{2,\veps}$, i.e., there is a constant $\alpha>0$ such that
  $\alpha\|\psi\|_{2,\veps}^2 \le a_\veps(\psi,\psi)$ for all
  $\veps>0$, $\psi\in\Psi$.
\end{theorem}

\begin{proof}
Without loss of generality, we may assume the characteristic height
$\bar{H}$ of the domain was chosen so that $0<h(x)\le1$ for $0\le
x\le1$.  We now use a standard Poincar\'{e}--Friedrichs
argument \cite{braess}.  Suppose $\psi\in
C^2(\overline{\Omega})\cap\Psi$.  Then
\begin{align}
   |\psi(x,y)|^2 = \left|\int_0^y \psi_y(x,\eta)\,d\eta\right|^2 &\le
  y\int_0^{h(x)}|\psi_y(x,\eta)|^2\,d\eta, \\
   |\psi(x,y)|^2 = \left|\int_0^y (y-\eta)
  \psi_{yy}(x,\eta)\,d\eta\right|^2 &\le
  \frac{y^3}{3}\int_0^{h(x)}|\psi_{yy}(x,\eta)|^2\,d\eta.
\end{align}
Integrating over $\Omega$,\vspace{5pt}
\begin{equation}
  \|\psi\|_0^2\le\frac{h_{\max}^2}{2} \int_\Omega|\psi_y|^2\,dA \le
  \frac{1}{2}|\psi|_{1,\veps}^2, \qquad
  \|\psi\|_0^2\le\frac{h_{\max}^4}{12} \int_\Omega|\psi_{yy}|^2\,dA \le
  \frac{1}{12}|\psi|_{2,\veps}^2.
\end{equation}
Repeating this argument on the derivatives of $\psi$ yields
\begin{equation}
  |\psi|_{1,\veps}^2 = \|\psi_y\|_0^2 + \|\veps\psi_x\|_0^2
  \le \frac{1}{2}\left( |\psi_y|_{1,\veps}^2 + 
  |\veps\psi_x|_{1,\veps}^2 \right)
  = \frac{1}{2}|\psi|_{2,\veps}^2 = \frac{1}{2} a_\veps(\psi,\psi)
\end{equation}
so that $\|\psi\|_{2,\veps}^2\le \frac{19}{12}a_\veps(\psi,\psi)$.  Since
$C^2(\overline{\Omega})\cap\Psi$ is dense in $\Psi$, we conclude that
$(12/19)\|\psi\|_{2,\veps}^2\le a_\veps(\psi,\psi)$ for all
$\psi\in\Psi$ as claimed.\qquad\end{proof}

\begin{theorem} \label{thm:psi:bound}
A weak solution $\psi$ of the boundary value problem
{\rm(\ref{eqn:psi:bvp})--(\ref{eqn:int:psi:yyy})} exists and is unique.
Moreover, the following estimate holds:
\begin{align}
  \label{eqn:inhom:est}
  \|\psi\|_{2,\veps} \,&\le\,\, 
  \frac{19}{12}\|\mb{F}\|_{-1,\veps} \,\,\, \\
  \notag
  &\quad+\left(\hspace*{-2pt}72 + 860\left(\hspace*{-2pt}
  \begin{gathered}\veps^2\|h_x\|^2_\infty 
    + \veps^4\|h_x\|^4_\infty \\
  +\, 4\veps^4\left\|\frac{1}{2} \jd hh_{xx}\right\|^2_\infty
  \end{gathered}\hspace*{-2pt}\right)\hspace*{-3pt}\right)^{1/2}\hspace*{-6pt}
  \left(\left\|h^{-1/2}g_0\right\|_{1/2,\veps} \hspace*{-2pt}+
  \left\|h^{-1/2}g_1\right\|_{1/2,\veps}\right)\hspace*{-2pt}.
\end{align}
In particular, if\, $\veps\le \frac{r_0}{3}$ with $r_0^{-1}=\max(\|h_x\|_\infty,
\|\frac{1}{2}hh_{xx}\|_\infty^{1/2})$, then
\begin{equation}
  \label{eqn:inhom:est2}
  \|\psi\|_{2,\veps} \le 
  \frac{19}{12}\|\mb{F}\|_{-1,\veps} +
  15\left(\left\|h^{-1/2}g_0\right\|_{1/2,\veps} +
  \left\|h^{-1/2}g_1\right\|_{1/2,\veps}\right).
\end{equation}
\end{theorem}
\unskip

\begin{proof}
We begin by constructing a function $\psi_0\in H^2(\Omega)$ that
satisfies the boundary conditions (\ref{eqn:psi:bcs}) with $Q=0$.
First, we map the domain $\Omega$ to the $x$-periodic unit square
$R=T\times(0,1)$ via the transformation
\begin{equation}
  \wtil{\psi}_0(x,y) = h(x)^{-3/2}\psi_0(x, h(x)y),
  \qquad 0\le x \le 1, \; 0<y<1.
\end{equation}
We include $h^{-3/2}$ here to avoid powers of $h_0^{-1}$
in (\ref{eqn:psi0:est}), where $h_0=\min_{0\le x\le1}h(x)$.
We require $\wtil{\psi}_0(x,0)=0$,\,
$\wtil{\psi}_{0,y}(x,0)=h(x)^{-1/2}g_0(x)$,\,
$\wtil{\psi}_0(x,1)=0$, and
$\wtil{\psi}_{0,y}(x,1)=h(x)^{-1/2}g_1(x)$.
To construct such a function, we define $\zeta\in C^1(\mbb{R})$ via
%
\begin{equation*}
  \zeta(y) = \left\{\begin{array}{lr@{\,\,}l}
     0, &   y& \le-1 \\
     y+2y^2+y^3, &  -1& \le y\le 0 \\
     y-2y^2+y^3, &  \phm 0& \le y \le 1 \\
     0, &  1& \le y
    \end{array}
  \right\} \;\;
  \parbox[c][.7in][b]{2.58in}{
  \includegraphics[width=2.55in]{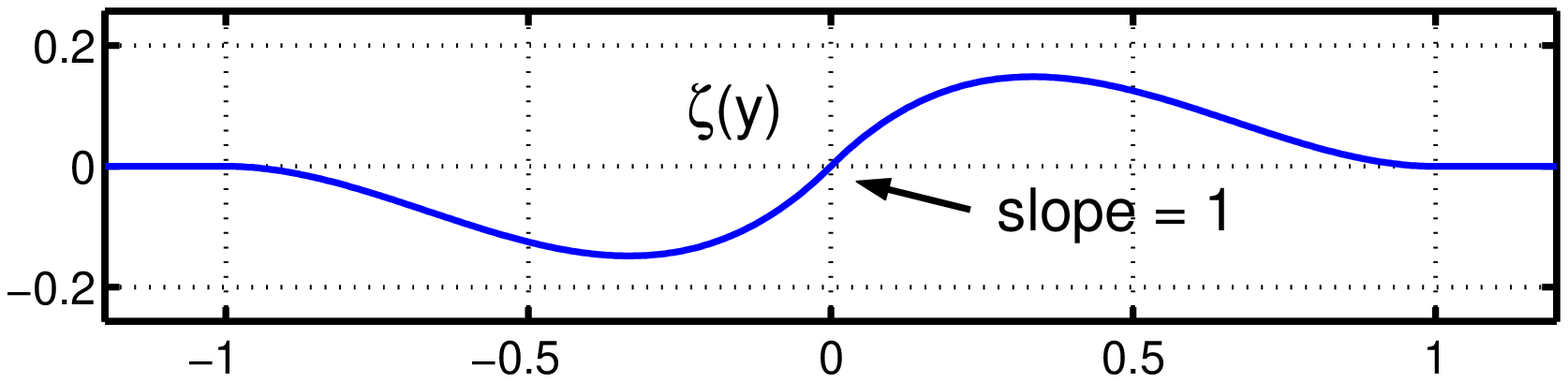}}
\end{equation*}
and set
%
\begin{equation}
  \label{eqn:extension}
  \wtil{\psi}_0(x,y) = \sum_{k=-\infty}^\infty
  \left(c_k\frac{\zeta(\langle k\rangle y)}{\langle k\rangle} +
  d_k\frac{\zeta(\langle k\rangle (y-1))}{\langle k\rangle}\right)
  e^{2\pi ikx},
\end{equation}
where
%
\begin{equation}
  \langle k\rangle = \left[1+(2\pi k\veps)^2\right]^{1/2}, \qquad
  [c_k,d_k] = \int_0^1 [g_0,g_1](x)h(x)^{-1/2}e^{-2\pi ikx}\,dx.
\end{equation}
The value and slope of $\zeta$ at $y=0$ and $|y|\ge1$ ensure that
$\wtil\psi_0$ satisfies the
desired boundary conditions.  Assume for the moment that each $d_k$ is
zero (i.e., $g_1\equiv0$).  Let $S$ be the strip $T\times\mbb{R}$.  We
may use (\ref{eqn:extension}) to define $\wtil\psi_0$ on all of $S$
and take its Fourier transform
\begin{equation}
  \label{eqn:ft:psi0}
  \left(\wtil\psi_0\right)^\wedge(k,\eta) = \int_0^1 \int_{-\infty}^\infty
  \wtil\psi_0(x,y)e^{-2\pi i(kx+\eta y)}\,dy\,dx =
    \frac{c_k\hat\zeta(\eta/\langle k\rangle)}{\langle k\rangle^2}.
\end{equation}
Since $\zeta$ is antisymmetric and supported on $[-1,1]$, we have
\begin{equation}
  \begin{aligned}
  2\left\|\wtil\psi_0\right\|^2_{2,\veps,R} &=
  \left\|\wtil\psi_0\right\|^2_{2,\veps,S} \,\le\, \left\|\wtil\psi_0\right\|^2_{0,S}
  + 2\left|\wtil\psi_0\right|^2_{1,\veps,S} + \left|\wtil\psi_0\right|^2_{2,\veps,S} \\
  &= \sum_k\int_{-\infty}^\infty\left[1+(2\pi k\veps)^2 + (2\pi \eta)^2\right]^2\,
  \left|\left(\wtil\psi_0\right)^\wedge(k,\eta)\right|^2\,d\eta \\
  &= 
  \sum_k \langle k\rangle|c_k|^2
  \int_{-\infty}^\infty\left[1+(2\pi t)^2\right]^2
  \left|\hat{\zeta}(t)\right|^2\,dt = \frac{898}{105}
  \left\|h^{-1/2}g_0\right\|^2_{1/2,\veps}.
  \end{aligned}
\end{equation}
A similar argument works if we assume $g_0\equiv0$, but $g_1\not\equiv0$.
Thus, on $R$, we have
\begin{equation} \label{eqn:wtil:est}
   \left\|\wtil\psi_0\right\|_{2,\veps} \le \js \sqrt{\dst\frac{449}{105}} \jd
  \left(\left\|h^{-1/2}g_0\right\|_{1/2,\veps} +
  \left\|h^{-1/2}g_1\right\|_{1/2,\veps}\right).
\end{equation}
Next we use the formula $\psi_0(x,y) =
h(x)^{3/2}\,\wtil\psi_0(x,\frac{y}{h(x)})$ to obtain
\begin{equation}
  \label{eqn:psi0:derivs}
  \begin{aligned}
  \psi_{0,y} =& h^{1/2}\wtil\psi_{0,y}, \qquad
  \psi_{0,yy} = h^{-1/2}\wtil\psi_{0,yy}, \\
  \psi_{0,x} =& h^{3/2}\wtil\psi_{0,x}
  - yh^{-1/2}h_x\wtil\psi_{0,y}
  + \frac{3}{2}h^{1/2}h_x\wtil\psi_0, \\
  \psi_{0,xy} =&
    h^{1/2}\wtil\psi_{0,xy}
  - yh^{-3/2}h_x\wtil\psi_{0,yy}
  + \frac{1}{2}h^{-1/2}h_x\wtil\psi_{0,y}, \\
  \psi_{0,xx} =&
    h^{3/2} \wtil\psi_{0,xx}
  + 3 h^{1/2} h_x \wtil\psi_{0,x} 
  - 2yh^{-1/2}h_x\wtil\psi_{0,xy}
  + \frac{3}{2} h^{1/2} h_{xx}\wtil\psi_0 \\ &
  + y^2h^{-5/2}h_x^2\wtil\psi_{0,yy}
  - yh^{-3/2}h_x^2\wtil\psi_{0,y}
  - yh^{-1/2}h_{xx}\wtil\psi_{0,y}
  + \frac{3}{4} h^{-1/2}h_x^2\wtil\psi_0.
  \end{aligned}
\end{equation}
Using Lemma \ref{lem:sumsq} below and $0<h(x)\le1$, we find that
\begin{align}
  \int_\Omega \psi_0^2\,dA &=
  \int_0^1\int_0^1 \psi_0\left(x,h(x)y\right)^2h(x)\,dx\,dy
  = \int_R h^4\wtil\psi_0^2\,dA
  \le\int_R \wtil\psi_0^2\,dA, \\[5pt]
\notag  \int_\Omega \psi_{0,y}^2\,dA &\le
  \int_R \wtil\psi_{0,y}^2\,dA, \qquad
  \int_\Omega \psi_{0,yy}^2\,dA =
  \int_R \wtil\psi_{0,yy}^2\,dA, \\[5pt]
  \notag
  \int_\Omega\veps^2\psi_{0,x}^2\,dA &\le
  2\int_R\veps^2\wtil\psi_{0,x}^2\,dA
  \,+\, 4\veps^2\|h_x^2\|_\infty\left[
    \int_R \wtil\psi_{0,y}^2\,dA +
    \frac{9}{4}\int_R\wtil\psi_0^2\,dA\right], \\[5pt]
  \notag
  \int_\Omega 2\veps^2\psi_{0,xy}^2\,dA &\le
  2\int_R 2\veps^2\wtil\psi_{0,xy}^2\,dA
  \,+\, 4\veps^2\|h_x^2\|_\infty\left[
    2\int_R \wtil\psi_{0,yy}^2\,dA +
    \frac{1}{2}\int_R \wtil\psi_{0,y}^2\,dA\right], \\[3pt] \notag
  \int_\Omega \veps^4\psi_{0,xx}^2\,dA &\le \,\,
  \frac{5}{2}\int_R \veps^4\wtil\psi_{0,xx}^2\,dA
  \,+\, 30 \veps^2\|h_x^2\|_\infty\int_R \veps^2\wtil\psi_{0,x}^2\,dA
  \\[-6pt]
  \notag
  & \quad+\, 30\veps^2\|h_x^2\|_\infty\int_R2\veps^2\wtil\psi_{0,xy}^2\,dA
  \,+\, 30 \veps^4\|h^2h^2_{xx}\|_\infty\int_R \wtil\psi_0^2\,dA\\\notag
&\quad+\, 30\veps^4\|h_x^4\|_\infty\int_R\wtil\psi_{0,yy}^2\,dA 
+ 30\veps^4\|h_x^4\|_\infty\int_R\wtil\psi_{0,y}^2\,dA\\\notag
&  \quad +\, 30\veps^4\|h^2h_{xx}^2\|_\infty\int_R\wtil\psi_{0,y}^2\,dA
  \,+\, 30\veps^4\|h_x^4\|_\infty\int_R\wtil\psi_0^2\,dA \jd.
\end{align}
Note that the inverse powers of $h$ in (\ref{eqn:psi0:derivs}) are
canceled when we change variables from $\Omega$ to $R$ by the
factors of $h$ that arise from the substitution $y\rightarrow hy$ and
from the Jacobian of the transformation.  Collecting terms and
majorizing, we obtain
\begin{equation}
  \label{eqn:psi0:est}
  \begin{aligned}
  \|\psi_0\|_{2,\veps} &\le \left(
  5/2 \,+\, 30\veps^2\|h_x\|^2_\infty \,+\, 30\veps^4\|h_x\|^4_\infty \,+\,
  30\veps^4\|hh_{xx}\|^2_\infty\right)^{1/2}\left\|\wtil\psi_0\right\|_{2,\veps}.
  \end{aligned}
\end{equation}
Finally, we correct $\psi_0$ by a function in $\Psi$ to obtain the weak
solution $\psi$, which must satisfy
\begin{equation} \label{eqn:weak:inhom}
  \psi-\psi_0\in\Psi, \qquad
  a_\veps(\psi-\psi_0,\phi)=\langle l,\phi\rangle :=
  \langle\mb{F},\nabla\times\phi\rangle
  - a_\veps(\psi_0,\phi)   \quad\mbox{for all } \phi\in\Psi.
\end{equation}
Since $l$ is a bounded linear functional on $\Psi$, the Lax--Milgram
theorem implies existence and uniqueness of the solution $\psi$ of
(\ref{eqn:weak:inhom}) and gives the error bound
\begin{equation}
\|\psi-\psi_0\|_{2,\veps} \le \alpha^{-1}\|l\|_{2,\veps} \le
\frac{19}{12}\left(\|\mb{F}\|_{-1,\veps}+\|\psi_0\|_{2,\veps}\right).
\end{equation}
Combining this with (\ref{eqn:wtil:est}) and (\ref{eqn:psi0:est}) and
using the triangle inequality gives (\ref{eqn:inhom:est}), where we
note that $\frac{5}{2}(\frac{31}{12})^2
(\frac{449}{105})\le 72$ and
$30(\frac{31}{12})^2
(\frac{449}{105})\le 860$.\qquad\end{proof}

The following lemma was used to balance the coefficients in the
terms of (\ref{eqn:psi0:derivs}) as much as possible.

\begin{lemma}
\label{lem:sumsq}
For any $a_1,\ldots,a_8\in\mbb{R}$,
\begin{align}
  (a_1+a_2+a_3)^2 &\le\,
    2a_1^2\,+\,4a_2^2 \,+\, 4a_3^2, \\
    \notag
  (a_1+\cdots+a_8)^2  &\le\,
  \frac{5}{2}a_1^2\,+\,\frac{10}{3}a_2^2\,+\,15a_3^2\,+
  \,\frac{40}{3}a_4^2\,+\,
  30\left(a_5^2+a_6^2+a_7^2\right)\,+\,\frac{160}{3}a_8^2.
\end{align}
\end{lemma}
\unskip

\begin{proof}
In general, given positive real numbers $\gamma_1, \ldots, \gamma_n$
such that $\sum_1^n \gamma_j^{-1}\le 1$, then for all $a\in\mbb{R}^n$, 
we have $\left(\sum_1^n a_j\right)^2 \le \sum_1^n \gamma_j a_j^2$.
This is a consequence of the Cauchy--Schwarz inequality:
\begin{equation} 
  \left(\sum_j a_j\right)^2 =
  \left(\sum_j \left(\gamma_j^{-1/2}\right)\left(\gamma_j^{1/2}a_j\right)
  \right)^2
  \le \left(\sum_j \gamma_j^{-1}\right)\left(\sum_j\gamma_j a_j^2\right).
\end{equation}
One readily checks that $\frac{1}{2}+\frac{1}{4}+\frac{1}{4}=1$
and $(\frac{2}{5}+\cdots+\frac{3}{160})=\frac{461}{480}\le1$.\qquad\end{proof}

\subsection{Truncation error equation}

In section~\ref{sec:expand}, we showed how to construct successive
terms in the stream function expansion by solving the recursion
\mbox{(\ref{eqn:M:recur})--(\ref{eqn:int:psi2k:yyy})}.
Theorem~\ref{thm:struc} guarantees that derivatives of $h$ higher than
$2k$ do not appear in the formulas for $\psi^\e0,\ldots,\psi^\e{2k}$;
hence, if $h\in C^{2k}(T)$, these functions satisfy
(\ref{eqn:M:recur})--(\ref{eqn:int:psi2k:yyy}) in the classical sense
(with $k$ replaced by $\ell$ and running from $0$ to $k$ instead of
\mbox{$0$ to $\infty$}).  Thus, if $h\in C^{2k+4}$,
$\psi^\e{2k}_\tem{approx}=\psi^\e{0}+\veps^2\psi^\e{2}+
\cdots+\veps^{2k}\psi^\e{2k}$ satisfies
\begin{equation}
  \label{eqn:approx:psi}
  \Delta_\veps^2\psi^\e{2k}_\tem{approx} =
  \veps^{2k+2}\left(2\psi^\e{2k}_{xxyy}+\psi^\e{2k-2}_{xxxx}\right)
  + \veps^{2k+4}\psi^\e{2k}_{xxxx}.
\end{equation}
The truncation error
$\psi^\e{2k}_\tem{err}=\psi_\tem{exact}-\psi^\e{2k}_\tem{approx}$ then
satisfies
$\Delta_\veps^2\psi_\tem{err}=-\Delta_\veps^2\psi_\tem{approx}$ with
$O(\veps^{2k+2})$ boundary data.  Since the right-hand side of
(\ref{eqn:approx:psi}) and the boundary data are known in terms of
$h$, we are able to estimate the size of $\psi_\tem{err}^\e{2k}$ using
Theorem~\ref{thm:psi:bound} above.  However, to use this theorem, we
need to formulate (\ref{eqn:approx:psi}) weakly.

We begin by showing that the $\psi^\e{2\ell}$ satisfy a weak version
of the recursion (\ref{eqn:M:recur}).  Suppose $k\ge0$ and $h\in
C^{2k}(T)$.  Let $\phi\in\Psi$ and denote the constant value of $\phi$
on $\Gamma_1$ by $q$.  We multiply both sides of (\ref{eqn:M:recur})
by $\phi$ and integrate by parts using
\begin{align}
 \label{eqn:by:parts:phi:q}
  \int_\Omega\phi\,\psi^\e{2\ell}_{yyyy}\,dA &=
  \int_\Omega\phi_{yy}\psi^\e{2\ell}_{yy}dA +
  q\int_0^1 \psi^\e{2\ell}_{yyy}(x,h(x))\,dx, \\
  \notag
  2 \int_\Omega\phi\,\psi^\e{2\ell}_{xxyy}\,dA &=
  2\int_\Omega\phi_{xy}\psi^\e{2\ell}_{xy}dA +
  q\int_0^1\left[\psi^\e{2\ell}_{xxy}(x,h(x)) -
    \psi^\e{2\ell}_{xyy}(x,h(x))h_x\right]\,dx, \\
\notag  \int_\Omega\phi\,\psi^\e{2\ell}_{xxxx}\,dA &=
  \int_\Omega\phi_{xx}\psi^\e{2\ell}_{xx}dA -
  q\int_0^1 \psi^\e{2\ell}_{xxx}(x,h(x))h_x\,dx
\end{align}
to obtain the recursion
\begin{align}
  \notag
  a^\e0\left(\psi^\e0,\phi\right) &= 0, \\
  \label{eqn:a:recur}
  a^\e0\left(\psi^\e2,\phi\right) &= -a^\e2\left(\psi^\e0,\phi\right), \\
  \notag
  a^\e0\left(\psi^\e{2\ell},\phi\right) &= -a^\e2\left(\psi^\e{2\ell-2},\phi\right)
  -a^\e4\left(\psi^\e{2\ell-4},\phi\right), \quad \ell=2,3,\ldots,k.
\end{align}
By (\ref{eqn:uvp:from:psi}), the boundary terms in
(\ref{eqn:by:parts:phi:q}) combine to form
\begin{equation}
  q\int_0^1\left[p_x^\e{2\ell}(x,h(x))+
  p_y^\e{2\ell}(x,h(x))h_x\right]\,dx = 0, \qquad 0\le\ell\le k,
\end{equation}
when substituted into (\ref{eqn:M:recur}).  Other boundary terms do
not arise in (\ref{eqn:by:parts:phi:q}), since $\phi=0$ on $\Gamma_0$
and $\phi_x=\phi_y=0$ on $\Gamma_0$ and $\Gamma_1$.

Now suppose $h\in C^{2k+1,1}(T)$, i.e., $h(x)$ has
$2k+1$ continuous periodic derivatives and $\partial_x^{2k+1}h$ is
Lipschitz continuous so that $\partial_x^{2k+2}h\in L^\infty(T)$.
Let $(\psi_\tem{exact},Q_\tem{exact})$ be the weak solution of
\begin{equation}
  \Delta_\veps^2\psi=0, \qquad B\psi=(0,g_0,Q,g_1),
\end{equation}
with $g_0$, $g_1$ given in (\ref{eqn:gdef}), and define the truncation
errors and approximate solutions
\begin{equation}
  \label{eqn:ek:def}
  \begin{aligned}
    \psi_\tem{err}^\e{2k} &= \psi_\tem{exact}-\psi_\tem{approx}^\e{2k},
    & \quad \psi_\tem{approx}^\e{2k} &=
    \psi^\e0 + \veps^2\psi^\e2 + \cdots + \veps^{2k}\psi^\e{2k}, \\
    Q_\tem{err}^\e{2k} &= Q_\tem{exact}-Q_\tem{approx}^\e{2k},
    & \quad Q_\tem{approx}^\e{2k} &=
    Q^\e0 + \veps^2Q^\e2 + \cdots + \veps^{2k}Q^\e{2k}.
  \end{aligned}
\end{equation}
Since $\psi^\e0,\ldots,\psi^\e{2k}$ satisfy (\ref{eqn:a:recur}) and
(\ref{eqn:B:recur}) while $a_\veps(\psi_\tem{exact},\phi)=0$ for every
$\phi\in\Psi$, we may expand $a_\veps(\psi_\tem{err}^\e{2k},\phi)$ in
powers of $\veps$ to obtain the truncation error equation
\begin{equation}
  \label{eqn:error:bvp}
  \begin{aligned}
    a_\veps\left(\psi_\tem{err}^\e{2k},\phi\right) &= 
    -\veps^{2k+2}\langle \mb{F}_k,\nabla\times\phi\rangle,
    \qquad \phi\in\Psi, \\
  B\psi_\tem{err}^\e{2k} &=
  \left(0,0,Q_\tem{err}^\e{2k},\veps^{2k+2}\gamma_k\right),
  \end{aligned}
\end{equation}
where 
\begin{align}
  \gamma_k &= \veps^{-2k-2}\left(
  g_1 - \left[g_1^\e0 + \veps^2 g_1^\e2 + \cdots +
                   \veps^{2k} g_1^\e{2k}\right]\right)
  \intertext{and}
  \notag
  \left\langle
    \mb{F}_k, \begin{pmatrix}\phi_y\\-\phi_x\end{pmatrix}
  \right\rangle &=
  \begin{cases}
    a^\e2\left(\psi^\e0,\phi\right) +
    \veps^2 a^\e4\left(\psi^\e0,\phi\right), & k=0, \\
    a^\e2\left(\psi^\e{2k},\phi\right) +
    \veps^2 a^\e4\left(\psi^\e{2k},\phi\right) +
    a^\e4\left(\psi^\e{2k-2},\phi\right),\; & k\ge1.
  \end{cases}
\end{align}
There are many functionals $\mb{F}_k\in H^{-1}(\Omega)^2$ that have
this action on the subspace
\begin{equation}
  V = \{\nabla\times\phi\;:\;\phi\in\Psi\}
  =\left\{(u,v)\in H^1_0(\Omega)^2\;:\; u_x+v_y=0\right\}.
\end{equation}
For example, $\mb{F}_k=\veps^{-2k-2}[\nabla p_\tem{approx}^\e{2k}
  - (\Delta_\veps\mb{u}_\tem{approx}^\e{2k})_\veps]$ satisfies
$\Delta_\veps^2\psi_\tem{approx}^\e{2k}=\veps^{2k+2}\nabla\times
\mb{F}_k$ classically and, using (\ref{eqn:uvp:from:psi}), may be
shown to have following action on $H^1_0(\Omega)^2$:
\begin{equation}
\label{eqn:f2k:def0}
\begin{aligned}
  \langle \mb{F}_k,(u;v)\rangle &=
  \int_\Omega \left(\psi_{xx}^\e0\right) \left(u_y - \veps^2 v_x\right) -
  \left(\veps^{-1}\psi^\e0_{xy}\right)(\veps v_y)\,dA,
  & &k=0, \\
  \langle \mb{F}_k,(u;v)\rangle &=
  \int_\Omega \left(\begin{gathered}
    \left(\psi_{xx}^\e{2k}\right) \left(u_y - \veps^2 v_x\right)
  -\left(\veps^{-1}\psi^\e{2k}_{xy}\right)(\veps v_y) \\
  -\left(\veps^{-2}\psi_{xx}^\e{2k-2}\right)\left(\veps^2 v_x\right)\end{gathered}\right)
  dA, & \;\; &k\ge1.
\end{aligned}
\end{equation}
This choice is suboptimal because the
terms $\veps^{-1}\psi_{xy}^\e{2k}$ and $\veps^{-2}\psi_{xx}^\e{2k-2}$
diverge as $\veps\rightarrow0$.
We grouped $\veps$ with $v_y$ and $\veps^2$ with $v_x$
due to the definition
(\ref{eqn:neg:norms}) of $\|\mb{F}_k\|_{-1,\veps}$.  Instead, we will
use the following functional, which agrees with (\ref{eqn:f2k:def0})
on $V$:
\begin{align}
  \notag
  \langle \mb{F}_k,(u;v)\rangle &=
  \int_\Omega \left(\psi_{xx}^\e0\right) \left(2u_y - \veps^2 v_x\right)\,dA,
  \hspace*{1.196in}\;\;\;k=0, \\
\label{eqn:f2k:def}
  \langle \mb{F}_k,(u;v)\rangle &=
  \int_\Omega \left(\psi_{xx}^\e{2k}\right) \left(2u_y - \veps^2 v_x\right)
  +h^2\psi_{xxxx}^\e{2k-2}\widetilde\phi[u]\,dA, \quad\;\;\; k\ge1, \\
  \notag
  &\widetilde\phi[u](x,y) :=
  -\int_y^{h(x)}\frac{u(x,\eta)}{h(x)^2}\,d\eta =
  \int_y^{h(x)}\frac{(\eta-y)u_y(x,\eta)}{h(x)^2}\,d\eta.
\end{align}
Note that if $\phi\in\Psi$ and $q=\phi\vert_{\Gamma_1}$, then
$\widetilde\phi[\phi_y]=(\phi-q)h^{-2}$.  The purpose of the $h^{-2}$
here is to be able to include an $h^2$ with $\psi_{xxxx}^\e{2k-2}$ in
the error estimates below (to properly consolidate terms).  Another
alternative to (\ref{eqn:f2k:def}) that would lead to similar estimates
below is $\brak{\mb{F}_k,(u;v)}=\int_\Omega [-\veps^2\psi_{xx}^\e{2k}
v_x-\psi_{yy}^\e{2k+2}u_y]\,dA$.

\subsection{Error estimates}
\label{sec:error:est}

Let us assume from now on that $\veps\le\frac{r_0}{3}$ with 
$r_0=\max(\|h_x\|_\infty,
\|\frac{1}{2}hh_{xx}\|_\infty^{1/2})^{-1}$.
Then by
(\ref{eqn:error:bvp}) and Theorem~\ref{thm:psi:bound}, we have
\begin{equation}
  \label{eqn:e2k0:bound1}
  \left\|\psi_\tem{err}^\e{2k}\right\|_{2,\veps} \le \veps^{2k+2}\left(
  \frac{19}{12}\|\mb{F}_k\|_{-1,\veps} +
  15\left\|h^{-1/2}\gamma_k\right\|_{1/2,\veps}\right).
\end{equation}
It remains to bound the norms of $\mb{F}_k$ and $\gamma_k$
in terms of $h$.  From (\ref{eqn:f2k:def}), we have
\begin{equation}
  \label{eqn:phitil1}
  \left|\wtil\phi[u](x,y)\right|^2
  \le\left(\int_y^{h(x)}\frac{(\eta-y)^2}{h(x)^4}\,d\eta\right)
  \left(\int_0^{h(x)}|u_y(x,\eta)|^2\,d\eta\right),
\end{equation}
where the first integral is $\frac{(h-y)^3}{3h^4}$ and the second
is independent of $y$.  Hence
\begin{equation}
  \label{eqn:phitil2}
  \left\|\wtil\phi[u]\right\|_0^2 =
  \int_0^1\int_0^{h(x)}\left|\wtil\phi[u](x,y)\right|^2\,dy\,dx
  \le\frac{1}{12}\|u_y\|_0^2.
\end{equation}
From (\ref{eqn:f2k:def}), we then have
\begin{align}
\notag
   \left|\langle\mb{F}_k,(u;v)\rangle\right| &\le
 \left(2a + \frac{b}{\sqrt{12}}\right) \|u_y\|_0 + a\left\|\veps^2 v_x\right\|_0 \\
\label{eqn:fk:uv:bd}
 &\le
  \left(5a^2 + \frac{4ab}{\sqrt{12}} + \frac{b^2}{12}\right)^{1/2}
  \left(\|u\|_{1,\veps}^2+\|\veps v\|_{1,\veps}^2\right)^{1/2},
\end{align}
where $a=\|\psi_{xx}^\e{2k}\|_0$ and
$b=\|h^2\psi_{xxxx}^\e{2k-2}\|_0$.  Using $\frac{4ab}{\sqrt{12}}\le
\frac{5}{4}a^2 + \frac{4}{15}b^2$, we find that
\begin{equation}
  \left\|\mb{F}_k\right\|_{-1,\veps}\le\sqrt{\frac{25}{4}a^2 + 
\frac{7}{20}b^2}\le\frac{5}{2}a+\frac{3}{5}b,\qquad
  \frac{19}{12}\left\|\mb{F}_k\right\|_{-1,\veps}\le 4a+b.
\end{equation}
Finally, by (\ref{eqn:e2k0:bound1}), we obtain
\begin{equation}
  \label{eqn:e2k0:bound2}
  \left\|\psi_\tem{err}^\e{2k}\right\|_{2,\veps} \le \veps^{2k+2}\left(
  4\left\|\psi_{xx}^\e{2k}\right\|_0 +
  \left\|h^2\psi_{xxxx}^\e{2k-2}\right\|_0 +
  15\left\|h^{-1/2}\gamma_k\right\|_{1/2,\veps}\right),
\end{equation}
where the fourth derivative term should be omitted when $k=0$.
In section~\ref{sec:pressure}, the following bound will also
prove useful:\vspace{5pt}
\begin{equation}
  \label{eqn:e2k0:bound3}
  \left\|\psi_\tem{err}^\e{2k}\right\|_{2,\veps} 
  +\veps^{2k+2}\left\|\psi^\e{2k}_{xx}\right\|_0
  \le \veps^{2k+2}\left(
  5\left\|\psi_{xx}^\e{2k}\right\|_0 +
  \left\|h^2\psi_{xxxx}^\e{2k-2}\right\|_0 +
  15\left\|\frac{\gamma_k}{h^{1/2}}\right\|_{1/2,\veps}\right).
\end{equation}
The truncation error in the flux expansion satisfies
\begin{equation}
  Q_\tem{err}^\e{2k} = \psi_\tem{err}^\e{2k}(x,h(x)) =
  \int_0^{h(x)}(h(x)-\eta)\der{^2\psi_\tem{err}^\e{2k}}{y^2}(x,\eta)\,d\eta
\end{equation}
for any $x\in T$. Using estimates similar to (\ref{eqn:phitil1}) and
(\ref{eqn:phitil2}), we find that
\begin{equation}
  \left|Q_\tem{err}^\e{2k}\right| \le \frac{1}{\sqrt{3}}
  \left\|\psi_\tem{err}^\e{2k}\right\|_{2,\veps}.
\end{equation}
Thus, we have reduced the problem of bounding the truncation errors to
that of checking the norms of three quantities that can be computed
explicitly in closed form.

We begin by attacking the boundary term
$\|h^{-1/2}\gamma_k\|_{1/2,\veps}$.  This norm was
defined in (\ref{eqn:H12}) above.  Recall that
\begin{equation}
  \label{eqn:gamma:2k0}
  \gamma_k=V_1\left(\left[1+\veps^2h_x^2\right]^{-1/2} -
  \sum_{\ell=0}^{k}{-1/2 \choose \ell}\left(\veps^2 h_x^2\right)^\ell\right)\veps^{-2k-2},
\end{equation}
where ${-1/2 \choose 0}=1$ and for $\ell\ge1$,
\begin{equation}
  \label{eqn:binomial}
  {-1/2 \choose \ell} = \frac{(-1/2)(-3/2)\cdots(-[\ell-1/2])}{(1)(2)\cdots(\ell)}
  = (-1)^\ell\,\frac{1}{2}\cdot\frac{3}{4}\cdot\frac{5}{6}\cdots\frac{2\ell-1}{2\ell}.
\end{equation}
Taking the logarithm of this product and its inverse, one may show that
\begin{equation}
  \frac{1}{\sqrt{4\ell+1}}\le\left|{-1/2\choose \ell}\right|
  \le\frac{1}{\sqrt{3\ell+1}}, \qquad \ell=0,1,2,\ldots.
\end{equation}
Since $h\in C^{2k+1,1}(T)\subset C^{1,1}(T)$, we know
$\gamma_k$ is at least Lipschitz continuous on $T$ and so
belongs to $H^1(T)$.  As a result,
\begin{align}
  \notag
  \left\|\frac{\gamma_k}{h^{1/2}}\right\|_{1/2,\veps}^2 &\le
  \left\|\frac{\gamma_k}{h^{1/2}}\right\|_{1,\veps}^2 =
  \int_0^1 h^{-1}\gamma_k^2 + \veps^2\left(
    -\frac{h^{-3/2}}{2} h_x\gamma_k
    + h^{-1/2}\gamma_{k,x}
    \right)^2\,dx \\
    \label{eqn:gamma2k:bound}
    &\le \int_0^1 h^{-1}\gamma_k^2 +
    \frac{5}{4}\veps^2 h^{-3}h_x^2 \gamma_k^2
    + \frac{5}{4}\veps^2h^{-1}\gamma_{k,x}^2
    \, dx.
\end{align}
Since $\veps\|h_x\|_\infty\le1/3<1$, the binomial expansion of
$[1+\veps^2h_x^2]^{-1/2}$ converges uniformly on $T=[0,1]_p$.  As the
terms in this expansion alternate in sign, the error in truncating the
series is smaller (pointwise) than the first omitted term.  Therefore,
for each $x\in T$,
\begin{equation}
\label{eqn:gamma2k:bound2}
  |\gamma_k(x)|\le V_1\left|{-1/2 \choose k+1}\right|h_x(x)^{2k+2}
  \le \frac{V_1}{\sqrt{3k+4}}h_x(x)^{2k+2}.
\end{equation}
Since $(-2\ell){-1/2 \choose \ell} = {-3/2
\choose \ell-1}$, by differentiating (\ref{eqn:gamma:2k0}) we obtain
\begin{equation}
  \gamma_{k,x} = V_1\left(\left[1+\veps^2h_x^2\right]^{-3/2} -
  \sum_{\ell=0}^{k-1}{-3/2\choose \ell}\left(\veps^2h_x^2\right)^\ell\right)
  \left(-\veps^2h_x h_{xx}\right)\veps^{-2k-2}.
\end{equation}
The terms in the expansion of $[1+\veps^2h_x^2]^{-3/2}$ also alternate
in sign, so
\begin{equation}
  \label{eqn:gamma2k:bound3}
  |\gamma_{k,x}(x)|\le V_1\left|{-3/2 \choose k}
  h_x(x)^{2k+1}h_{xx}(x)\right|
  \le V_1\frac{2(k+1)}{\sqrt{3k+4}}\left|h_x(x)^{2k+1}h_{xx}(x)\right|.
\end{equation}
Combining (\ref{eqn:gamma2k:bound}), (\ref{eqn:gamma2k:bound2}), and
(\ref{eqn:gamma2k:bound3}) and using
$\frac{(k+1)^2}{3k+4}\le\frac{1}{4}+\frac{1}{3}k$, we find that
\begin{align}
\label{eqn:bc:effect}
  &\left\|\frac{\gamma_k}{h^{1/2}}\right\|_{1/2,\veps}^2 \le
  \frac{V_1^2}{3k+4}
  \int_0^1 \frac{h_x^{4k+4}}{h} + \frac{5}{4}\frac{\veps^2h_x^2}{h^3}
  \left[h_x^{4k+4} + 16(k+1)^2 h_x^{4k}\left(
    \frac{hh_{xx}}{2}\right)^2\right]dx \\
\notag
  &\qquad\le \,\,\frac{V_1^2}{4} I_1 \wtil{E}^\e{2k+2}_{1,1} +
    V_1^2\veps^2\|h_x\|_\infty^2\left[
    \frac{5}{16}I_3\wtil{E}^\e{2k+2}_{3,1} +
    \left(5 + \frac{20}{3}k\right)
     I_3\wtil{E}^\e{2k+2}_{3,2}\right] \\
\notag
    & \qquad \le \,\, V_1^2 I_1 \left[\frac{1}{4} +
      \frac{\veps^2\|h_x\|_\infty^2}{I_1/I_3} \left(
      \frac{85}{16} + \frac{20}{3}k\right)\right]
    \max_{(m,j)\in\{(1,1),(3,1),(3,2)\}} \wtil{E}^\e{2k+2}_{m,j},
\end{align}
where
\begin{equation}
  I_m = \int_0^1 \frac{1}{h(x)^m}\,dx, \qquad
  \wtil{E}^\e{2k}_{m,j} = \frac{1}{I_m} \int_0^1
  \frac{\left(\varphi^\e{2k}_j(x)\right)^2}{h(x)^m}\,dx.
\end{equation}
Recall\vspace*{-2pt} that $\Phi_{2k}=\{
\varphi_1^\e{2k},\ldots,\varphi_{d_{2k}}^\e{2k} \}=
\{h_x^{2k},\,\frac{1}{2}hh_x^{2k-2}h_{xx},\ldots,
\frac{1}{(2k)!}h^{2k-1}\partial_x^{2k}h\}$ is a basis\vspace*{-1pt} for the
space $\mc{H}_{2k}$ defined in (\ref{eqn:Hk:def}).  $\wtil{E}^\e{2k}_{m,j}$
is the square of the 2-norm (or second moment)\vspace*{-2pt} of
$\varphi^\e{2k}_j(x)$ with respect to the probability measure
$I_m^{-1}h(x)^{-m}\,dx$, whereas $E^\e{2k}_{m,j}$ in
(\ref{eqn:Eklm:def}) is the expected value of $\varphi^\e{2k}_j(x)$
with respect to this measure.
Using (\ref{eqn:bc:effect}) in
(\ref{eqn:e2k0:bound2}) gives a bound on the error caused by failing to
satisfy the boundary conditions in the stream function expansion.  It
is perhaps surprising that this bound can be expressed in terms of
three simple integrals involving $h$ and its derivatives.

\begin{remark} \upshape
$I_1$ and $I_3$ are dimensionless quantities in (\ref{eqn:bc:effect})
--- if $h$ and $x$ still carried dimensions of length, an extra length
scale (e.g., $h_\text{max}$, which is currently set to $1$) would need
to be included in the Sobolev norms to allow $\|\psi\|_0^2$,
$|\psi|_{1,\veps}^2$, and $|\psi|_{2,\veps}^2$ to be added together;
this length scale would also appear in (\ref{eqn:bc:effect}) to
nondimensionalize $I_1/I_3$ in the denominator.
\end{remark}

Finally, we estimate $\|\psi^\e{2k}_{xx}\|_0$ and
$\|h^2\psi^\e{2k-2}_{xxxx}\|_0$ in (\ref{eqn:e2k0:bound2})
in terms of $h$.
It will be useful in our analysis to define
\begin{equation}
  \label{eqn:r:def}
  r_k = \left(\,\max_{1\le \ell\le 2k+2}
  \left\{\left\|\frac{1}{\ell!}h^{\ell-1}
  \partial_x^\ell h\right\|_\infty^{1/\ell}\right\}\right)^{-1}
\end{equation}
so that for $0\le\ell\le k+1$,
$1\le j\le d_{2\ell}$, and $m=1,2,3$, we have
\begin{equation}
  \label{eqn:r:prop}
  \left|\varphi_j^\e{2\ell}(x)\right|\le r_k^{-2\ell}, \qquad
  \left|E^\e{2\ell}_{m,j}\right|\le r_k^{-2\ell}, \qquad
  0\le \wtil{E}^\e{2\ell}_{m,j}\le r_k^{-4\ell}.
\end{equation}
If $h$ is real analytic as well as periodic, a standard contour
integral argument shows that there is an $r>0$ such that
$\|\partial_x^k h\|_\infty\le k!\,r^{-k}$ for all $k\ge0$.  Such
an $r$ serves as a common lower bound for $r_k$ in (\ref{eqn:r:def})
for all $k\ge0$.  If $h$ is a constant function, the results below
hold with $r_k=\infty$ and $r_k^{-1}=0$ (i.e., the lubrication
approximation is exact).

Our first task will be to bound the growth of the terms $Q^\e{2k}$ in
the expansion of the flux.  By Theorem~\ref{thm:struc} and
Remark~\ref{rk:mat:rep}, there are rational matrices $A_0^\e{2k}$,
$A_1^\e{2k}$, $B^\e{2k}$ with rows indexed from 0 to $2k+3$ and
columns indexed from $1$ to $d_{2k}$ such that
\begin{equation}
  \label{eqn:Q:form2}
  Q^\e{2k} = \frac{I_2}{I_3} a^\e{2k} +
  \sum_{\ell=0}^{k-1}Q^\e{2\ell}b^\e{2k-2\ell}, \qquad\quad
  k\ge0,
\end{equation}
where $a^\e{2k}=\frac{1}{2}[V_0A_0^\e{2k}(3,:) + V_1A_1^\e{2k}(3,:)]
E_2^\e{2k}$ and $b^\e{2k}=\frac{1}{2}B^\e{2k}(3,:)E_3^\e{2k}$.  See
Example~\ref{exa:thm} above for a reminder of how this works.  We now
use (\ref{eqn:r:prop}) together with the fact that $|v\cdot w|\le
\|v\|_1\|w\|_\infty$ for $v,w\in\mbb{R}^d$ to conclude that
\begin{equation}
  \left|a^\e{2\ell}\right|\le
  \left(|V_0|\kappa_0^\e{2\ell} +
  |V_1|\kappa_1^\e{2\ell}\right)r_k^{-2\ell}, \quad\;
  \left|b^\e{2\ell}\right| \le \kappa_2^\e{2\ell}r_k^{-2\ell},
  \quad\; 0\le \ell\le k,
\end{equation}
where\vspace{5pt} 
\begin{equation}
  \kappa_i^\e{2k}=\frac{1}{2}
  \sum_{j=1}^{d_{2k}}\left|A_i^\e{2k}(3,j)\right|, \quad
  i=0,1, \qquad
  \kappa_2^\e{2k}=\frac{1}{2}
  \sum_{j=1}^{d_{2k}}\left|B^\e{2k}(3,j)\right|, \qquad k\ge0.
\end{equation}
It follows
from (\ref{eqn:Q:form2}) that if we increase $\kappa_0^\e{2k}$,
$\kappa_1^\e{2k}$ via the  loop
\begin{align}
\notag
\mbox{{\bf for }} &k=1,2,3,\ldots \\
\label{eqn:kappa:loop}
&\kappa_i^\e{2k} = \kappa_i^\e{2k} + \sum_{\ell=0}^{k-1}
\kappa_i^\e{2\ell}\kappa_2^\e{2k-2\ell},\qquad\quad (i=0,1), \hspace*{.75in}
\end{align}
then
\begin{equation}
  \label{eqn:Q:bound}
  \left|Q^\e{2\ell}\right| \le \frac{I_2}{I_3}\left( |V_0| \kappa_0^\e{2\ell}
  + |V_1| \kappa_1^\e{2\ell} \right) r_k^{-2\ell}, \qquad (k\ge0,\;
  0\le \ell\le k).
\end{equation}
The constants $\kappa_0^\e{2k}$, $\kappa_1^\e{2k}$, $\kappa_2^\e{2k}$
do not depend on $h$ and may be computed once and for all along with
the matrices $A_i^\e{2k}$ and $B^\e{2k}$; see Table~\ref{tbl:kappa}.

\begin{table}[t]
\begin{center}
\caption{$\kappa_i^\e{2k}$ before and after loop
  {\rm (\ref{eqn:kappa:loop})}.} \label{tbl:kappa}
\scriptsize
\begin{tabular}{rccccc}
  \hline \\[-8pt]
 $k$ & $\kappa_0^\e{2k}$\rule{0pt}{8pt} before & $\kappa_1^\e{2k}$ before &
  $\kappa_2^\e{2k}$ & $\kappa_0^\e{2k}$ after & $\kappa_1^\e{2k}$ after
\\[1pt]
  \hline \\[-8pt]
 0 & $5.00\times10^{-01}$\rule{0pt}{8pt} & $5.00\times10^{-01}$ & $1.00\times10^{+00}$ & $5.00\times10^{-01}$ & $5.00\times10^{-01}$ \\
 1 & $3.00\times10^{-01}$ & $5.83\times10^{-01}$ & $8.00\times10^{-01}$ & $7.00\times10^{-01}$ & $9.83\times10^{-01}$ \\
 2 & $5.30\times10^{-01}$ & $7.05\times10^{-01}$ & $1.73\times10^{+00}$ & $1.96\times10^{+00}$ & $2.36\times10^{+00}$ \\
 3 & $2.72\times10^{+00}$ & $3.73\times10^{+00}$ & $6.74\times10^{+00}$ & $8.87\times10^{+00}$ & $1.07\times10^{+01}$ \\
 4 & $1.83\times10^{+01}$ & $3.32\times10^{+01}$ & $4.14\times10^{+01}$ & $5.43\times10^{+01}$ & $7.32\times10^{+01}$ \\
 5 & $2.00\times10^{+02}$ & $3.69\times10^{+02}$ & $4.55\times10^{+02}$ & $5.28\times10^{+02}$ & $7.30\times10^{+02}$ \\
 6 & $3.41\times10^{+03}$ & $6.32\times10^{+03}$ & $7.22\times10^{+03}$ & $8.00\times10^{+03}$ & $1.13\times10^{+04}$ \\
 7 & $7.77\times10^{+04}$ & $1.66\times10^{+05}$ & $1.54\times10^{+05}$ & $1.68\times10^{+05}$ & $2.63\times10^{+05}$ \\
 8 & $2.69\times10^{+06}$ & $5.23\times10^{+06}$ & $4.69\times10^{+06}$ & $5.31\times10^{+06}$ & $7.98\times10^{+06}$ \\
 9 & $1.26\times10^{+08}$ & $2.31\times10^{+08}$ & $1.94\times10^{+08}$ & $2.31\times10^{+08}$ & $3.40\times10^{+08}$ \\
10 & $6.51\times10^{+09}$ & $1.45\times10^{+10}$ & $9.97\times10^{+09}$ & $1.18\times10^{+10}$ & $2.00\times10^{+10}$ \\
\hline
\end{tabular}
\end{center}
\end{table}

Now that the terms $Q^\e{2\ell}$ have been bounded, we are ready to
estimate $\|\psi^\e{2k}_{xx}\|_0$ and
$\|h^2\psi^\e{2k-2}_{xxxx}\|_0$.  Recall that
\begin{equation}
  \label{eqn:psi:form2}
  \psi^\e{2k}(x,y) = \frac{I_2}{I_3}\alpha^\e{2k}(x,y) +
  \sum_{\ell=0}^k Q^\e{2\ell} \beta^\e{2k-2\ell}(x,y),
\end{equation}
where $\alpha^\e{2k}(x,y)$ and $\beta^\e{2k}(x,y)$ have the matrix
representations (\ref{eqn:mat:rep}).  A slight modification of the
proof of Theorem~\ref{thm:struc} shows that there are also matrices
$\dot{A}_0^\e{2k}$, $\dot{A}_1^\e{2k}$, $\dot{B}^\e{2k}$ 
and $\ddot{A}_0^\e{2k}$, $\ddot{A}_1^\e{2k}$, $\ddot{B}^\e{2k}$ of
dimension $(2k+4)\times d_{2k+2}$ and $(2k+4)\times d_{2k+4}$,
respectively, such that
\begin{align}
  \notag
    \frac{I_2}{I_3}\alpha^\e{2k}_{xx}(x,y) &= h(x)^{-1}
    \left(Y_{2k}(x,y)\right)^T
    \left[V_0\dot{A}_0^\e{2k} + V_1\dot{A}_1^\e{2k} \right]
      \Phi_{2k+2}(x), \\
      \notag
    \beta^\e{2k}_{xx}(x,y) &= h(x)^{-2}\left(Y_{2k}(x,y)\right)^T
    \dot{B}^\e{2k}\Phi_{2k+2}(x), \\
    \notag
    \frac{I_2}{I_3}\alpha^\e{2k}_{xxxx}(x,y) &=
    h(x)^{-3}\left(Y_{2k}(x,y)\right)^T
    \left[V_0\ddot{A}_0^\e{2k} + V_1\ddot{A}_1^\e{2k} \right]
    \Phi_{2k+4}(x), \\
  \label{eqn:A:dot}
    \beta^\e{2k}_{xxxx}(x,y) &= h(x)^{-4}\left(Y_{2k}(x,y)\right)^T
    \ddot{B}^\e{2k}\Phi_{2k+4}(x),
\end{align}
where $Y_{2k}= (1,\frac{y}{h},\ldots,
(\frac{y}{h})^{2k+3})^T$.  We can achieve significantly
sharper estimates of $\|\psi_{xx}^\e{2k}\|_0$ and
$\|h^2\psi_{xxxx}^\e{2k-2}\|_0$ by expressing the dependence
of $\psi$ on $y$ using orthogonal polynomials.  Let
\begin{equation}
  \wtil{Y}_{2k} =\left(
  \wtil{P}_0\left(\frac{y}{h}\right),
  \wtil{P}_1\left(\frac{y}{h}\right),  \ldots,
  \wtil{P}_{2k+3}\left(\frac{y}{h}\right)\right)^T
\end{equation}
be the vector of shifted Legendre polynomials \cite{abramowitz}, which
satisfy
\begin{equation}
  \label{eqn:shifted:Leg1}
  \int_0^1 \wtil{P}_{m}(x)\wtil{P}_n(x)\,dx=\frac{\delta_{mn}}{2n+1},
  \qquad m,n\ge0.
\end{equation}
The first several are
\begin{equation}
  \label{eqn:shifted:Leg2}
  \wtil{P}_0 = 1, \quad
  \wtil{P}_1 = 2x-1, \quad
  \wtil{P}_2 = 6x^2-6x+1, \quad
  \wtil{P}_3 = 20x^3-30x^2+12x-1.
\end{equation}
The well-known recurrence \cite{abramowitz}
\begin{equation}
  \wtil{P}_n(x) = \frac{2n-1}{n} (2x-1)
  \wtil{P}_{n-1}(x) -
\frac{n-1}{n}\wtil{P}_{n-2}(x), \qquad n\ge 2,
\end{equation}
can be used to construct a nested family of lower triangular matrices
$R_{2k}$ of dimension $(2k+4)\times(2k+4)$ with indices
starting at zero such that
\begin{equation}
  \wtil{Y}_{2k} = R_{2k}Y_{2k}, \qquad
  Y_{2k}^T = \wtil{Y}_{2k}^{T}R_{2k}^{-T}.
\end{equation}
For example,
\begin{equation}
  R_0 = \begin{pmatrix}
    \phm1 & \phm0 & \phm0 & 0 \\
    -1 & \phm2 & \phm0 & 0 \\
    \phm1 & -6 & \phm6 & 0 \\
    -1 & \phm12 & -30 & 20
    \end{pmatrix}, \qquad
   R_0^{-T} = \begin{pmatrix}
     1 & 1/2 & 1/3 & 1/4 \\
     0 & 1/2 & 1/2 & 9/20 \\
     0 &  0  & 1/6 & 1/4 \\
     0 &  0  &  0  & 1/20
    \end{pmatrix}.
\end{equation}
The entries of $R_{2k}^{-T}$ are nonnegative and have unit column sums
for all $k\ge0$.  Next we renormalize the shifted Legendre polynomials
and define
\begin{align}
  \dot{P}_n(x,y) &= h(x)^{-1/2}\sqrt{2n+1}\wtil{P}_n(y/h(x)), \\
  \notag
  \dot{Y}_{2k} = \left(\dot{P}_0,\ldots,\dot{P}_{2k+3}\right)^T
  &= h^{-1/2}D_{2k}
  \wtil{Y}_{2k}, \qquad
  D_{2k} = \diag\left(\sqrt{1},\sqrt{3},\ldots,\sqrt{4k+7}\right)
\end{align}
so that (\ref{eqn:A:dot}) becomes
\begin{align}
\label{eqn:A:dot2}\quad
  \frac{I_2}{I_3}\alpha^\e{2k}_{xx} &=
  \left(\dot{Y}_{2k}(x,y)\right)^T
  D_{2k}^{-1}R_{2k}^{-T}
  \left[V_0\dot{A}_0^\e{2k} + V_1\dot{A}_1^\e{2k} \right]
  \left(h(x)^{-1/2} \Phi_{2k+2}(x)\right), \\
  \notag
  \beta^\e{2k}_{xx} &= \left(\dot{Y}_{2k}(x,y)\right)^T
  D_{2k}^{-1}R_{2k}^{-T}
  \dot{B}^\e{2k}\left(h(x)^{-3/2}\Phi_{2k+2}(x)\right), \\
  \notag
  \frac{I_2}{I_3}h^2\alpha^\e{2k}_{xxxx} &=
  \left(\dot{Y}_{2k}(x,y)\right)^T
  D_{2k}^{-1}R_{2k}^{-T}
  \left[V_0\ddot{A}_0^\e{2k} + V_1\ddot{A}_1^\e{2k} \right]
  \left(h(x)^{-1/2} \Phi_{2k+4}(x)\right), \\
\notag
  h^2\beta^\e{2k}_{xxxx} &= \left(\dot{Y}_{2k}(x,y)\right)^T
  D_{2k}^{-1}R_{2k}^{-T}
  \ddot{B}^\e{2k}\left(h(x)^{-3/2}\Phi_{2k+4}(x)\right).
\end{align}

\begin{example} \upshape
From (\ref{eqn:Q0}) and (\ref{eqn:psi0}), we have
\begin{align}
  \notag
  Q^\e{0} = \,&\frac{V_0+V_1}{2}\frac{I_2}{I_3}, \\
  \psi^\e0_{xx} =
  \,&\left[
    \left(V_0+\frac{V_1}{2}\right)\left(\frac{-4h_x^2 + 2hh_{xx}}{h}\right)
    +Q^\e{0}\left(\frac{18h_x^2-6hh_{xx}}{h^2}\right)
    \right]\frac{y^2}{h^2} \\
  \notag
  &+ \left[
    (V_0+V_1)\left(\frac{6h_x^2 - 2hh_{xx}}{h}\right)
    +Q^\e{0}\left(\frac{-24h_x^2+6hh_{xx}}{h^2}\right)
  \right]\frac{y^3}{h^3},
\end{align}
which has the form described in (\ref{eqn:psi:form2}) and
(\ref{eqn:A:dot}) with
\begin{equation}
  \left[\dot{A}_0^\e0,\dot{A}_1^\e0,\dot{B}_0^\e0\right] =
  \left(\begin{array}{cc|cc|cc}
    \phm0 & \phm 0 & \phm0 & \phm 0 & \phm0 & \phm 0 \\
    \phm0 & \phm 0 & \phm0 & \phm 0 & \phm0 & \phm 0 \\
    -4 & \phm 4 & -2 & \phm 2 & \phm 18 & -12 \\
    \phm 6 & -4 & \phm 6 & -4 & -24 & \phm 12
  \end{array} \right)
\end{equation}
and the form described in (\ref{eqn:psi:form2}) and
(\ref{eqn:A:dot2}) with
\begin{equation*}
  D_0^{-1}R_0^{-T}
  \left[\dot{A}_0^\e0,\dot{A}_1^\e0,\dot{B}_0^\e0\right] =
  \left(\begin{array}{cc|cc|cc}
    \frac{1}{6} & \frac{2}{6} & \frac{5}{6} & \frac{-2}{6} & 0 &
    \frac{-6}{6} \\[3pt]
    \frac{7}{10\sqrt{3}} & \frac{2}{10\sqrt{3}} &
	\frac{17}{10\sqrt{3}} & \frac{-8}{10\sqrt{3}} &
	\frac{-18}{10\sqrt{3}} & \frac{-6}{10\sqrt{3}} \\[5pt]
    \frac{5}{6\sqrt{5}} & \frac{-2}{6\sqrt{5}} & \frac{7}{6\sqrt{5}} &
    \frac{-4}{6\sqrt{5}} & \frac{-18}{6\sqrt{5}} & \frac{6}{6\sqrt{5}} \\[5pt]
    \frac{3}{10\sqrt{7}} & \frac{-2}{10\sqrt{7}} & \frac{3}{10\sqrt{7}} &
    \frac{-2}{10\sqrt{7}} & \frac{-12}{10\sqrt{7}} &  \frac{6}{10\sqrt{7}}
  \end{array} \right),
\end{equation*}
where we recall that $\Phi_2(x)=(h_x^2, \frac{1}{2}h h_{xx})^T$.  The
matrices $\ddot{A}_0^\e{0}$, $\ddot{A}_1^\e{0}$, $\ddot{B}^\e{0}$
representing $\psi^\e{0}_{xxxx}$ are each $4\times5$ matrices while
$\Phi_4(x)$ was given
in (\ref{eqn:Phi:024}).
\end{example}

To compute $\|\psi^\e{2k}_{xx}\|_0$ and
$\|h^2\psi^\e{2k-2}_{xxxx}\|_0$, we note that each of the
expressions in (\ref{eqn:A:dot2}) is of the form $\sum_{n=0}^{2k+3}
\dot{P}_n(x,y)w_n(x)$, where $w(x)=Sz(x)$, $S$ is a constant matrix,
$z(x)=h^{-\frac{m}{2}} \Phi_{2k+2j}(x)$, $j=1$ or $2$, and $m=1$ or
$3$.  For fixed $x$, we have $\int_0^{h(x)}\dot{P}_m(x,y)
\dot{P}_n(x,y)\,dy=\delta_{mn}$.  It follows that
\begin{equation}
  \label{eqn:dotP:w}
  \int_0^1\int_0^{h(x)}\left(\sum_n
  \dot{P}_n(x,y)w_n(x)\right)^2\,dy\,dx = 
  \int_0^1 \vspace{6pt}\sum_n w_n(x)^2\,dx.
\end{equation}
Moreover, $\sum_n w_n(x)^2=\|w(x)\|_2^2 \le \|S\|_F^2 \|z(x)\|_2^2$,
where $\|\cdot\|_F$ and $\|\cdot\|_2$ are the Frobenius and 2-norms of
a matrix and vector, respectively.  Integrating, we have $\int_0^1
\|z(x)\|_2^2\,dx=I_m\|\wtil{E}_m^\e{2k+2j}\|_1$. Since
$\|\wtil{E}_m^\e{2\ell}\|_1\le d_{2\ell}
r_k^{-4\ell}$ for $0\le \ell\le k+1$, we define
\begin{equation}
\begin{aligned}
  K_i^\e{2k} &= \sqrt{d_{2k+2}}
  \left\| D_{2k}^{-1} R_{2k}^{-T}\dot{A}_i^\e{2k}\right\|_F,
  & \qquad i=0,1, \;\; k\ge0&, \\
  K_2^\e{2k} &= \sqrt{d_{2k+2}}
  \left\| D_{2k}^{-1} R_{2k}^{-T}\dot{B}^\e{2k}\right\|_F,
  & k\ge0&,\\
  \wtil{K}^\e{2k}_i &= \sqrt{d_{2k+4}}
  \left\| D_{2k}^{-1} R_{2k}^{-T}
  \ddot{A}_i^\e{2k}\right\|_F, & \quad i=0,1, \;\; k\ge0&, \\
  \wtil{K}^\e{2k}_2 &= \sqrt{d_{2k+4}}
  \left\| D_{2k}^{-1} R_{2k}^{-T} \ddot{B}^\e{2k}\right\|_F,
  & k\ge0&,
\end{aligned}
\end{equation}
so that
\begin{equation}
\begin{aligned}
  \left\|\frac{I_2}{I_3}\alpha^\e{2\ell}_{xx}\right\|_0 &\le
  \sqrt{I_1}\left(|V_0| K_0^\e{2\ell} + |V_1| K_1^\e{2\ell}\right)
  r_k^{-2\ell-2}, & \;\; &0\le \ell\le k, \\
  \left\|\beta^\e{2\ell}_{xx}\right\|_0 &\le
  \sqrt{I_3} \, K_2^\e{2\ell} r_k^{-2\ell-2}, & &0\le \ell\le k, \\
  \left\|\frac{I_2}{I_3}h^2\alpha^\e{2\ell}_{xxxx}\right\|_0 &\le
  \sqrt{I_1}\left(|V_0| \wtil{K}_0^\e{2\ell} +
  |V_1| \wtil{K}_1^\e{2\ell}\right) r_k^{-2\ell-4}, & &0\le \ell\le k-1, \\
  \left\|h^2\beta^\e{2\ell}_{xxxx}\right\|_0 &\le
  \sqrt{I_3} \, \wtil{K}_2^\e{2\ell} r_k^{-2\ell-4}, & &0\le \ell\le k-1.
\end{aligned}
\end{equation}
From the bound (\ref{eqn:Q:bound}) on $|Q^\e{2\ell}|$ and the formula
(\ref{eqn:psi:form2}) for $\psi^\e{2k}$ in terms of $\alpha^\e{2k}$
and $\beta^\e{2k-2\ell}$, we see that after increasing $K_0^\e{2k}$,
$K_1^\e{2k}$, $\wtil{K}_0^\e{2k}$, $\wtil{K}_1^\e{2k}$ via
\begin{equation}
\label{eqn:K:loop}
  \begin{aligned}
    K_i^\e{2k} &= K_i^\e{2k} + \sum_{\ell=0}^k
    \kappa_i^\e{2\ell} K_2^\e{2k-2\ell},
    \qquad i=0,1,\;\; k\ge0, \\
    \wtil{K}_i^\e{2k} &= \wtil{K}_i^\e{2k} +
    \sum_{\ell=0}^k \kappa_i^\e{2\ell}
    \wtil{K}_2^\e{2k-2\ell}, \qquad i=0,1,\;\; k\ge0,
  \end{aligned}
\end{equation}
we have
\begin{equation}
  \label{eqn:interior:effect}
  \begin{aligned}
    \left\|\psi^\e{2k}_{xx}\right\|_0 &\le
    \sqrt{I_1}\,\left(|V_0| K_0^\e{2k} + |V_1| K_1^\e{2k}\right) r_k^{-2k-2},
  & \;\; &k\ge0, \\
    \left\|h^2\psi^\e{2k-2}_{xxxx}\right\|_0 &\le
    \sqrt{I_1}\,\left(|V_0| \wtil{K}_0^\e{2k-2} +
    |V_1| \wtil{K}_1^\e{2k-2}\right) r_k^{-2k-2}, & \;\; &k\ge0.
  \end{aligned}
\end{equation}
In each term $Q^\e{2\ell}\beta^\e{2k-2\ell}$, we have used
$I_2^2/I_1I_3\le1$ (which follows from the Cauchy--Schwarz
inequality) to majorize $I_2/\sqrt{I_3}$ by $\sqrt{I_1}$.  The
constants $K_i^\e{2k}$, $\wtil{K}_i^\e{2k}$ do not depend on $h$ and
may be computed once and for all along with the constants
$\kappa_i^\e{2k}$ and the matrices $A_i^\e{2k}$ and $B^\e{2k}$;
see Tables~\ref{tbl:K} and~\ref{tbl:wtil:K}.

\begin{table}[t]
\begin{center}
\caption{$K_i^\e{2k}$ before and after loop {\rm (\ref{eqn:K:loop})}.} \label{tbl:K}
\footnotesize
\begin{tabular}{rccccc}
  \hline \\[-8pt]
 $k$ & $K_0^\e{2k}$ before & $K_1^\e{2k}$ before &
  $K_2^\e{2k}$ & $K_0^\e{2k}$ after & $K_1^\e{2k}$ after \\[1pt]
  \hline \\[-8pt]
 0 & $9.95\times10^{-01}$ & $2.17\times10^{+00}$ & $2.99\times10^{+00}$ & $2.49\times10^{+00}$ & $3.67\times10^{+00}$ \\
 1 & $2.33\times10^{+00}$ & $4.70\times10^{+00}$ & $7.99\times10^{+00}$ & $8.41\times10^{+00}$ & $1.16\times10^{+01}$ \\
 2 & $7.42\times10^{+00}$ & $1.58\times10^{+01}$ & $2.51\times10^{+01}$ & $3.14\times10^{+01}$ & $4.32\times10^{+01}$ \\
 3 & $4.29\times10^{+01}$ & $8.62\times10^{+01}$ & $1.19\times10^{+02}$ & $1.62\times10^{+02}$ & $2.21\times10^{+02}$ \\
 4 & $4.58\times10^{+02}$ & $8.71\times10^{+02}$ & $1.03\times10^{+03}$ & $1.34\times10^{+03}$ & $1.87\times10^{+03}$ \\
 5 & $7.21\times10^{+03}$ & $1.52\times10^{+04}$ & $1.62\times10^{+04}$ & $1.85\times10^{+04}$ & $2.77\times10^{+04}$ \\
 6 & $1.87\times10^{+05}$ & $3.54\times10^{+05}$ & $3.51\times10^{+05}$ & $4.06\times10^{+05}$ & $5.90\times10^{+05}$ \\
 7 & $6.57\times10^{+06}$ & $1.25\times10^{+07}$ & $1.08\times10^{+07}$ & $1.28\times10^{+07}$ & $1.92\times10^{+07}$ \\
 8 & $2.74\times10^{+08}$ & $6.17\times10^{+08}$ & $4.64\times10^{+08}$ & $5.31\times10^{+08}$ & $8.87\times10^{+08}$ \\
 9 & $1.75\times10^{+10}$ & $3.28\times10^{+10}$ & $2.52\times10^{+10}$ & $3.12\times10^{+10}$ & $4.70\times10^{+10}$ \\
10 & $1.32\times10^{+12}$ & $2.40\times10^{+12}$ & $1.69\times10^{+12}$ & $2.22\times10^{+12}$ & $3.34\times10^{+12}$ \\
\hline
\end{tabular}
\end{center}
\end{table}

\begin{table}[t]
\begin{center}
\caption{$\wtil{K}_i^\e{2k}$ before and after loop {\rm (\ref{eqn:K:loop})}.} \label{tbl:wtil:K}
\footnotesize
\begin{tabular}{rccccc}
  \hline \\[-8pt]
 $k$ & $\wtil{K}_0^\e{2k}$ before & $\wtil{K}_1^\e{2k}$ before &
  $\wtil{K}_2^\e{2k}$ & $\wtil{K}_0^\e{2k}$ after & $\wtil{K}_1^\e{2k}$ after
  \\[1pt]
  \hline \\[-8pt]
 0 & $1.23\times10^{+02}$ & $2.08\times10^{+02}$ & $4.84\times10^{+02}$ & $3.65\times10^{+02}$ & $4.50\times10^{+02}$ \\
 1 & $7.81\times10^{+02}$ & $1.50\times10^{+03}$ & $3.35\times10^{+03}$ & $2.79\times10^{+03}$ & $3.65\times10^{+03}$ \\
 2 & $3.55\times10^{+03}$ & $7.74\times10^{+03}$ & $1.56\times10^{+04}$ & $1.47\times10^{+04}$ & $2.00\times10^{+04}$ \\
 3 & $1.73\times10^{+04}$ & $3.88\times10^{+04}$ & $6.63\times10^{+04}$ & $7.22\times10^{+04}$ & $1.00\times10^{+05}$ \\
 4 & $1.85\times10^{+05}$ & $3.76\times10^{+05}$ & $4.38\times10^{+05}$ & $5.37\times10^{+05}$ & $7.68\times10^{+05}$ \\
 5 & $3.97\times10^{+06}$ & $8.59\times10^{+06}$ & $8.12\times10^{+06}$ & $9.04\times10^{+06}$ & $1.40\times10^{+07}$ \\
 6 & $1.45\times10^{+08}$ & $2.66\times10^{+08}$ & $2.46\times10^{+08}$ & $2.81\times10^{+08}$ & $4.08\times10^{+08}$ \\
 7 & $6.62\times10^{+09}$ & $1.26\times10^{+10}$ & $1.00\times10^{+10}$ & $1.19\times10^{+10}$ & $1.81\times10^{+10}$ \\
 8 & $3.50\times10^{+11}$ & $8.00\times10^{+11}$ & $5.50\times10^{+11}$ & $6.36\times10^{+11}$ & $1.09\times10^{+12}$ \\
 9 & $2.81\times10^{+13}$ & $5.24\times10^{+13}$ & $3.81\times10^{+13}$ & $4.76\times10^{+13}$ & $7.22\times10^{+13}$ \\
10 & $2.59\times10^{+15}$ & $4.70\times10^{+15}$ & $3.12\times10^{+15}$ & $4.18\times10^{+15}$ & $6.31\times10^{+15}$ \\
\hline
\end{tabular}
\end{center}
\end{table}

Finally, we combine the boundary estimate (\ref{eqn:bc:effect}) with
the interior estimate (\ref{eqn:interior:effect}) to bound the
truncation error via (\ref{eqn:e2k0:bound2}).  In terms of $r_k$,
(\ref{eqn:bc:effect}) gives
\begin{equation}
  \label{eqn:bc:effect2}
  \left\|\frac{\gamma_{k}}{h^{1/2}}\right\|_{1/2,\veps} \le
  |V_1|\sqrt{I_1}\left[\frac{1}{2} + \frac{\veps}{r_k}
    \sqrt{\frac{I_3}{I_1}}\sqrt{\frac{85}{16}+\frac{20}{3}k}
      \right]r_k^{-2k-2}, \qquad k\ge0.
\end{equation}\vspace*{-24pt}\pagebreak

\noindent We now define\vspace{5pt}
\begin{equation}
  \label{eqn:rhok:def}
  \rho_k = \begin{cases}
    \left[\max\left(5K_0^\e{2k},\,\,
    5K_1^\e{2k}+\frac{15}{2}\right)
    \right]^{-\frac{1}{2k+2}}, & k = 0, \\
    \left[\max\left(5K_0^\e{2k}+\wtil{K}_0^\e{2k-2},\,\,
    5K_1^\e{2k}+\wtil{K}_1^\e{2k-2}+\frac{15}{2}, \,\,
    \rho_{k-1}^{-2k}\right)
    \right]^{-\frac{1}{2k+2}}, & k \ge 1,
    \end{cases}
\end{equation}
so that (\ref{eqn:e2k0:bound3}), (\ref{eqn:interior:effect}), and
(\ref{eqn:bc:effect2}) imply
\begin{equation}
\begin{aligned}
  \label{eqn:thm:bound0}
  &\left\|\psi_\tem{err}^\e{2k}\right\|_{2,\veps} +
  \veps^{2k+2}\left\|\psi^\e{2k}_{xx}\right\|_0 \\
  &\qquad\;\;\le \sqrt{I_1}\left( |V_0| + |V_1|\right)
  \left[\rho_k^{-2k-2} + 15\frac{\veps}{r_k}\sqrt{\frac{I_3}{I_1}}
    \sqrt{\frac{85}{16}+\frac{20}{3}k}\right]
  \left( \frac{\veps}{r_{k}} \right)^{2k+2}.
\end{aligned}
\end{equation}
To simplify this expression, we define
\begin{equation}
  \label{eqn:thetak:def}
  \theta_{k} = 15\rho_{k}^{2k+2}\sqrt{
    \frac{85}{16}+\frac{20}{3}k}
\end{equation}
and summarize the main result of this section as a theorem.\pagebreak

\begin{theorem} \label{thm:bound}
Suppose $k\ge0$, $h\in C^{2k+1,1}(T)$, $0< h(x)\le1$ for $0\le x\le
1$, and $\veps\le r_0/3$.
Then the truncation
errors $\psi_\tem{err}^\e{2k}$ and $Q_\tem{err}^\e{2k}$ in
{\rm(\ref{eqn:ek:def})} satisfy the bound
\begin{equation}
  \label{eqn:thm:bound}
  \begin{aligned}
  \sqrt{3}\left|Q_\tem{err}^\e{2k}\right|\le
  \left\|\psi_\tem{err}^\e{2k}\right\|_{2,\veps} & \le
  \left\|\psi_\tem{err}^\e{2k}\right\|_{2,\veps} +
  \veps^{2k+2}\left\|\psi^\e{2k}_{xx}\right\|_0 \\
  & \le \sqrt{I_1}\left( |V_0| + |V_1|\right)
  \left[1 + \theta_{k} \frac{\veps}{r_{k}}\sqrt{\frac{I_3}{I_1}}
  \right] \left( \frac{\veps}{\rho_{k} r_{k}} \right)^{2k+2},
\end{aligned}
\end{equation}
where $I_m$, $r_k$, $\rho_k$, and $\theta_k$ were defined in
{\rm(\ref{eqn:Q}), (\ref{eqn:r:def}), (\ref{eqn:rhok:def}),} and
{\rm(\ref{eqn:thetak:def}).}
\end{theorem}

\begin{remark} \upshape
The constants in this estimate have been organized to be either (1)
given in the problem statement or easily computable from $h$; or (2)
difficult to compute but universal (independent of $h$).  The first 26
constants in the latter category ($\rho_k$ and $\theta_k$) are given
in Table~\ref{tbl:rho:theta}.  We have therefore identified the
features of $h$ that are most likely to affect the validity of the
lubrication approximation.  In particular, higher derivatives are
allowed to be large in regions where $h$ is small (since $r_k$
depends on the uniform norms of the products
$\frac{1}{\ell!}h^{\ell-1}\partial_x^\ell h$ rather than on
$\partial_x^\ell h$ alone).
\end{remark}

\begin{figure}
\centerline{\includegraphics[width=.58\linewidth]
{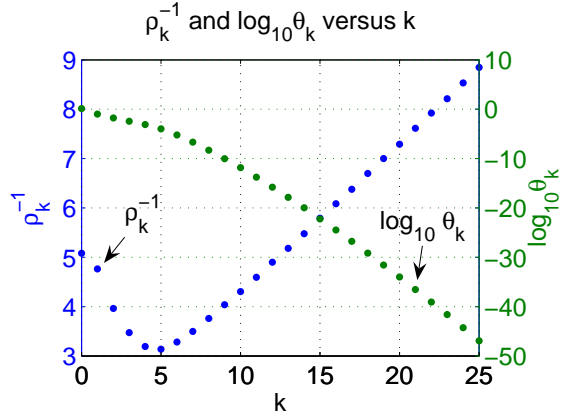}}
\caption{Plot of $\rho_k^{-1}$ and $\log_{10}\theta_k$ versus $k$.
Note that $\rho_k^{-1}$ initially decreases but eventually grows
almost linearly, indicating that the lubrication expansion is probably
an asymptotic series rather than a convergent series.  The term
involving $\theta_k$ in {\rm (\ref{eqn:thm:bound})} is only important when
$k$ is small, since $\theta_k$ converges rapidly to zero as\,
$k\rightarrow\infty$.} \label{fig:rho:theta}
\end{figure}

\begin{remark} \upshape
\label{rk:rho:increasing}
The term $\rho_{k-1}^{-2k}$ in the definition of $\rho_k$ ensures
that $\rho_k^{2k+2}$ is a nonincreasing function of $k$.  This
assumption is useful in the following section when deriving a bound on
the truncation error of the pressure.  Note that $\rho_k$ itself is
allowed to increase as long as $\rho_k^{2k+2}$ does not.  It is
probably not necessary to include this term in the definition of
$\rho_k$ since it is not the argmax for $1\le k\le25$, and by that
point $\theta_k$ (and hence $\rho_k^{2k+2}$) appears to be decreasing
rapidly without it; see Figure~\ref{fig:rho:theta}.
\end{remark}

\begin{table}[t]
\begin{center}
\caption{$\rho_k$ and $\theta_k$.} \label{tbl:rho:theta}
\footnotesize
\begin{tabular}{rcc}
  \hline
 $k$ & $\rho_k$ & $\theta_k$  \\
  \hline \\[-8pt]
 0 & 0.197 & $1.34\times10^{+00}$ \\
 1 & 0.210 & $1.01\times10^{-01}$ \\
 2 & 0.252 & $1.67\times10^{-02}$ \\
 3 & 0.288 & $3.58\times10^{-03}$ \\
 4 & 0.313 & $7.73\times10^{-04}$ \\
 5 & 0.319 & $1.03\times10^{-04}$ \\
 6 & 0.305 & $5.96\times10^{-06}$ \\
 7 & 0.286 & $2.15\times10^{-07}$ \\
 8 & 0.266 & $5.10\times10^{-09}$ \\
 9 & 0.248 & $9.15\times10^{-11}$ \\
10 & 0.232 & $1.43\times10^{-12}$ \\
11 & 0.218 & $1.69\times10^{-14}$ \\
12 & 0.204 & $1.58\times10^{-16}$ \\
13 & 0.193 & $1.42\times10^{-18}$ \\
14 & 0.183 & $1.04\times10^{-20}$ \\
15 & 0.173 & $5.98\times10^{-23}$ \\
16 & 0.164 & $3.46\times10^{-25}$ \\
17 & 0.157 & $1.75\times10^{-27}$ \\
18 & 0.149 & $6.86\times10^{-30}$ \\
19 & 0.143 & $2.72\times10^{-32}$ \\
20 & 0.137 & $1.02\times10^{-34}$ \\
21 & 0.131 & $2.94\times10^{-37}$ \\
22 & 0.126 & $8.36\times10^{-40}$ \\
23 & 0.122 & $2.41\times10^{-42}$ \\
24 & 0.117 & $5.40\times10^{-45}$ \\
25 & 0.113 & $1.15\times10^{-47}$ \\
\hline
\end{tabular}
\end{center}
\end{table}


\subsection{Velocity, vorticity, and pressure}
\label{sec:pressure}

We now show how to use the error bound we have obtained for the stream
function to bound the error in the velocity, vorticity, and pressure.
We define $u^\e{2k}$, $v^\e{2k}$, $\omega^\e{2k}$, and $p^\e{2k}$ in
terms of $\psi^\e{2k}$ as in (\ref{eqn:uvp:from:psi}) and define,
e.g.,
\begin{align}
  \omega_\tem{err}^\e{2k} &= \omega_\tem{exact} - \omega_\tem{approx}^\e{2k},
  &\quad \omega_\tem{approx}^\e{2k} &=
  \omega^\e{0}+\veps^2 \omega^\e{2} + \cdots + \veps^{2k}\omega^\e{2k}.
\label{eqn:uvwp:err:def}
\end{align}
From (\ref{eqn:uvp:from:psi}), we then have
\begin{align}
  \label{eqn:u:err:from:psi}
  \left(u_\tem{err}^\e{2k}, v_\tem{err}^\e{2k}\right)
  &= \nabla\times\psi_\tem{err}^\e{2k}, &
  \omega_\tem{err}^\e{2k} &= -\Delta_\veps \psi_\tem{err}^\e{2k}
  - \veps^{2k+2}\psi^\e{2k}_{xx}, \\
  \label{eqn:p:err:from:psi}
  \partial_x p_\tem{err}^\e{2k} &= -\partial_y\omega_\tem{err}^\e{2k}, &\;\;
  \partial_y p_\tem{err}^\e{2k} &= \begin{cases}
    \partial_x\left(\veps^2\omega_\tem{exact}\right), & k=0, \\[1pt]
  \partial_x\left(\veps^2\omega_\tem{err}^\e{2k-2}\right), & k\ge1, \end{cases}
\end{align}
which immediately gives bounds on the error in velocity and vorticity:\vspace{5pt}
\begin{equation}
  \begin{aligned}
    \left(\left\|u_\tem{err}^\e{2k}\right\|_{1,\veps}^2 +
  \left\|\veps v_\tem{err}^\e{2k}\right\|_{1,\veps}^2\right)^{1/2} &\le
  \left\|\psi_\tem{err}^\e{2k}\right\|_{2,\veps}\hspace*{75.2pt}\le (*), \qquad k\ge0, \\
  \left\|\omega_\tem{err}^\e{2k}\right\|_0 & \le
  \left\|\psi_\tem{err}^\e{2k}\right\|_{2,\veps} +
  \veps^{2k+2}\left\|\psi^\e{2k}_{xx}\right\|_0 \le 
  (*), \qquad k\ge0,
  \end{aligned}
\end{equation}
where $(*)$ is the right-hand side of (\ref{eqn:thm:bound}).
Obtaining a bound on the error in the pressure is somewhat more
difficult, as it relies on the fact that the gradient is an isomorphism
from $L^2_\#(\Omega)$ (the space of square integrable functions with
zero mean) onto the polar set
\begin{equation}
  V^0=\left\{\mb{f}\in H^{-1}(\Omega)^2\;:\;\langle \mb{f},\mb{u}\rangle=0
 \text{ whenever } \mb{u}\in V\right\},
\end{equation}
where $V=\{\mb{u}\in H^1_0(\Omega)^2\;:\;\nabla\cdot\mb{u}=0\}$.
Given $\mb{f}\in V^0$, there is a unique $p\in L^2_\#(\Omega)$
such that $\nabla p=\mb{f}$; moreover, $p$ satisfies
\begin{equation}
  \label{eqn:beta:def1}
  \|p\|_0 \le \beta^{-1}|\mb{f}|_{-1}, \qquad
  \beta = \inf_p \sup_\mb{u} \frac{|(p,\nabla\cdot\mb{u})|}{
    \|p\|_0\,|\mb{u}|_1}.
\end{equation}
Here we use a standard (unweighted) Sobolev norm for $\mb{f}$.
More precisely, as $\|\cdot\|_1$ and $|\cdot|_1$ are equivalent on
$H^1_0(\Omega)$, the negative norms
\begin{equation}
  |\mb{f}|_{-1} = 
  \sup_{|u|_{1}^2 + |v|_{1}^2=1}
  |\langle\mb{f},(u;v)\rangle|,\qquad
  \|\mb{f}\|_{-1} = 
  \sup_{\|u\|_{1}^2 + \|v\|_{1}^2=1}
  |\langle\mb{f},(u;v)\rangle|
\end{equation}
are equivalent on $H^{-1}(\Omega)^2$.  Since $|u|_{1}\le\|u\|_{1}$ for
all $u\in H^1_0(\Omega)$, we have $\|\mb{f}\|_{-1}\le |\mb{f}|_{-1}$
for all $\mb{f}\in H^{-1}(\Omega)^2$.

Explicit estimates \cite{chiz00,stoyan01,dobro03}
for the LBB constant $\beta$ in (\ref{eqn:beta:def1}) 
have been obtained for rectangular domains (with no periodicity),
e.g.,
\begin{equation}\label{eqn:chiz:beta}
  \frac{1}{\ell}\sin\frac{\pi}{8}\le \beta \le \frac{\pi}{2\sqrt{3}\ell},
  \qquad \ell = \max\left(\frac{L_1}{L_2},\frac{L_2}{L_1}\right),
  \qquad \Omega = (0,L_1)\times(0,L_2).
\end{equation}
The lower bound here also works for an $x$-periodic rectangle as the
condition that $\mb{u}\vert_{x=0}=0$ and $\mb{u}\vert_{x=L_1}=0$ is
more restrictive than $\mb{u}\vert_{x=0}=\mb{u}\vert_{x=L_1}$.
Explicit estimates are also known for domains that are star-shaped
with respect to each point in a ball of radius $R$ contained inside
$\Omega$; see \cite{galdi,stoyan01}.  Our interest in the present work
is in $x$-periodic domains with the upper boundary given by a function
$h(x)$.  Such domains are not in general star-shaped, so the
previously known results do not apply.  In \cite{infsup}, we improve
the estimate for the lower bound on $\beta$ in (\ref{eqn:chiz:beta})
for an $x$-periodic rectangle by a factor of about 3.5 and show how to
avoid invoking Rellich's theorem in the change of variables to the
case that $\Omega$ is $x$-periodic with the upper boundary given by
$h(x)$.  It is shown that
\begin{equation*}
  \beta^{-1}\le \frac{9}{4}\left(1+r_0^{-2}\right)\left(\frac{h_1}{h_0}\right)^{1/2}
  \max\left(4,\frac{L}{h_0},\frac{h_1}{h_0}\right), \qquad
  \begin{array}{l}
    h_0=\min_{0\le x\le L}h(x), \\ h_1=\max_{0\le x\le L}h(x),
  \end{array}
\end{equation*}
where $r_0=\max(\|h_x\|_\infty,
  \|\frac{1}{2}hh_{xx}\|_\infty^{1/2})^{-1}$ and $L$ is the
period of $h(x)$.  In the current case, the length scales $\bar{H}$
and $\bar{W}$ were chosen so that $L=1$ and $h_1\le1$ in the
dimensionless problem.  Thus, solving $\nabla p=\mb{f}$ yields
\begin{equation}
  \label{eqn:beta:def}
  \|p\|_0 \le \beta^{-1}|\mb{f}|_{-1}, \qquad
  \beta^{-1} \le \max\left(9h_0^{-1/2},\frac{9}{4}h_0^{-3/2}\right)
  \left(1+r_0^{-2}\right).
\end{equation}
The dependence on gap thickness $h_0$ occurs because $p$ can change
rapidly in the gap without a large penalty from $\mb{u}$.  In
\cite{infsup}, an example is given to show that the factor of
$h_0^{-3/2}$ in the formula (\ref{eqn:beta:def}) for $\beta^{-1}$
cannot be improved.

We have reduced the problem of bounding $p_\tem{err}^\e{2k}$ to
that of bounding the functional on the right-hand side of $\nabla
p_\tem{err}^\e{2k}=\mb{f}_k$ in (\ref{eqn:p:err:from:psi}), namely,
\begin{equation}
  \label{eqn:Fk:def}
  \langle \mb{f}_k, (u;v) \rangle =
  \begin{cases}
    \int_\Omega \omega_\tem{err}^\e0\, u_y -
    \veps^2\omega_\tem{exact}\, v_x\,dA, & k=0, \\[1pt]
    \int_\Omega \omega_\tem{err}^\e{2k} u_y -
    \veps^2\omega_\tem{err}^\e{2k-2} v_x\,dA, & k\ge1.
  \end{cases}
\end{equation}
First, we check that $\mb{f}_k$ belongs to $V^0$.  If $u,v\in
H^1_0(\Omega)$ and $u_x+v_y=0$, the function $\phi(x,y)=\int_0^y
u(x,\eta)\,d\eta$ satisfies $\phi_y=u$, $\phi_x=-v$ and so belongs to
$\Psi$.  As a result,
\begin{equation}
  \begin{aligned}
  \langle \mb{f}_0, (u;v) \rangle &=
  \int_\Omega -\omega^\e{0}\phi_{yy} + \omega_\tem{exact}\left(\phi_{yy}+\veps^2
  \phi_{xx}\right)\,dA \\
  &= \int_\Omega \psi^\e{0}_{yy}\phi_{yy} +
  (\Delta_\veps\psi_\tem{exact})(\Delta_\veps\phi)\,dA=0,
 \end{aligned}
\end{equation}
and for $k\ge1$,
\begin{equation*}
\begin{aligned}
  \langle \mb{f}_k, (u;v) \rangle &=
  \overbrace{\langle \mb{f}_{k-1}, (u;v) \rangle}^0 +
  \int_\Omega -\veps^{2k}\omega^\e{2k}\phi_{yy} -
  \veps^{2k} \omega^\e{2k-2}\phi_{xx}\,dA \\
  &= \veps^{2k}\int_\Omega \psi_{yy}^\e{2k}\phi_{yy} +
  \psi^\e{2k-2}_{xx}\phi_{yy} + \psi^\e{2k-2}_{yy}\phi_{xx} +
  \underbrace{\psi^\e{2k-4}_{xx}\phi_{xx}}_\text{omit if $k=1$}\,dA = 0.
 \end{aligned}
\end{equation*}
Next, we bound the norm of $\mb{f}_k$.
From (\ref{eqn:Fk:def}), we see that
$$
  \notag
  |\mb{f}_0|_{-1}\le \left\|\omega_\tem{err}^\e{0}\right\|_0 +
      \left\|\veps^2\omega_\tem{exact}\right\|_0, \qquad
  |\mb{f}_k|_{-1}\le \left\|\omega_\tem{err}^\e{2k}\right\|_0 +
      \left\|\veps^2\omega_\tem{err}^\e{2k-2}\right\|_0, \quad k\ge1.
$$
Denoting the right-hand side of (\ref{eqn:thm:bound}) by
$(*)$, we claim that
\begin{equation}
  \label{eqn:vort:bound0}
  \frac{\veps^2}{r_0^2}\|\omega_\tem{exact}\|_0 \le (*), \quad (k=0), \qquad
  \frac{\veps^2}{r_k^2}\left\|\omega_\tem{err}^\e{2k-2}\right\|_0 \le
  (*), \quad (k\ge1).
\end{equation}
Once this is shown to be true, we will have the following
bound on $|\mb{f}_k|_{-1}$:
\begin{equation}
  \label{eqn:Fk:bound}
  |\mb{f}_k|_{-1} \le \left(1+r_k^2\right)(*),  \qquad (k\ge0).
\end{equation}
Together with (\ref{eqn:beta:def}) and the fact that
$r_k\le r_0$ for $k\ge0$, this will give
\begin{equation}
  \left\|p_\tem{err}^\e{2k}\right\|_0\le
  \max\left(9h_0^{-1/2},\frac{9}{4}h_0^{-3/2}\right)
  \left(r_k+r_k^{-1}\right)^2(*).
\end{equation}
Note that there is at least a power of $r_k^{-2}$ in $(*)$ to prevent
this bound from diverging as $r_k\rightarrow\infty$.  It may be
possible to improve the bound in this regime by replacing
$\|\veps^2\omega_\tem{exact}\|_0$ in (\ref{eqn:vort:bound0}) by
$\|\veps^2\partial_x\omega_\tem{exact}\|_{-1}$, but this seems very
difficult.  At any rate, if $r_k=\infty$, then $(*)=0$, $h(x)$ is a
constant function, the exact and approximate vorticity are constants,
$\mb{f}_k$ is the zero functional, and
$p_\tem{exact}=p_\tem{err}^\e{2k}=0$.  Let us now prove
(\ref{eqn:vort:bound0}).  For $k\ge1$, this follows from
\begin{equation}
  \begin{aligned}
  &\left[\rho_{k-1}^{-2k} + 15\frac{\veps}{r_{k-1}}
    \sqrt{\frac{I_3}{I_1}}\sqrt{\frac{85}{16} + \frac{20}{3}(k-1)}
    \right]\left(\frac{\veps}{r_{k-1}}\right)^{2k}
  \left(\frac{\veps}{r_{k}}\right)^{2} \\ &\hspace*{1.25in} \le
  \left[\rho_{k}^{-2k-2} + 15\frac{\veps}{r_{k}}
    \sqrt{\frac{I_3}{I_1}}\sqrt{\frac{85}{16} + \frac{20}{3}k}
    \right]\left(\frac{\veps}{r_{k}}\right)^{2k+2},
\end{aligned}
\end{equation}
which holds because $\rho_k^{2k+2}$ and $r_k$ are nonincreasing
functions of $k$; see Remark~\ref{rk:rho:increasing} and the
definition of $r_k$ in (\ref{eqn:r:def}).  For $k=0$, we use
Theorem~\ref{thm:psi:bound} to conclude that
\begin{equation}
  \|\omega_\tem{exact}\|_0 \le \|\psi\|_{2,\veps} \le 15
  \left(\left\|h^{-1/2}g_0\right\|_{1/2,\veps} +
  \left\|h^{-1/2}g_1\right\|_{1/2,\veps}\right),
\end{equation}
where $g_0(x)=V_0$ and $g_1(x)=(1+\veps^2h_x(x)^2)^{-1/2}$.  Now
\begin{equation}
  \begin{aligned}
  \left\| h^{-1/2}g_0 \right\|_{1/2,\veps}^2 &\le
  V_0^2\left\|h^{-1/2}\right\|_{1,\veps}^2 =
  V_0^2\int_0^1 h^{-1}+\frac{1}{4}(\veps h_x)^2 h^{-3}\,dx \\
  &\le V_0^2 \left[ I_1 + \frac{1}{4} \left(\frac{\veps}{r_0}\right)^2 I_3
    \right]
  \end{aligned}
\end{equation}
and
\begin{align}
  &\left\| h^{-1/2}g_1 \right\|_{1/2,\veps}^2 \le
  V_1^2\left\|h^{-1/2}\left(1+\veps^2h_x\right)^{-1/2}\right\|_{1,\veps}^2 \\
  \notag
  &\qquad = V_1^2\int_0^1 \left(h^{-1}\right)(\cdot)^{-1} + \left[-\frac{1}{2}h^{-3/2}
    (\veps h_x)(\cdot)^{-1/2} - \left(h^{-1/2}\right)(\cdot)^{-3/2}\veps^3 h_x h_{xx}
    \right]^2 dx \\
  \notag
  &\qquad \le V_1^2 \int_0^1 h^{-1} + \frac{1}{2} h^{-3}(\veps h_x)^2
  + 8 h^{-3}(\veps h_x)^2 \left(\veps^2 \frac{hh_{xx}}{2}\right)^2 \,dx \\
  &\qquad\le
  V_1^2\left[ I_1 + \left(\frac{1}{2} + 8\frac{\veps^4}{r_0^4}\right)
    \left(\frac{\veps}{r_0}\right)^2 I_3\right].
\end{align}
Since we have assumed that $\veps\le r_0/3$, we conclude that
\begin{equation}
  \frac{\veps^2}{r_0^2}\|\omega_\tem{exact}\|_0
  \le\sqrt{I_1}(|V_0|+|V_1|)\left[
    15 + \left(15\sqrt{\frac{1}{2} + \frac{8}{81}}\right)
    \frac{\veps}{r_0} \sqrt{\frac{I_3}{I_1}}
    \right]\frac{\veps^2}{r_0^2}.
\end{equation}
Comparing this to (\ref{eqn:thm:bound0}) with $k=0$ and noting from
Table~\ref{tbl:rho:theta} that $\rho_0^{-2}\ge15$, we obtain
(\ref{eqn:vort:bound0}) as claimed.  Thus, we have proved the
following theorem. 

\begin{theorem} \label{thm:bound2}
Suppose $k\ge0$, $h\in C^{2k+1,1}(T)$, $0< h_0\le h(x)\le1$ for $x\in
T$, and $\veps\le r_0/3$.
Then the truncation errors of the stream function, flux, velocity,
vorticity, and pressure satisfy the bounds
\begin{align}
\notag
    \left\|\psi_\tem{err}^\e{2k}\right\|_{2,\veps} &\le (*), \hspace*{18pt}
    |Q_\tem{err}^\e{2k}|\le\frac{(*)}{\sqrt{3}}, \hspace*{18pt}
    \left(\left\|u_\tem{err}^\e{2k}\right\|_{1,\veps}^2 +
    \left\|\veps v_\tem{err}^\e{2k}\right\|_{1,\veps}^2\right)^{1/2} \le (*),\hspace*{-12pt} \\
  \label{eqn:thm:bound2}
    \left\|\omega_\tem{err}^\e{2k}\right\|_0 &\le (*), \hspace*{18pt}
  \left\|p_\tem{err}^\e{2k}\right\|_0 \le
  \max\left(9h_0^{-1/2},\frac{9}{4}h_0^{-3/2}\right)
  \left(r_k+r_k^{-1}\right)^2 (*),
\end{align}
where $T=[0,1]_p$ is the periodic unit interval
\begin{gather}
\label{eqn:star:def2}
  (*) = \sqrt{I_1}\left( |V_0| + |V_1|\right)
  \left[1 + \theta_{k} \frac{\veps}{r_{k}}\sqrt{\frac{I_3}{I_1}}
  \right] \left( \frac{\veps}{\rho_{k} r_{k}} \right)^{2k+2}, \\
  r_k = \left(\max_{1\le \ell\le 2k+2}
  \left\{\left\|\frac{1}{\ell!}h^{\ell-1}
  \partial_x^\ell h\right\|_\infty^{1/\ell}\right\}\right)^{-1}, \qquad
  I_m = \int_0^1 h(x)^{-m}\,dx,
\end{gather}
and $\rho_k$, $\theta_k$ are constants independent of $h$ that can
be computed once and for all as described in section~{\rm\ref{sec:error:est}}
and listed in Table~{\rm\ref{tbl:rho:theta}.}
\end{theorem}

\section{Finite element validation}
\label{sec:fe}

In this section, to test the error bounds of
Theorem~\ref{thm:bound2}, we compute
$u_\tem{err}^\e{2k}$, $v_\tem{err}^\e{2k}$, $p_\tem{err}^\e{2k}$,
$\omega_\tem{err}^\e{2k}$ numerically for the simple geometry
described by
\begin{equation}
  \label{eqn:h:sine}
  h(x) = \frac{1+a}{2} + \frac{1-a}{2}\sin(2\pi x),
  \qquad \begin{aligned}
    &\text{Case 1: }\; a=1/5, \\
    &\text{Case 2: }\; a=1/100,
  \end{aligned}
\end{equation}
with boundary conditions $V_0 = -0.5$, $V_1 = 1$.  We do this by
comparing $u_\tem{approx}^\e{2k}$, $v_\tem{approx}^\e{2k}$,
$p_\tem{approx}^\e{2k}$, $\omega_\tem{approx}^\e{2k}$ in
(\ref{eqn:uvwp:err:def}) to finite element solutions of the Stokes
equations on appropriately rescaled geometries.

\begin{figure}[htbp]
\begin{center}
\includegraphics[width=25pc]
{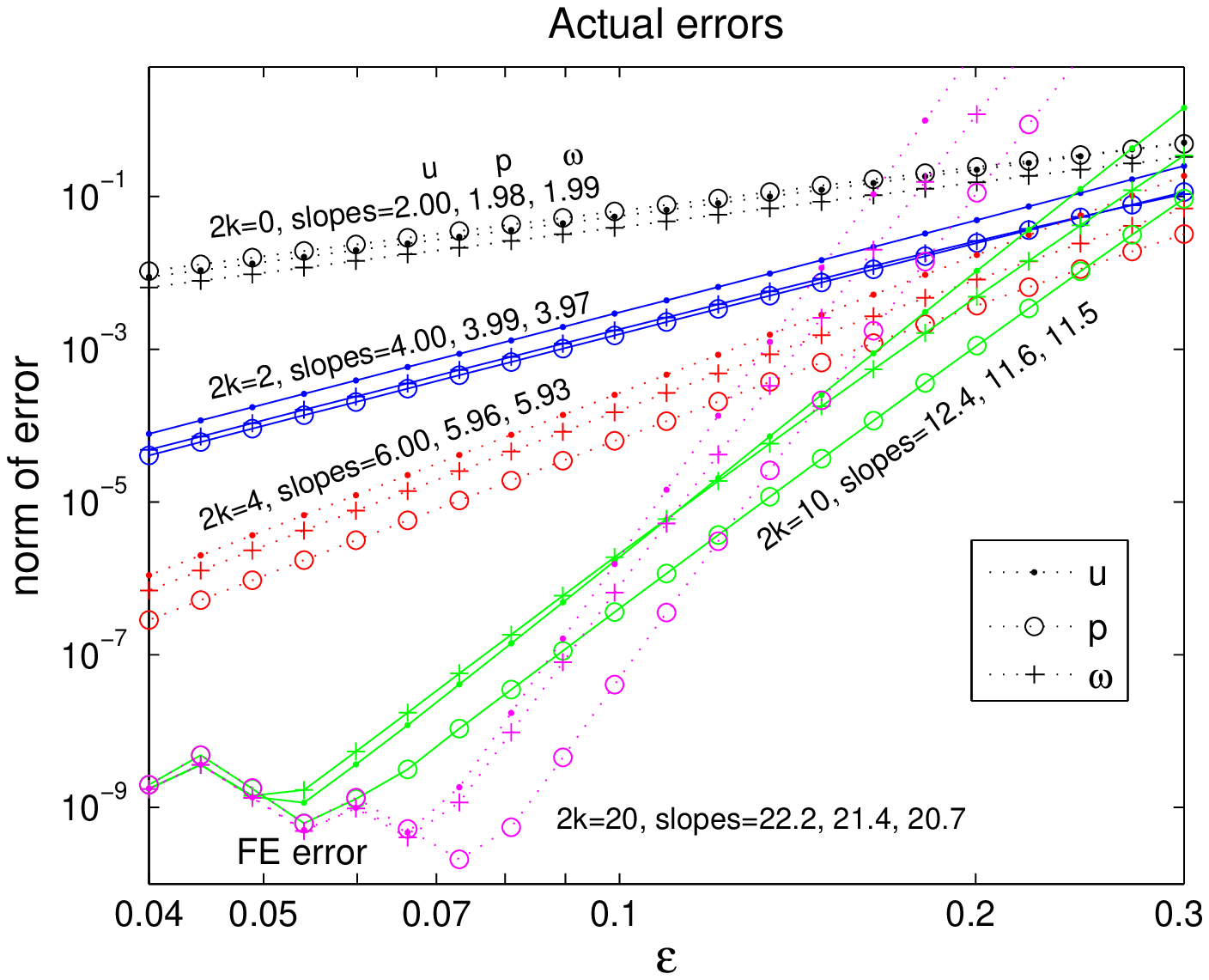} \\[8pt]
\includegraphics[width=25pc]
{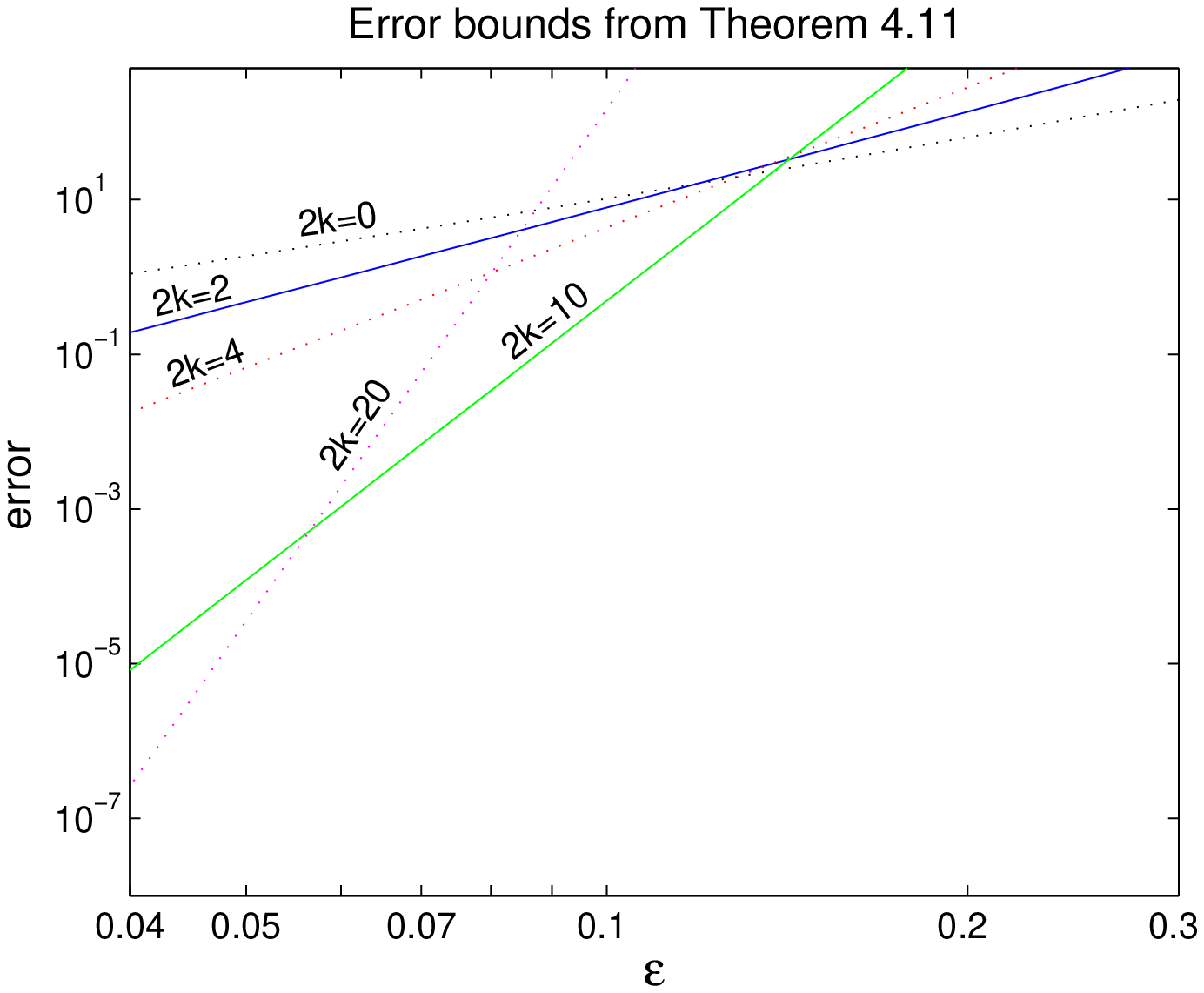}
\end{center}
\caption{Top: plot of
$\|\mb{u}_\tem{err}^\e{2k}\|_{1,\veps}=
(\|u_\tem{err}^\e{2k}\|_{1,\veps}^2
+\|\veps v_\tem{err}^\e{2k}\|_{1,\veps}^2)^{1/2}$,
$\|p_\tem{err}^\e{2k}\|_0$, and
$\|\omega_\tem{err}^\e{2k}\|_0$ for $a=1/5$,
$0.04\le\veps\le0.3$, and $2k\in\{0,2,4,10,20\}$.
The slopes of the lines were computed via
linear regression using the 
smallest $10$ values of $\veps$ for which the finite element solution
is trusted ($\veps\ge0.066$ for $2k=10$
and $\veps\ge0.09$ for $2k=20$).
As expected, for fixed $k$, the
error is $O(\veps^{2k+2})$.
Bottom: plot of the error bound $(*)$ in
{\rm (\ref{eqn:star:def2})}, using $V_0=-.5$, $V_1=1$, $I_1=2.236$,
$I_3=24.60$, and $r_k=0.3559$ for $k\ge0$, as appropriate for $h(x)$
in {\rm (\ref{eqn:h:sine})} with $a=1/5$.  }\label{fig:errorsAB}
\vspace{-12pt}
\end{figure}

\begin{figure}[p]
\begin{center}
\mypsdraft
\includegraphics[width=27pc]
{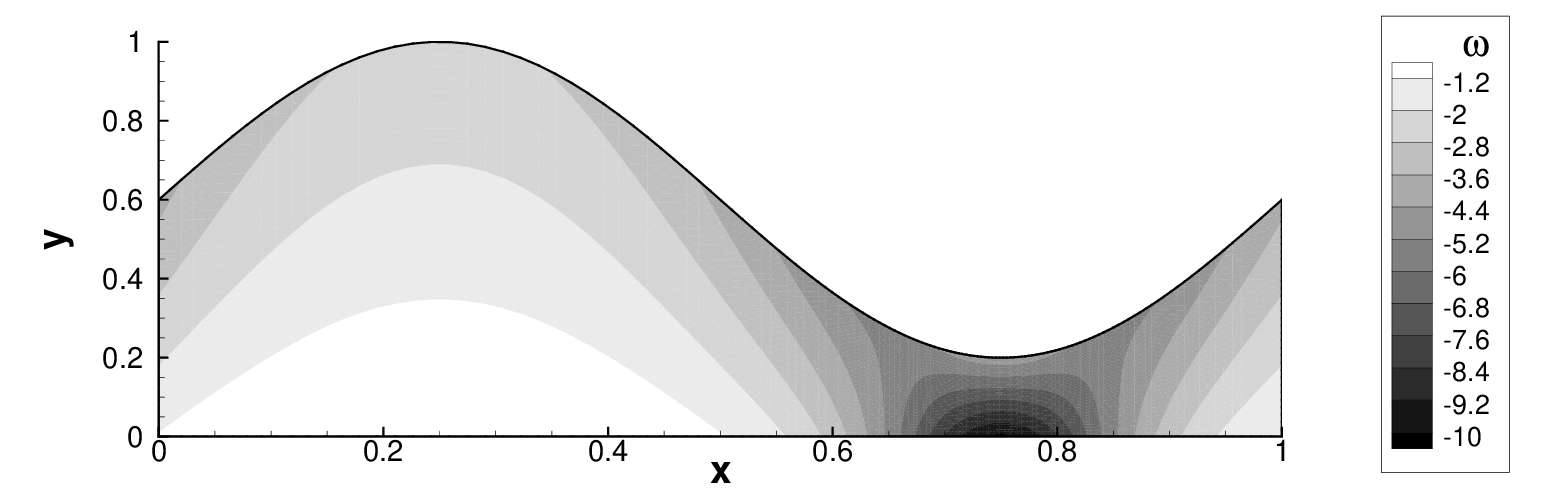}\vspace{8pt}
\includegraphics[width=27pc]
{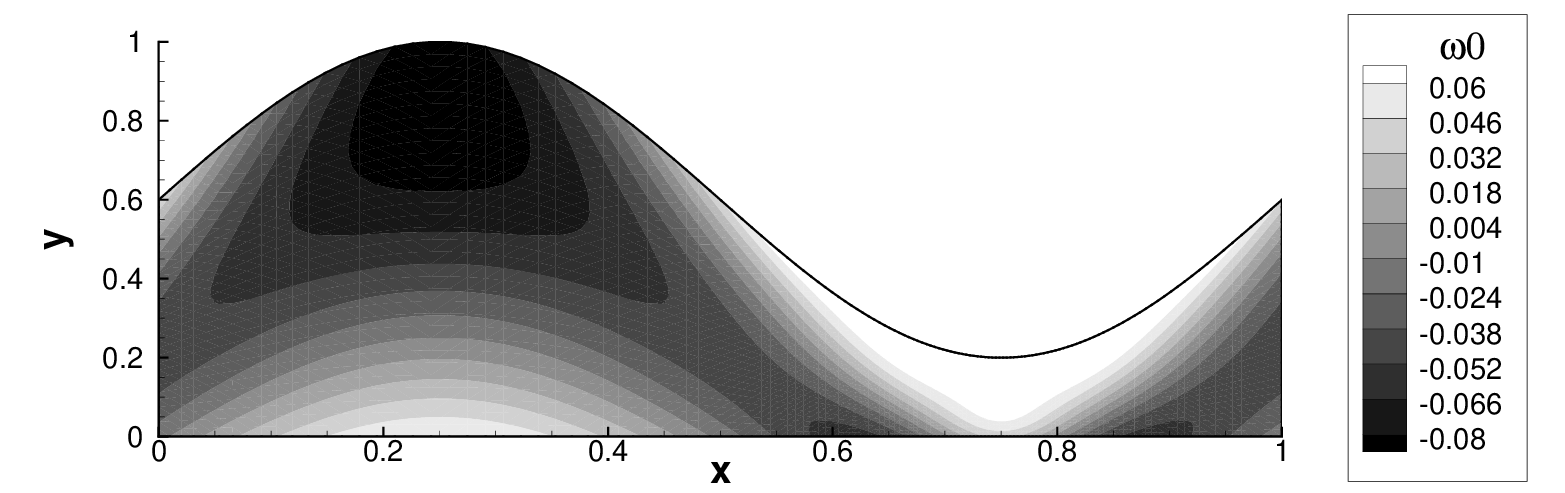}\vspace{8pt}
\includegraphics[width=27pc]
{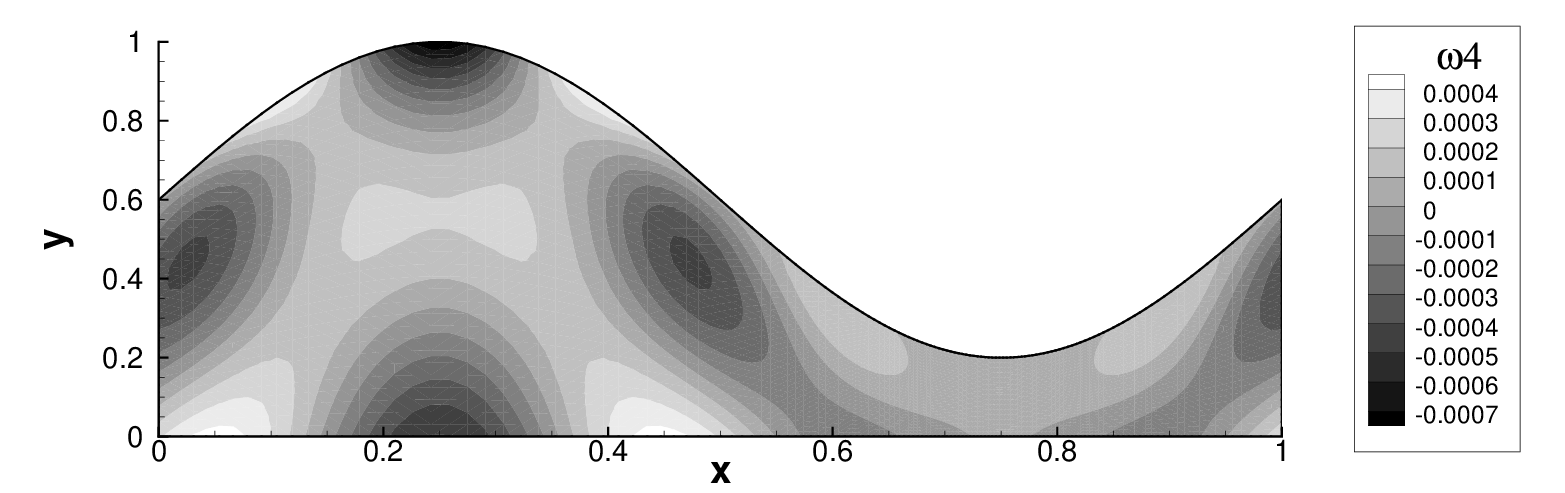}\vspace{8pt}
\includegraphics[width=27pc]
{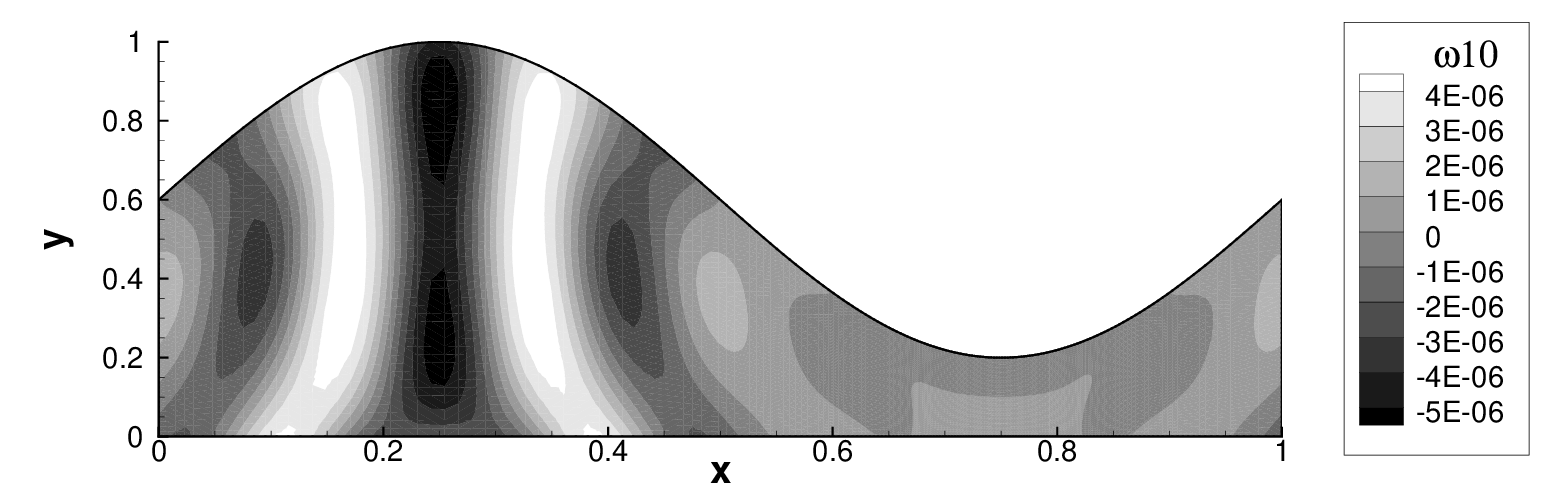}\vspace{8pt}
\includegraphics[width=27pc]
{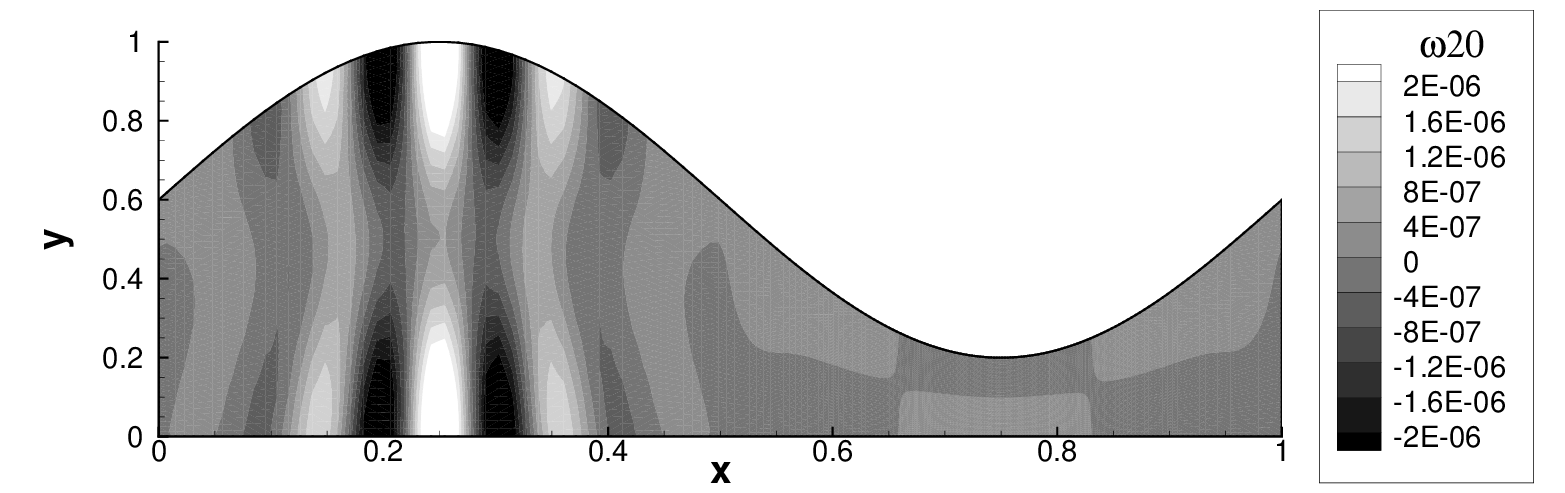}
\mypsfull
\end{center}
\caption{Contour plots of $\omega_\tem{exact}$,
$\omega_\tem{err}^\e{0}$,
$\omega_\tem{err}^\e{4}$,
$\omega_\tem{err}^\e{10}$, and
$\omega_\tem{err}^\e{20}$ for $h(x)$ in {\rm (\ref{eqn:h:sine})} 
with $a=1/5$, $V_0=-0.5$, $V_1=1.0$, and $\veps=0.099$.
Each of these plots corresponds to one of the markers in
Figure~{\rm \ref{fig:errorsAB}}.
}\label{fig:tec:w}
\end{figure}

\begin{figure}[p]
\begin{center}
\mypsdraft
\includegraphics[width=27pc]
{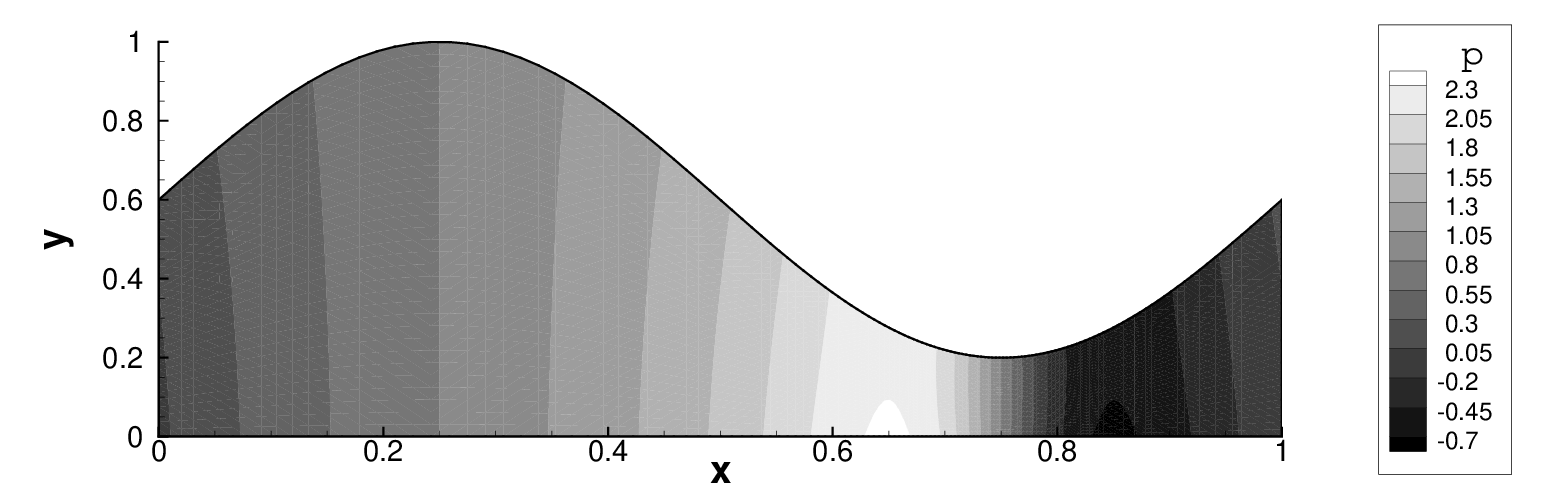}\vspace{8pt}
\includegraphics[width=27pc]
{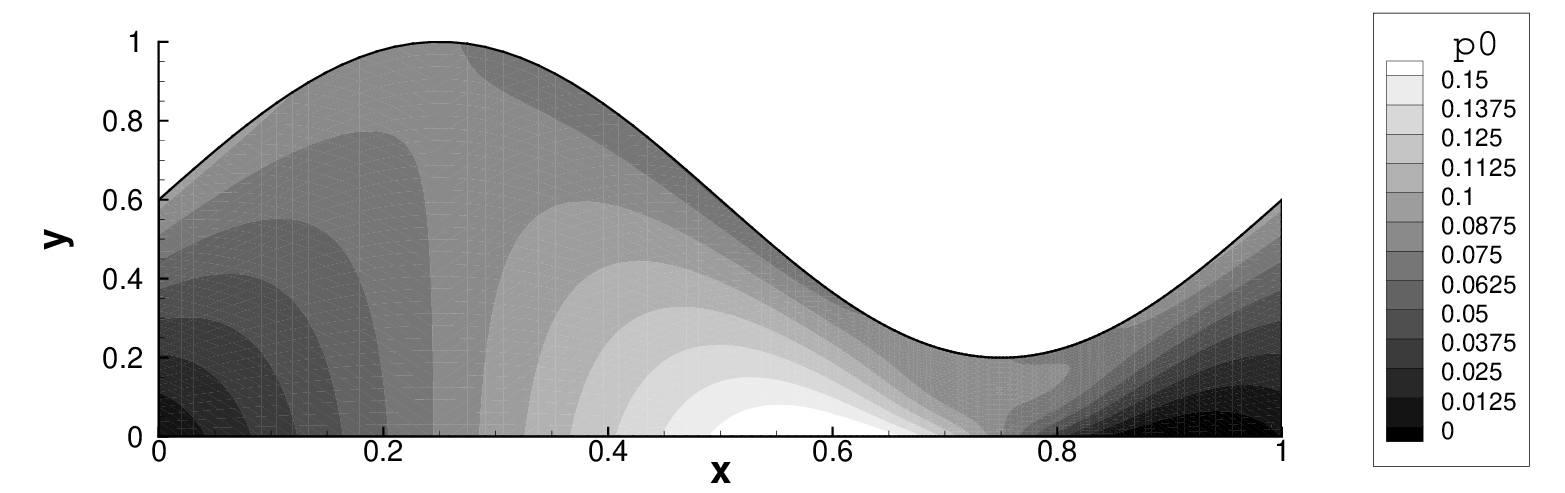}\vspace{8pt}
\includegraphics[width=27pc]
{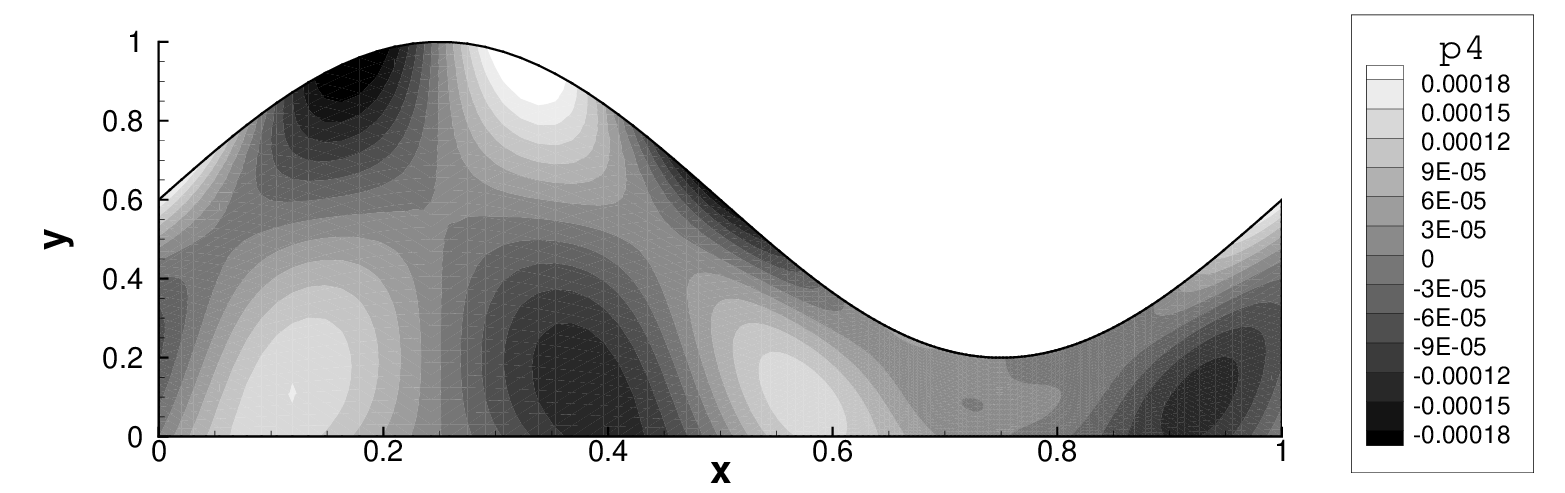}\vspace{8pt}
\includegraphics[width=27pc]
{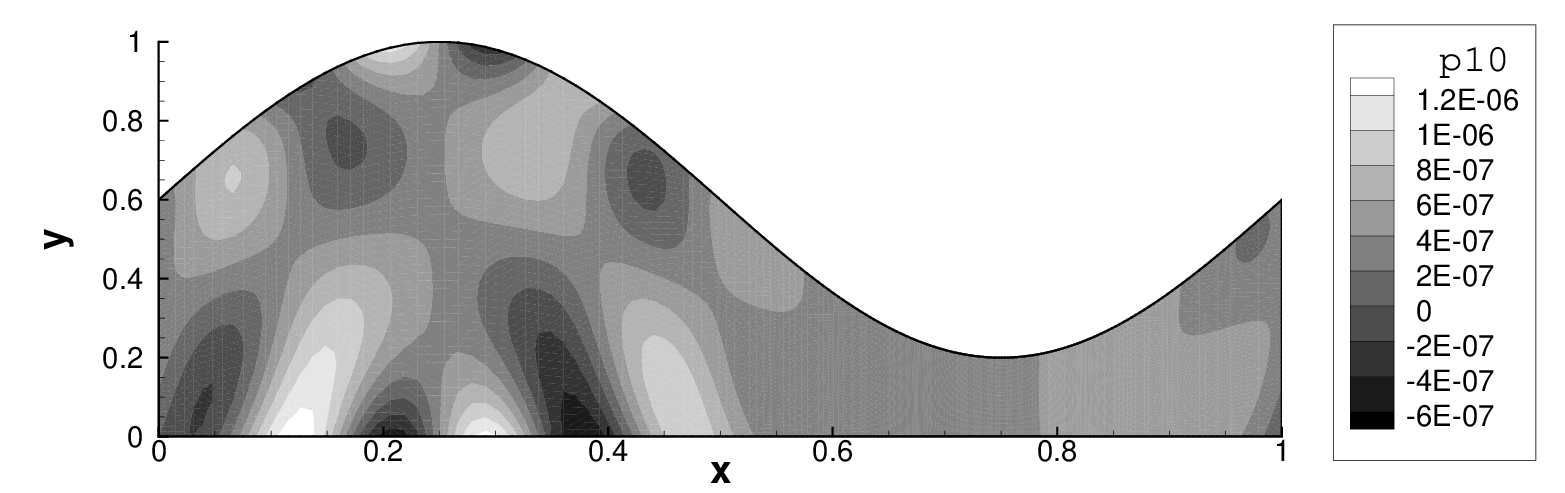}\vspace{8pt}
\includegraphics[width=27pc]
{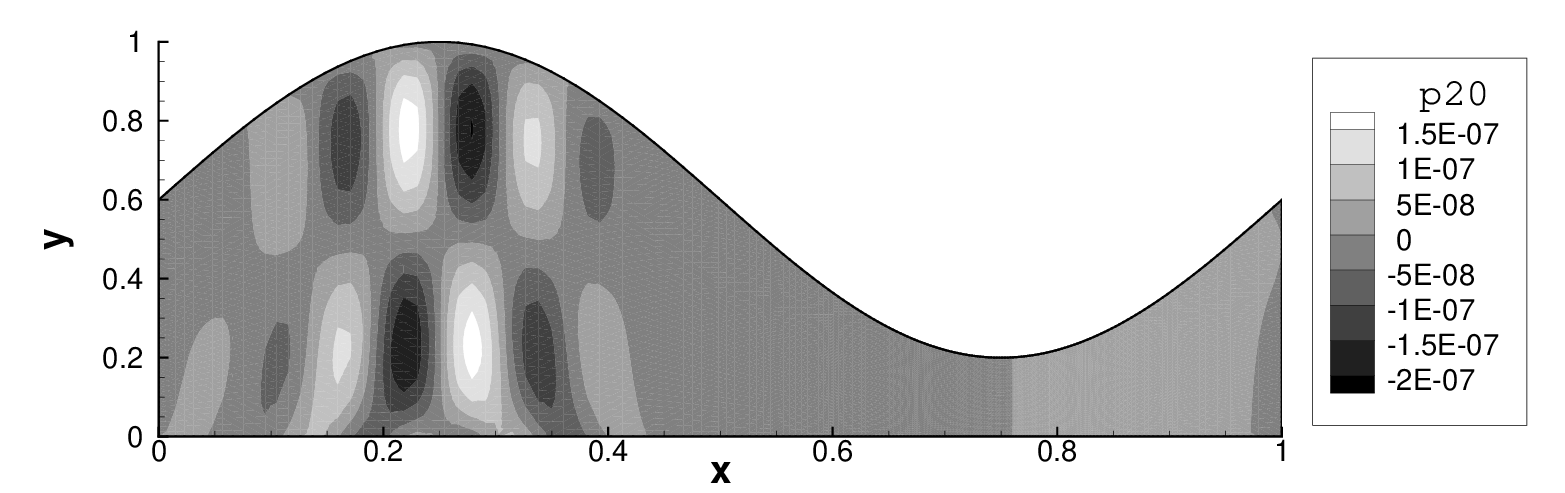}
\mypsfull
\end{center}
\caption{Contour plots of $p_\tem{exact}$,
$p_\tem{err}^\e{0}$,
$p_\tem{err}^\e{4}$,
$p_\tem{err}^\e{10}$, and
$p_\tem{err}^\e{20}$ for $h(x)$ in (\ref{eqn:h:sine})
with $a=1/5$
and $\veps=0.099$.
The ``exact'' solution was computed using a least squares
finite element method with $15$ node quartic triangular
elements on a $2208\times 96$ grid.
}\label{fig:tec:p}
\end{figure}

\begin{figure}[p]
\begin{center}
\includegraphics[height=2.72in]{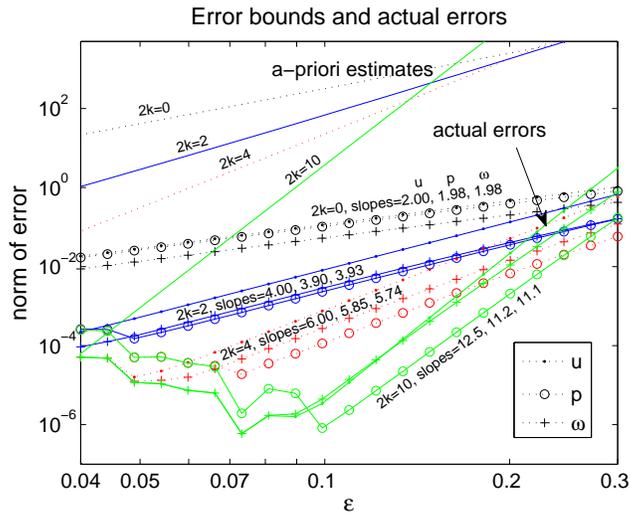}
\end{center}
\caption{Error estimates and actual errors with $a=1/100$.}
\label{fig:error495}
\end{figure}

\begin{figure}[p]
\begin{center}
\mypsdraft
\includegraphics[height=1.44in]{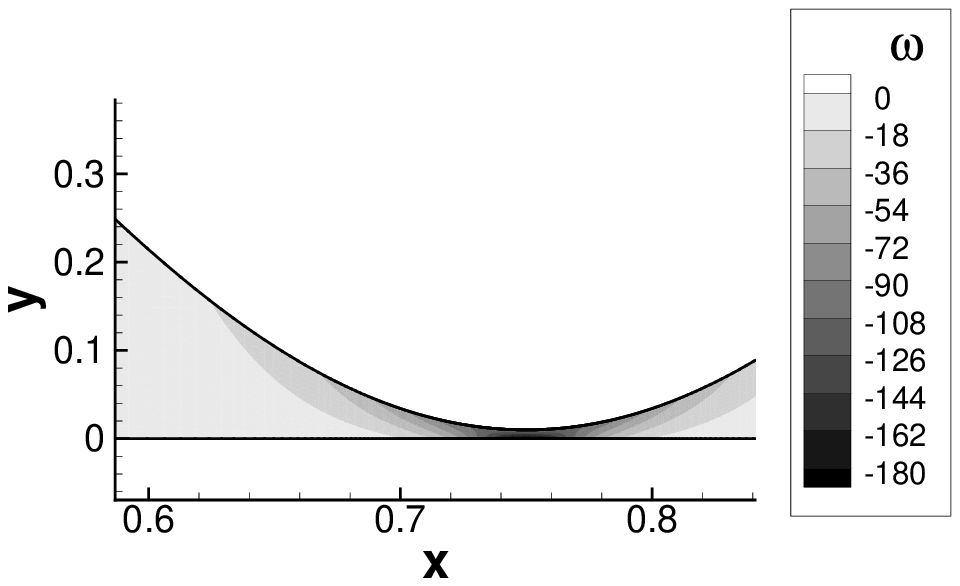}
\includegraphics[height=1.44in]{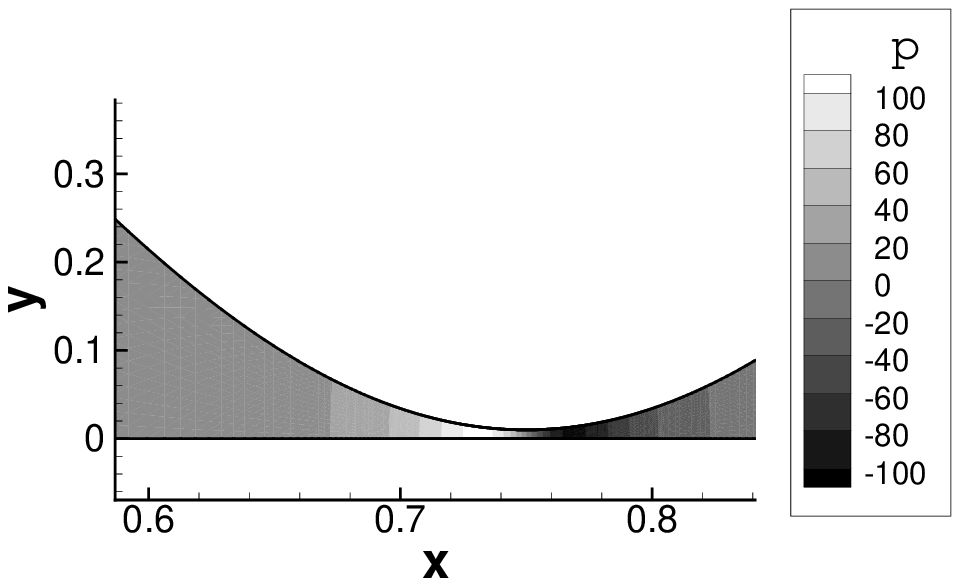}\vspace{6pt}
\includegraphics[height=1.44in]{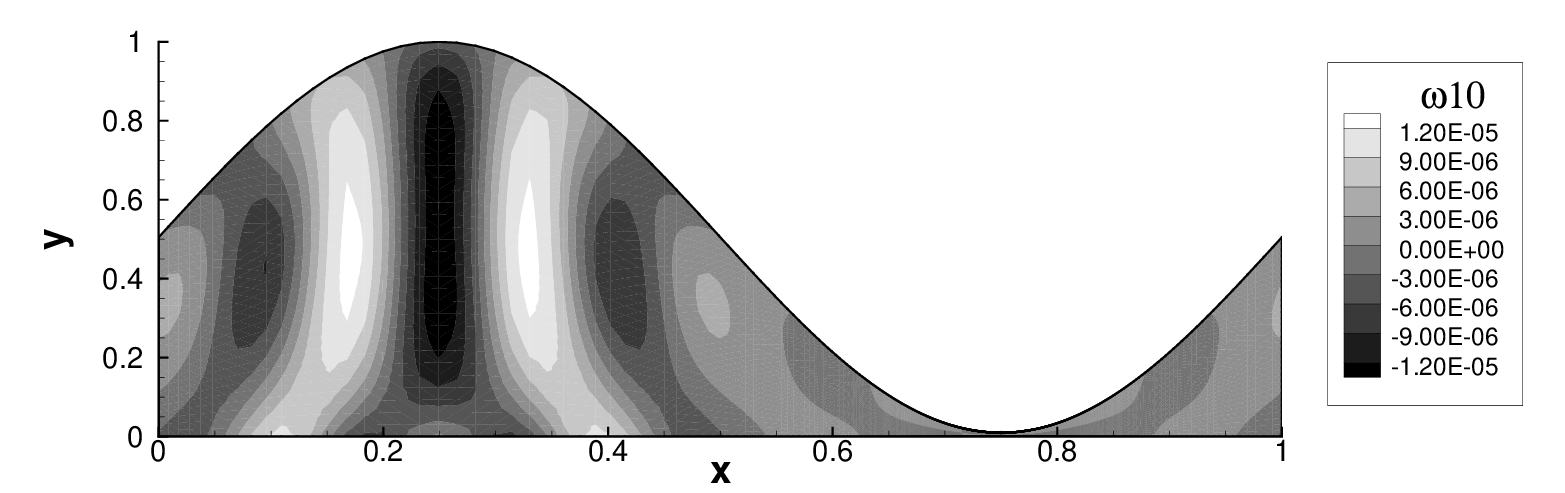}\vspace{6pt}
\includegraphics[height=1.44in]{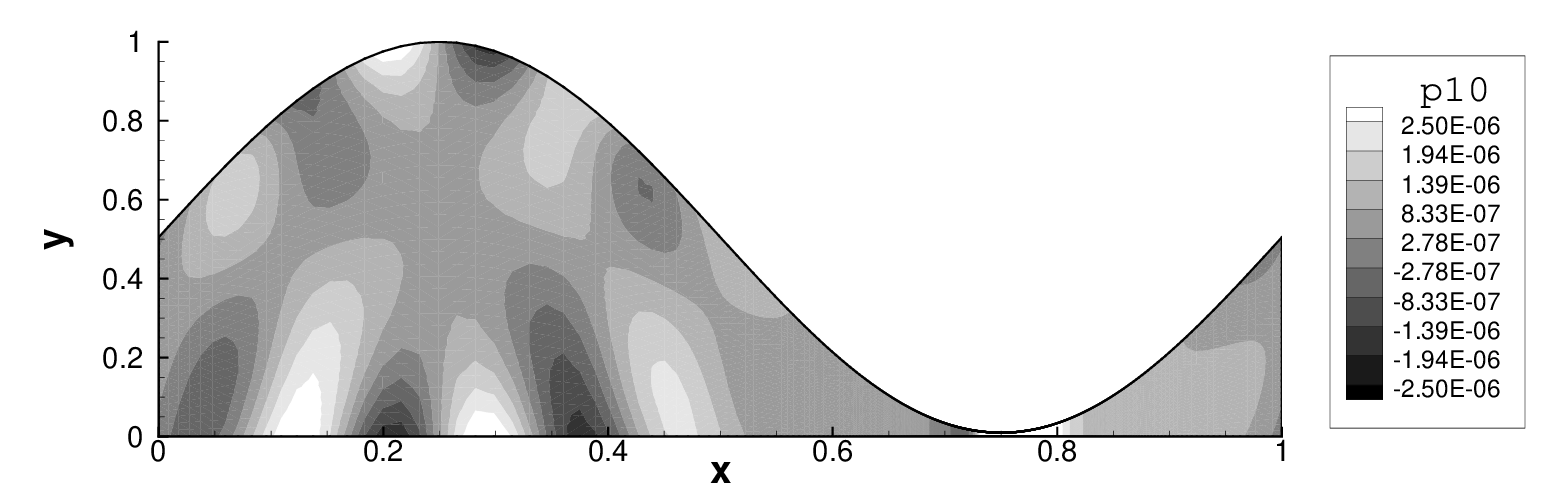}
\mypsfull
\end{center}
\caption{Plots of $\omega_\tem{exact}$, $p_\tem{exact}$,
$\omega_\tem{err}^\e{10}$, and
$p_\tem{err}^\e{10}$ with $a=1/100$, $N=64$, $M=4800$.
}\label{fig:tec2}
\end{figure}

The results are summarized in Figures~\ref{fig:errorsAB}--\ref{fig:tec2}.
For 21 values
of $\veps$ spaced exponentially between $\veps_0=0.04$ and
$\veps_{20}=0.3$, we set up a logically rectangular, $M\times N$
finite element mesh on the domain
\begin{equation}
  \Omega_\veps = \{(x,y) \; : \; 0\le x\le 1, \;\; 0<y<\veps h(x)\}.
\end{equation}
The mesh points are aligned vertically with equal spacing $\Delta
y=h(x)/N$, while the grid spacing in the $x$-direction is chosen to
keep the aspect ratios of the grid cells as close to 1 as possible; we
do this by solving an ODE to enforce $\Delta x\approx h(x)/N$, which
also determines $M$.  For $a=1/5$, we use $N=96$ with $M$ ranging from
768 to 5376 as $\veps$ ranges from $0.3$ to $0.04$; for $a=1/100$, we use
$N=64$ with $M$ ranging from 1600 to 10368.  Four-by-four blocks of
neighboring grid cells are merged and cut into two 15 node triangles.
Interior nodes of the triangles are adjusted to keep the edges straight
except on the top boundary, where we use quartic isoparametric
elements.  We solve the Stokes equations on this mesh using a least
squares finite element method similar to~\cite{fosls:stokes} but
using quartic elements to model the velocity components $u$ and~$v$,
the pressure $p$, the vorticity $\omega=v_x-u_y$, and two strain rates
$\tau=u_y+v_x$ and $\gamma=v_y-u_x$.  We use multigrid to solve the
resulting system of equations, which takes from 3 to 15 minutes
on a 2.4 GHz desktop machine with 16 GB RAM.

Once the finite element solution is known at the grid points, we
normalize the velocity, pressure, and vorticity as described in
section~\ref{sec:derive} and rescale the domain from $\Omega_\veps$
to $\Omega$.  We then use the method described in
Appendix~\ref{sec:impl} to compute $\psi^\e{0}$, $\psi^\e{2}, \ldots,
\psi^\e{20}$ and their derivatives through order 3 at the grid
points.  Next, we use the formulas in (\ref{eqn:uvp:from:psi}) to
obtain $u^\e{2k}$, $v^\e{2k}$, $\omega^\e{2k}$, and $p^\e{2k}$ for
$k=0,\ldots,10$.  For pressure, we use 20 point Gaussian quadrature to
integrate $p_x^\e{2k}$ along the $x$-axis to determine $p^\e{2k}(x,0)$
at the mesh nodes.  The integration of $p_y^\e{2k}$ in the
$y$-direction is done analytically.  With the expansion coefficients
in hand, we evaluate
\begin{align}
  u_\tem{err}^\e{2k} &= u_\tem{exact} - u_\tem{approx}^\e{2k},
  &\quad u_\tem{approx}^\e{2k} &=
  u^\e{0}+\veps^2 u^\e{2} + \cdots + \veps^{2k}u^\e{2k},
\end{align}
etc.,~at the grid nodes, where we use the
finite element solution for $u_\tem{exact}$.  We then run through
the triangles and sum up the local contributions to the errors
\begin{equation}
  \left\|u_\tem{err}^\e{2k}\right\|_{1,\veps}^2 +
  \left\|\veps v_\tem{err}^\e{2k}\right\|_{1,\veps}^2, \qquad
  \left\|p_\tem{err}^\e{2k}\right\|_0^2, \qquad
  \left\|\omega_\tem{err}^\e{2k}\right\|_0^2
\end{equation}
by interpolating the values at the grid nodes and integrating the
resulting polynomials on the triangle; this step is very similar to
the assembly of the stiffness matrix.  Finally, we store the results
in a file for visualization
(see Figures~\ref{fig:tec:w}, \ref{fig:tec:p}, and~\ref{fig:tec2})
and record the norms of the truncation errors
for comparison with the error bounds of Theorem~\ref{thm:bound2}.

The results of this comparison are shown in Figures~\ref{fig:errorsAB}
and~\ref{fig:error495}.  As expected, for fixed $k$, the actual errors
decay as $O(\veps^{2k+2})$.  The a priori error bounds eventually 
decrease like $O(\veps^{2k+2})$ as well, but the term involving
$\theta_k$ in (\ref{eqn:star:def2}) is significant over this range
of $\veps$ in some of the cases, causing the slopes to be larger:
\begin{equation*}
  \frac{\theta_k}{r_k}\sqrt{\frac{I_3}{I_1}} =
  \left\{\begin{array}{c|c|c|c|c|c}
    & k=0 & k=1 & k=2 & k=5 & k=10 \\
    \hline
    a=1/5 & 12.5 & 0.94 & 0.16 & .00096 & 1.3\times10^{-11^{\phantom{1}}} \\
    a=1/100 & 257 & 19.3 & 3.2 & 0.020 & 2.7\times 10^{-10^{\phantom{1}}}
  \end{array}.\right.
\end{equation*}
This effect is much more pronounced when $a=1/100$ in (\ref{eqn:h:sine})
due to
\begin{equation}
  \sqrt{\frac{I_3}{I_1}} = \frac{1}{2}\sqrt{\frac{3}{2} +
    \frac{1}{a} + \frac{3}{2a^2}} = \begin{cases}
    3.32, & a=1/5, \\ 61.4, & a=1/100.
  \end{cases}
\end{equation}
The deviation from linearity in the plots of ``actual error'' for
small $\veps$ and large $k$ is due to error in the finite element
solutions, which are accurate to about 9 digits.  This occurs sooner
when $a=1/100$ since the pressure and vorticity of the exact solution
in the vicinity of the narrow gap increases as $a$ decreases, and also
because we were forced to use a coarser mesh with $a=1/100$ to avoid
running out of computer memory in the finite element simulations.  The
data points with $\veps=0.099$ in Figure~\ref{fig:errorsAB} correspond
to the contour plots in Figures~\ref{fig:tec:w} and~\ref{fig:tec:p},
where we plot $\omega_\tem{exact}$, $\omega_\tem{err}^\e{2k}$,
$p_\tem{exact}$, and $p_\tem{err}^\e{2k}$ for $2k=0,4,10,20$.  The data
points with $\veps=0.099$ in Figure~\ref{fig:error495} correspond to
the contour plots in Figure~\ref{fig:tec2}.  We remark that the
apparently large value of $p_\tem{err}^\e{10}$ in the narrow gap in
Figure~\ref{fig:tec2} is due to smoothing in the least squares finite
element solver; the expansion solution is more accurate than the finite
element solution in this region of the domain.  The error patterns
that emerge in all these cases are rather interesting, indicating that
the spaces $\mc{H}_{2k}$ in Theorem~\ref{thm:struc} (the structure
theorem) can be quite complicated even for simple curves $h(x)$.

Although our estimates for the error in pressure include an additional
factor of $h_0^{-3/2}(r_k+r_k^{-1})^2$, all our numerical experiments
(including complicated geometries 
in which the inf-sup constant $\beta^{-1}$ does exhibit
$h_0^{-3/2}$ behavior) indicate that $\|p_\tem{err}^\e{2k}\|_0$ is
comparable to $\|\omega_\tem{err}^\e{2k}\|_0$.  In fact, for large
$k$, pressure seems to be the most accurately computed variable; see
Figures~\ref{fig:errorsAB} and~\ref{fig:error495}.
We do not know how to explain this as the
pressure \emph{is} determined by solving (\ref{eqn:p:err:from:psi}),
which involves inverting the operator
$\nabla:L^2_\#(\Omega)\rightarrow H^{-1}(\Omega)^2$.  For some reason,
in lubrication-type problems, the right-hand side $\mb{f}_k$
belongs to a subspace of $H^{-1}(\Omega)^2$ that is not amplified by
$(\nabla)^{-1}$ when solving $\nabla p_\tem{err}^\e{2k}=\mb{f}_k$.

The following table shows the minimum ratio of the a priori error
estimate to the actual error $\|\mb{u}_\tem{err}^\e{2k}\|_{1,\veps}$
for the data points in Figures~\ref{fig:errorsAB} and~\ref{fig:error495}
that were used to compute the slopes of the best-fit lines:
\begin{equation*}
  \begin{array}{r|c|c|c|c|c}
   \multicolumn{1}{r}{k=} & 0 & 1 & 2 & 5 & 10 \\
    \hline
    \text{(min ratio, $a=1/5$)}^{1/(2k+2)}\rule{0pt}{12pt}   & 11.1 & 7.0 &
    5.0 & 2.8 & 2.3 \\
    \text{(min ratio, $a=1/100$)}^{1/(2k+2)}   & 34.0 & 8.6 &
    5.5 & 3.0 & -  \end{array}.
\end{equation*}
For example, in the 10 calculations (with $\veps$ ranging from
$0.04\le\veps\le0.099$) that were used to determine the slope of the
$2k=4$ line in Figure~\ref{fig:errorsAB}, the ratios of the a priori
errors to the exact errors ranged between $1.608\times10^4$ and
$1.617\times10^4$, so we recorded
$\sqrt[6]{1.608\times10^4}\approx5.0$.
This table gives information on how far the values~$\rho_k$ in
Table~\ref{tbl:rho:theta} are from their optimal values.
For example, if we increased $\rho_5$ by more than a factor of $2.8$
while holding $\theta_5$ fixed, the estimate (\ref{eqn:thm:bound2})
would fail to hold for this geometry.  Since $r_k^{-1}$ in
(\ref{eqn:r:def}) is used as a convenient upper bound on all the
integrals $|E_{m,j}^\e{2\ell}|^{1/2\ell}$ and
$|\wtil{E}_{m,j}^\e{2\ell}|^{1/4\ell}$ that arise in the
definition of $Q^\e{2k}$ and also in the bounds for
$\|\psi^\e{2k}_{xx}\|_0$ and
$\|h^2\psi^\e{2k-2}_{xxxx}\|_0$, it is remarkable
that the values of $\rho_k$ we computed are within a factor of 3 of
optimal for $k=5$, $k=10$, and perhaps all $k\ge5$.

\section{Discussion}
\label{sec:discuss}

Although we are able to estimate the effective radius of convergence
$\rho_kr_k$ quite closely, our estimates of
$\|\psi_\tem{err}^\e{2k}\|_{2,\veps}$,
$\|\omega_\tem{err}^\e{2k}\|_0$, etc.,~are likely to be several orders
of magnitude too large.  One shouldn't expect an a priori bound that
holds for all geometries alike to provide an exceptionally sharp bound
for any specific geometry.  Instead, our analysis provides a clear
picture of the features of $h(x)$ that cause the effective radii of
curvature $r_k\rho_k$ to become small, namely, large values of
$h^{k-1}\partial_x^kh$.  No previous study has ever described how the
constant hidden in the $O(\veps^{2k+2})$ depends on $h$; instead, $h$
has always been fixed at the outset and only the limit as
$\veps\rightarrow0$ has been considered.

\begin{figure}[t]
\begin{center}
\includegraphics[width=23pc]
{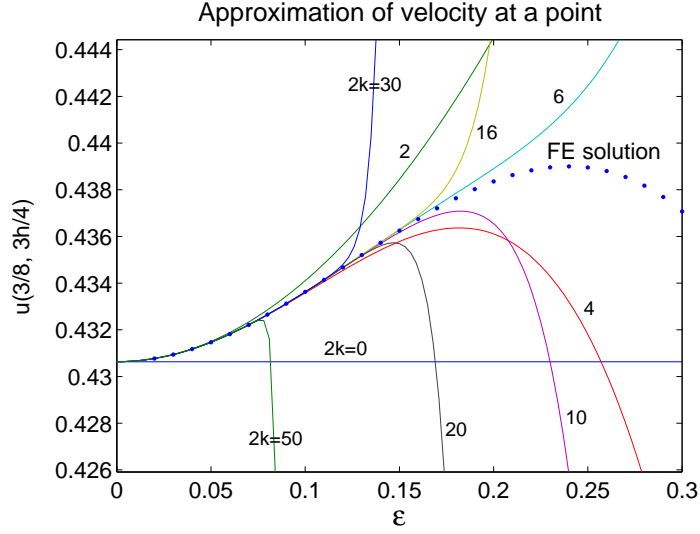}
\end{center}
\caption{Comparison of $u_\tem{approx}^\e{2k}$ (solid lines) to
$u_\tem{exact}$ (dots) at the point
$(x,y)=(\frac{3}{8},\frac{3}{4}h)$ for $2k=0,2,4,6,10,16,
20,30,50$.  Here $h(x)=\frac{3}{5}+\frac{2}{5}\sin(2\pi x)$,
$V_0=-0.5$, and $V_1=1$.  This function $h(x)$ is real analytic and
periodic, yet the expansion solution appears to be an asymptotic
series rather than a convergent series.  }\label{fig:u_at_p}
\vspace{-10pt}
\end{figure}

Another feature of this analysis is that it separates the constants
into two types: those that are (1) given in the problem statement or
easily computable from $h$; or (2) difficult to compute but universal
(independent of $h$).  We listed the first several constants in the
latter category ($\rho_k$ and $\theta_k$) in
Table~\ref{tbl:rho:theta}.  It is interesting that $\rho_k$ actually
increases until $2k=10$ and doesn't get as bad as $\rho_0$ again
until $2k=26$.  However, at that point it seems to be decreasing
steadily like $1/k$, indicating that the effective radius of curvature
in our a priori error bound will shrink to zero as
$k\rightarrow\infty$.  The reason for this is that the recurrences
(\ref{eqn:AiB:recur}) and (\ref{eqn:4th:derivs}) relating the matrices
$A_i^\e{2k}$ and $B^\e{2k}$ to their lower order counterparts cause
the norms of these matrices to grow like $k!$. Thus, although $\rho_k$
involves $k$th roots of these constants, these $k$th roots still grow
linearly in $k$.  On the other hand, if $h$ is real analytic as well
as periodic, a standard contour integral argument shows that there is
an $r>0$ such that $\|\partial_x^k h\|_\infty\le k!\,r^{-k}$ for
all $k\ge0$; thus, the constants $r_k$ will remain bounded away from
zero.  For example, if $h(x)$ is of the form (\ref{eqn:h:sine}),
one may show that if $a\in(0,2/3]$, then the largest
value of $\|\frac{1}{\ell!}h^{\ell-1}\partial_x^\ell
h\|_\infty^{1/\ell}$ occurs when $\ell=2$, so all the $r_k$ are
equal to $r_0=(\pi\sqrt{1-a})^{-1}$.  It is conceivable that when
$h$ is real analytic, the norms of the functions $\psi^\e{2k}$ grow
slowly enough that the stream function expansion converges in spite of
the fact that the matrices $A_i^\e{2k}$ and $B^\e{2k}$ in their
representation (\ref{eqn:mat:rep}) blow up like $k!$. This would
simply mean that we chose a bad basis in terms of which to represent
$\psi$. We used orthogonal polynomials in (\ref{eqn:A:dot2}) to
improve this basis, but there may be other improvements.
Figure~\ref{fig:u_at_p} shows that this is not the case.
Even when $h$ varies sinusoidally, the expansion solution appears to
be an asymptotic series rather than a convergent series: all the
variables, including the flux terms $Q^\e{2k}$, appear to grow like
$k!$ as $k$ becomes large.

Nevertheless, the expansion solutions can be extremely accurate
(almost exact) as long as they are used for a geometry that falls
within the effective radius of convergence of the truncated series.
It is hoped that the estimates in this paper will help to identify
these cases and provide practical a priori (as well as a posteriori)
error estimates for many interesting problems.

\appendix

\section{Implementation}
\label{sec:impl}
We have developed two methods for computing the higher order corrections
described in section~\ref{sec:struc} using a computer.
In the first, we use Mathematica to evaluate the derivatives and
antiderivatives in recursion~(\ref{eqn:psi:Q:recur}) and
Algorithm~\ref{alg:G:def} symbolically.  With this approach, the main
challenge occurs at the step where $Q^\e{2k_0}$ is defined as a
definite integral.  We do this through pattern matching and symbol
replacement.  At the stage where the definite integral is to be
evaluated, we replace all instances of $\partial_x^j h$ in the
integrand by $j!\,t_j/h^{j-1}$.  Each term in the result (call it $R$)
will contain a factor of $h^{-3}$ or $h^{-2}$, with no other dependence\vspace*{-1pt}
on $h$.  For each $k=k_0,\ldots,0$ and $j=1,\ldots,d_{2k}$, we find the
terms in $R$ that contain $\varphi^\e{2k}_j$ (left in the form
$t_1^{i_1}\cdots t_{2k}^{i_{2k}}$ described in
Algorithm~\ref{alg:Phi}) as a factor.  These terms are removed from
$R$ while their symbolic integrals (with $\varphi^\e{2k}_j/h^m$
replaced by $I_m E^\e{2k}_{m,j}$) are divided  by $2I_3$ and added to
the desired flux $Q^\e{2k_0}$.  By running  through the
$\varphi^\e{2k}_j$ in decreasing $k$ order, we convert higher order
products (e.g., $t_1^2 t_2/h^3$) into symbols (e.g., $I_3E^\e{4}_{3,2}$)
before one of their lower order factors can be converted incorrectly
(e.g., into $I_3 E^\e{2}_{3,1} t_2$).  This approach is effective
through 6th or 8th order but becomes rather slow as the complexity of
the expansion increases.

The second approach is much faster and can be implemented in any
modern programming language.  We have written a version in $C^{++}$
and a version in Mathematica.  Instead of representing the basis
functions $\varphi^\e{k}_j$ for $\mc{H}_k$ using a computer algebra
system, we represent them as $(k+1)$-tuples of integers.  For example,
the functions $1$, $h_x$, $\frac{1}{6}h^2 h_x^3 h_{xxx}$, and
$\frac{1}{48}h^4 h_x h_{xx} h_{xxxx}$ in $\mc{H}_0$, $\mc{H}_1$,
$\mc{H}_6$, and $\mc{H}_{7}$ are represented by $(0)$, $(0,1)$,
$(2,3,0,1,0,0,0)$, and $(4,1,1,0,1,0,0,0)$.  A tuple $(i_0,\ldots,i_k)$
represents a basis function for $\mc{H}_k$ iff
\begin{equation}
  \label{eqn:Hk:condition}
  i_1+2i_2+\cdots+ki_k=k, \qquad i_0 = i_2 + 2i_3 + \cdots + (k-1)i_k.
\end{equation}
We begin by constructing the basis sets $\Phi_k$ for $0\le k\le2k_0$
and storing them as $(k+1)\times d_k$ integer matrices with columns
corresponding to the $\varphi^\e{k}_j$.  This is done using
Algorithm~\ref{alg:Phi}, which returns the columns sorted
lexicographically from the last slot to the first slot
(e.g., $(3,0,3,0)^T<(2,3,0,1)^T<(3,1,1,1)^T$).
Sorted columns allow us to find the column index corresponding to
a given tuple in $\log_2 d_k$ time.

Next, for $0\le k\le2k_0-1$, we compute the operators $h\partial_x$
and $h_x\cdot$ from $\mc{H}_k$ to $\mc{H}_{k+1}$ and store them as
sparse integer matrices of dimension $d_{k+1}\times d_k$.  If column
$J$ of $\Phi_k$ contains the tuple $(i_0,\ldots,i_k)$, we define
$i_{k+1}=0$ and compute
\begin{align}
  h_x\cdot : (i_0,\ldots,i_k) &\mapsto (i_0,i_1+1,i_2,\ldots,i_{k+1}), \\
  \notag
  h\partial_x : (i_0,\ldots,i_k) &\mapsto
\sum_{\{r\,:\,i_r\ne0\}} i_r(r+1)(i_0+1,\ldots,i_r-1,i_{r+1}+1,\ldots,i_{k+1}),
\end{align}
where the omitted indices are unmodified and the $+1$ and $-1$ cancel
in the first slot when $r=0$ in the sum.  The factor of $(r+1)$ is due
to the factorials in the definition of the $\varphi^\e{k}_j$.  The
column index $l$ of each $(k+2)$-tuple in the result is found in
$\Phi_{k+1}$, and the corresponding coefficient (1 or $i_r(r+1)$) is
added to the $l$th row and $J$th column of the sparse matrix
representing $h\partial_x$ or $h_x\cdot$.  The entries of these sparse
matrices are positive, and the column sums (i.e., 1-norms) are all equal
to $1$ for $h_x\cdot$ and to $i_0+2i_1+\cdots+(k+1)i_k=2k$ for
$h\partial_x$ (by (\ref{eqn:Hk:condition})).

Once the operators $h\partial_x$ and $h_x\cdot$ are known, we use them
to recursively compute the matrices $A^\e{2k}=V_0 A_0^\e{2k} + V_1
A_1^\e{2k}$ and $B^\e{2k}$ in (\ref{eqn:mat:rep}).  We start by
setting $A_0^\e0=(0,1,-2,1)^T$, $A_1^\e0=(0,0,-1,1)^T$, and
$B^\e0=(0,0,3,-2)^T$ as in Example~\ref{exa:thm}.  For $1\le k\le k_0$,
we mimic the proof of Theorem~\ref{thm:struc} to build up $A^\e{2k}$
and $B^\e{2k}$ row by row.  For $4\le n\le 2k+3$ and $i=0,1$, we use
sparse matrix--vector multiplication to define the rows
\begin{equation}
\label{eqn:AiB:recur}
  \begin{aligned}
  A_i^\e{2k}(n,:) &= \left(\frac{
    -2[h\partial_x - (n-2)h_x][h\partial_x - (n-3)h_x]
    \left[A_i^\e{2k-2}(n,:)^T\right]}{n(n-1)}\right)^T,\\
  B^\e{2k}(n,:) &= \left(\frac{
    -2[h\partial_x - (n-1)h_x][h\partial_x - (n-2)h_x]
    \left[B^\e{2k-2}(n,:)^T\right]}{n(n-1)}\right)^T.
  \end{aligned}
\end{equation}
If $k\ge2$, then for $6\le n\le 2k+3$, we add the following vectors to
$A_i^\e{2k}(n,:)$ and $B^\e{2k}(n,:)$, respectively:\vspace{5pt}
\begin{equation}
  \label{eqn:4th:derivs}
 \begin{aligned}
    &\textstyle\hspace*{-6pt}\left(\hspace*{-2pt}\frac{
    -[h\partial_x - (n-2)h_x][h\partial_x - (n-3)h_x]
    [h\partial_x - (n-4)h_x][h\partial_x - (n-5)h_x]
    \left[A_i^\e{2k-4}(n,:)^T\right]}{n(n-1)(n-2)(n-3)}\hspace*{-1pt}\right)^T\hspace*{-3pt},\\
    &\textstyle\hspace*{-6pt}\left(\hspace*{-2pt}\frac{
    -[h\partial_x - (n-1)h_x][h\partial_x - (n-2)h_x]
    [h\partial_x - (n-3)h_x][h\partial_x - (n-4)h_x]
    \left[B^\e{2k-4}(n,:)^T\right]}{n(n-1)(n-2)(n-3)}\hspace*{-1pt}\right)^T\hspace*{-3pt}.
  \end{aligned}
\end{equation}
Next we zero out rows 0 and 1 of $A_0^\e{2k}$, $A_1^\e{2k}$, $B^\e{2k}$,
and set\vspace{5pt}
\begin{equation}
  \begin{aligned}
  A_i^\e{2k}(2,:) &= \sum_{n=4}^{2k+3} (n-3) A_i^\e{2k}(n,:), \quad
  A_i^\e{2k}(3,:) = \sum_{n=4}^{2k+3} (2-n) A_i^\e{2k}(n,:),\\
  B^\e{2k}(2,:) &= \sum_{n=4}^{2k+3} (n-3) B^\e{2k}(n,:), \quad
  B^\e{2k}(3,:) = \sum_{n=4}^{2k+3} (2-n) B^\e{2k}(n,:).
  \end{aligned}
\end{equation}
Finally, we subtract ${-1/2 \choose k}$ from $A_1^\e{2k}(2,1)$ and add
it to $A_1^\e{2k}(3,1)$ to account for the boundary data, where we
recall that the rows and columns are indexed starting at 0 and 1,
respectively.  Using this
approach, our $C^{++}$ code can compute these matrices through order
$2k=50$ using floating point arithmetic in a few seconds, while our
Mathematica code can compute through order $2k=30$ in exact rational
arithmetic in about an hour.  This allows us to explore the properties
of the stream function expansion and test our error estimates to quite
a high order.



\end{document}